\documentclass[12pt]{amsart}
\usepackage{amssymb}
\usepackage{times}
\usepackage{graphicx}

\theoremstyle{plain}

\numberwithin{equation}{section}

\newcommand{\calA}{\mathcal{A}}
\newcommand{\calB}{\mathcal{B}}

\newcommand{\calD}{\mathcal{D}}
\newcommand{\calE}{\mathcal{E}}

\newcommand{\calH}{\mathcal{H}}
\newcommand{\calK}{\mathcal{K}}
\newcommand{\calL}{\mathcal{L}}
\newcommand{\calM}{\mathcal{M}}
\newcommand{\calN}{\mathcal{N}}
\newcommand{\calO}{\mathcal{O}}

\newcommand{\calR}{\mathcal{R}}
\newcommand{\calS}{\mathcal{S}}

\newcommand{\calI}{\mathcal{I}}

\newcommand{\calW}{\mathcal{W}}

\newcommand{\bbF}{\mathbb{F}}
\newcommand{\bbC}{\mathbb{C}}

\newcommand{\bbP}{\mathbb{P}}
\newcommand{\bbQ}{\mathbb{Q}}
\newcommand{\bbR}{\mathbb{R}}

\newcommand{\bbZ}{\mathbb{Z}}

\newcommand{\bfe}{\mathbf{e}}

\newcommand{\bft}{\mathbf{t}}

\newcommand{\bfF}{\mathbf{F}}
\newcommand{\la}{\langle}
\newcommand{\ra}{\rangle}
\newcommand{\li}{^\textup{l}}

\def\SL{{\text{SL}}}
\def\Aut{{\text{Aut}}}
\def\cub{{{\text{cub}}}}
\def\Pic{{\text{Pic}}}
\def\Ker{{\text{Ker}}}
\def\O{{\text{O}}}
\def\Gal{{\text{Gal}}}
\def\ss{{\text{ss}}}
\def\st{{\text{s}}}
\def\ncub{{{\text{ncub}}}}

\begin{document}
\title [A Complex Ball Uniformization]
{A Complex Ball Uniformization of
the Moduli Space of Cubic Surfaces Via Periods of K3 Surfaces}
\author{I.\ Dolgachev}
\address{Department of Mathematics, University of Michigan, Ann Arbor, MI
48109,USA}
\email{idolga@umich.edu}
\thanks{Research of the first author is partially supported by NSF grant
DMS 9970460}

\author{B.\ van Geemen}
\address{Dipartimento di Matematica, Universit\`a di Milano, Via Saldini
50, I-20133 Milano, Italy}
\email{geemen@mat.unimi.it}

\author{S.\ Kond{$\bar{\rm o}$}}
\address{Graduate School of Mathematics, Nagoya University, Nagoya,
464-8602, Japan}
\email{kondo@math.nagoya-u.ac.jp}
\thanks{Research of the third author is partially supported by
Grant-in-Aid for Scientific Research A-14204001, Japan}

\begin{abstract}
In this paper we show that the moduli space  of nodal cubic surfaces
is isomorphic to a quotient of a
4-dimensional complex ball by an arithmetic subgroup of the
unitary group. This complex ball
uniformization uses the periods of certain $K3$ surfaces which are
naturally associated to cubic
surfaces. A similar uniformization is given for the covers of
the moduli space corresponding to geometric markings of
the Picard group or to the choice of a line on the surface.
We also give a detailed description of the
boundary components corresponding to singular surfaces.
\end{abstract}
\maketitle

CONTENTS
\begin{itemize}
\item[1.] Introduction
\item[2.] Nodal cubic surfaces
\item[3.] Cubic surfaces and 2+5 points on the line
\item[4.] The $K3$ surface associated to a cubic surface
\item[5.] The Picard lattice
\item[6.] The moduli space of $K3$ surfaces associated to a cubic surface
\item[7.] A complex ball uniformization
\item[8.] The geometry of the discriminant locus
\item[9.] Extension of the isomorphism to the boundary
\item[10.] Half twists
\end{itemize}

\section{Introduction}
There are two main approaches to the construction of moduli spaces in
algebraic geometry.
One uses geometric invariant theory which allows one to construct the
moduli space as a
quotient of an open subset of an appropriate Hilbert scheme, the other one
uses period maps to construct
the moduli space as a quotient of an open subset of a Hermitian symmetric
homogeneous domain by a discrete subgroup
of its group of holomorphic automorphisms. Both approaches suggest a way to
compactify the moduli space.
In the algebraic approach one adds the equivalence classes of semi-stable
points. In the transcendental approach
one considers the whole domain together with its boundary.

There are several remarkable cases where both approaches work. Comparing
the constructions gives a beautiful
interplay between the algebraic theory of invariants and the theory of
automorphic functions. The historically first
example of such an interplay is of course the moduli space of elliptic
curves which, on one hand, is the quotient of the
space of binary forms of degree 4 by the group $\SL(2)$ and, on the other
hand, is a natural quotient of the upper
half-plane by the modular group. Similarly, binary forms of degree 5, 6, 8
and 12 give the moduli spaces of Del Pezzo
surfaces of degree 4, and hyperelliptic curves of genus 2, 3 and 5,
respectively. Using the theory of hypergeometric
functions one can show that the corresponding domains are complex balls of
dimension 2, 3, 5 and 9, respectively.
Increasing the number of variables by one, one finds the ternary cubic
forms which leads again to the moduli space of
elliptic curves, the forms of degree 4 corresponding to the moduli space of
non-hyperelliptic curves of genus 3 (in this
case the domain is the Siegel upper half space of degree 3) and the forms
of degree 6 corresponding to $K3$ surfaces
with degree 2 polarization (the domain is of type IV in Cartan's
classification).

Using domains of type IV one can also give a uniformization of the moduli
space of cubic and quartic forms in
4 variables. The case of forms of degree 3 (cubic surfaces) was treated in
the work of K. Matsumoto, T. Sasaki and
M. Yoshida \cite{MSY}, and degree 4 ($K3$ surfaces with degree 4
polarization) much earlier by J.\ Shah \cite{Sha}.
Although cubic surfaces do not admit non-zero holomorphic 2-forms, so that
the periods are not defined, there are
identifications of this moduli space with other moduli spaces for which the
period map is defined.
In \cite{MSY} one uses the moduli space of $K3$ surfaces which have a
certain primitive sublattice of rank 16 in the
Picard group. Such a surface can be realized as a double cover of $\bbP^2$
branched along the union of 6 lines in a
general position. The blow-up of the dual set of 6 points in $\bbP^2$ is a
nonsingular cubic surface.
Recent work of D.\ Allcock, J.\ Carlson and D.\ Toledo \cite{ACT} gives a
different uniformization of the moduli space of
cubic surfaces where the domain of type IV is replaced by a complex ball.
This ball quotient is the moduli space of
principally polarized abelian varieties of dimension 5 with complex
multiplication in the Eisenstein ring
$\bbZ[\zeta_3]$. Each such variety can be realized as the
intermediate Jacobian of the triple cyclic cover of $\bbP^3$
branched over a nonsingular cubic surface.
Independently this construction was found by the second author
and B. Hunt.
Subsequently, Allcock and Freitag \cite{AF} found modular forms
on the ball quotient which embed it into a nine dimensional
projective space. Freitag \cite{F} later proved that the ideal
of the image is defined by cubic polynomials
and that the quotient ring is the full ring of modular forms.
The image variety turns out to be  isomorphic to a compactification of the moduli space of marked cubic surfaces.

A similar approach works for Del Pezzo surfaces of degree 2 and 1 which can
be realized as surfaces of degree 4 and 6
in weighted projective spaces $\bbP(1,1,1,3)$ and $\bbP(1,1,2,3)$,
respectively (see also \cite{HL} for another approach
to a complex ball uniformization of the moduli space of Del Pezzo surfaces
of degree 1). All of this is based on
the existence of an embedding of a complex ball into a Siegel domain. It is
also known that a complex ball can be
embedded into a type IV domain. For example a moduli space of lattice
polarized $K3$ surfaces admitting an automorphism
of order 3 or 4 which acts non-trivially on the lattice of transcendental
cycles is parametrized by an arithmetical
quotient of an open subset of a complex ball. This observation was used by
the third author \cite{Ko1} and independently
by the second author (unpublished) to construct a complex ball
uniformization of the moduli space of Del Pezzo surfaces
of degree 2. This moduli space is isomorphic to the moduli space of
non-hyperelliptic curves of genus 3 via the map which
associates to a Del Pezzo surface the fixed curve of the Geizer involution.
The $K3$ surface associated to such a
surface is its double cover branched along this fixed curve. In \cite{Ko2}
a similar description of the moduli spaces of
curves of genus 4 and of Del Pezzo surfaces of degree 1 is given.

In this paper we give a similar construction for the moduli space of cubic
surfaces. To each stable cubic surface
$S$ we associate a $K3$ surface $X_S$ with an automorphism of order 3. Its
periods are parametrized by a complex 4-ball
and we do in fact recover most of the results from \cite{ACT}. Our
construction is also closely related to the work of
K.\ Matsumoto and T.\ Terasoma \cite{MT} who associate to a line on a cubic
surface a certain curve $C$ of genus 10
which admits an involution $\sigma$ with two fixed points such that the
Prym($C,\sigma)$ is isomorphic to the
intermediate Jacobian of the triple cover of $\bbP^3$ branched along the
cubic surface. The curve $C$ also admits an
automorphism $\tau$ of order 6 such that $\sigma = \tau^3$. The $K3$
surface associated to the cubic is the minimal
nonsingular model of the quotient $(C\times E)/\langle\tau\rangle$, where
$E$ is an elliptic curve with an automorphism
of order 6. The branching of the map $C\to C/\langle\tau\rangle \cong
\bbP^1$ is very special, we have 7 branch points,
5 of which have ramification index $(3,3)$ and two have index $(6)$.
According to Deligne-Mostow \cite{DM} the moduli
space of such covers is isomorphic to an open subset of a complex ball
quotient $\calB/\Gamma$. We identify this moduli
space with the moduli space of $K3$ surfaces $X_S$ and interprete the
monodromy group $\Gamma$ in terms of the
orthogonal group of the lattice of transcendental cycles on the $K3$
surfaces. We also give an interpretation of a
compactification of the ball quotient in terms of $K3$ surfaces.

Here is the review of the contents of the paper.
In section \ref{nodalcubs} we study stable cubic surfaces.
Since these have at most nodes as singularities we refer to them
as nodal cubic surfaces. We define markings of these cubics and
we introduce the moduli space of marked nodal cubic surfaces
$\calM^m_{\ncub}$. The Weyl group $W(E_6)$ acts on $\calM^m_{\ncub}$
(the action can be described by Cremona transformations) and
the quotient variety is $\calM_{\ncub}$,
the moduli space of stable cubic surfaces.
It has a natural compactification $\overline{\calM}_{\ncub}$,
the moduli space of semi-stable cubic surfaces,
which is obtained by adding one point.
The moduli space $\calM^m_{\ncub}$ also admits a natural
compactification $\overline{\calM}^m_{\ncub}$
which is obtained by adding
40 points. It admits a $W(E_6)$-equivariant embedding
into $\bbP^9$.
We discuss different constructions of the moduli space
$\overline{\calM}^m_{\ncub}$.

For a nodal cubic surface and a line on it we define
in section \ref{sec2+5}
a pair of binary forms, of degree $2$ and $5$, modulo the
action of $SL(2)$. Using this, we
prove that the moduli space of cubic surfaces together
with a choice of a line on it is a rational variety.

In section \ref{sec2} we define a $K3$ surface $X_{S,l}$
associated to a nodal cubic surface $S$ together with the
choice of a line $l$ on $S$.
The surface $X_{S,l}$ admits a natural elliptic fibration
as well as an automorphism of order three.
We show that this $K3$ surface depends only on $S$
(and not on the choice of $l$) by defining a $K3$ surface
$X_{S,l,m}$, where $l$ and $m$ are skew lines on $S$,
which can be seen to be isomorphic to both $X_{S,l}$ and $X_{S,m}$.
We write $X_S$ for the (isomorphism class of such a) $K3$ surface
associated to $(S,l)$.
We relate $X_S$ to the $K3$ surface
associated to a cubic fourfold with a plane,
to the cubic threefold $V$
associated to $S$ by Allcock, Carlson and Toledo
and to the `Matsumoto-Terasoma curve' $C$.

In section \ref{piclat} we show that the
Picard lattice of a generic $X_S$ is isomorphic to the lattice
$M = U\oplus A_2^{\oplus 5}$. The lattice of transcendental
cycles is isomorphic to the lattice $T =A_2(-1)\oplus A_2^{\oplus 4}$.
This follows from the fact that the elliptic fibration
on the generic $X_S$ has 5 singular fibres of type $IV$ and 2 fibres of
type $II$ and some lattice theoretic considerations. We also compute the
Picard lattices of  the $K3$ surfaces associated to general nodal cubic
surfaces.

In section \ref{modulik3}
we study the moduli space of $M$-polarized $K3$ surfaces
$(X,\phi:M\to \Pic(X)$). If $\phi(M) = \Pic(X)$,
a $M$-marking $\phi$ is equivalent to the data
which consists of an elliptic fibration with a unique section, an order on
the 5 reducible fibres of type $IV$
or $I_3$, and an order on the set of irreducible components of each fibre
which do not meet the section. An $M$-marking
on the $K3$ surface $X_S$ associated to a smooth cubic surface $S$
is equivalent to a marking on $S$, that is, an order on the
set of 27 lines (or, equivalently, a choice of an ordered set of six skew
lines).
The $M$-polarized $K3$ surfaces
$(X_S,\phi)$ are distinguished from general $M$-polarized $K3$ surfaces by
the property that there exists an
automorphism $\sigma$ of order 3 which is the identity on $\phi(M)$ and
coincides with some explicitly described
isometry $\rho$ on the orthogonal complement of $\phi(M)$ in
$H^2(X_S,\bbZ)$ for smooth $S$.
The isometry $\rho$ fixes the period
$H^{2,0}(X_S)$ of $X_S$ so that the image of the period map of the
surfaces $X_S$ lies in the fixed locus of a certain
automorphism of order 3 on the period space of $M$-polarized $K3$ surfaces.
This fixed locus turns out to be
isomorphic to a 4-dimensional complex ball $\calB$. In this way we
construct the moduli space $\calK3_{M,\rho}^m$ of
$(M,\rho)$-marked $K3$ surfaces as a quotient of $\calB$. The Weyl group
$W(E_6)$ acts naturally on $\calK3_{M,\rho}^m$
by changing the markings.

In section \ref{cbu}
we establish a
natural $W(E_6)$-equivariant isomorphism from the moduli space of marked
nonsingular cubic surfaces $\calM_{\cub}^m$
onto an open subset $\calK3_{M,\rho}^m\setminus \Delta^m$ of
$\calK3_{M,\rho}^m$.
The moduli space of
isomorphism classes of pairs $(S,l)$ of cubic surfaces together with a
choice of a line is isomorphic to
the quotient of $\calK3_{M,\rho}^m\setminus \Delta^m$ by a subgroup of $
W(E_6)$ isomorphic to $W(D_5)$.
In this way we obtain an interpretation of a line on a general cubic
surface $S$ as a choice,
up to automorphisms of $X_S$, of an elliptic pencil with 5 fibres of type
$IV$ on the associated $K3$ surface $X_S$.

In section \ref{discriminant} we study in detail the geometry
of the discriminant locus $\Delta^m$.
We show that each point $[(X,\phi)]\in \Delta^m$ admits
an automorphism $\sigma$ of order 3 such that
$H^2(X,\bbZ)^{\sigma^*}$ contains $\phi(M)\oplus R$,
where $R$ is spanned by all $(-2)$-vectors in
$\phi(M)^\perp\cap \Pic(X)$.
The lattice $R$ is isomorphic to $r$ ($\le 4$)
copies of the root lattice $A_2$.
The marking $\phi$ defines an elliptic fibration on $X$
and we describe its possible singular fibres.
We also prove that $\Delta^m$ consists of 36 irreducible components
on which $W(E_6)$ acts transitively.
The cubic surfaces with Eckardt points define another divisor
in $\calK3_{M,\rho}^m$ and we prove that it consists of
45 irreducible components permuted transitively by $W(E_6)$.
Finally we show that the Satake-Baily-Borel compactification of
$\calK3_{M,\rho}^m$ contains 40 cusps,
again transitively permuted by $W(E_6)$.
This agrees with the results obtained in \cite{ACT}.

In section \ref{extiso} we show that the
$W(E_6)$-equivariant isomorphism from $\calM_{\cub}^m$
onto $\calK3_{M,\rho}^m\setminus \Delta^m$ can be extended
to an equivariant isomorphism
from the moduli space of marked nodal cubics $\calM^m_{\ncub}$
to $\calK3_{M,\rho}^m$.
We also show that
the quotient $\calK3_{M,\rho}^m/W(D_5)$ and the moduli space
of nodal cubic surfaces together with
a choice of a line $\calM_{\ncub}\li = \calM^m_{\ncub}/W(D_5)$
are isomorphic.
Moreover,
the latter space is naturally isomorphic to  the GIT-quotient
$P_1(2^5,1,1)/S_5\times S_2 = (\bbP^1)^7/\!/\SL(2)\times (g)$,
where the
linearization of $\SL(2)$ is defined by weighting the first five
factors with weight 2 and the last two factors with weight 1.
Here $S_5$ acts by permutation of the first five factors and
$S_2$ acts
by permutations of  the last two factors.

The configuration space
$P_1(2^5,1,1)/(g) = (\bbP^1)^7/\!/\SL(2)\times (g)$
occurs in the work of Deligne and Mostow \cite{DM} and
we show that the group $\Gamma$ is isogenous to the
reflection group $\Pi$ acting on $\calB$ which is generated by
the reflection group $\Pi'$ of the hypergeometric function
defined by the multi-valued form
$ \omega = z^{-1/6}[(z-1)(z-a_1)(z-a_2)(z-a_3)(z-a_4)]^{-1/3}dz $
and an involution $g$. Moreover, we match the types of degeneration
of the elliptic fibration corresponding to the marking and the type
of degeneration of a stable point set through this morphism.

Finally, in section \ref{halftwist}, we
compare the Hodge structure on the $K3$ surface $X_S$ with the
principally polarized Hodge structure on $H^1(P,\bbZ)$,
where $P$ is the intermediate Jacobian of a cubic threefold
associated to the cubic surface $S$.

\section{Nodal cubic surfaces} \label{nodalcubs}

\subsection{Nodal cubics and points in $\bbP^2$}\label{2.1}
A \emph{nodal cubic surface}
is a  surface of degree 3 in $\bbP^3$ which has
at most ordinary double points as singularities.
Let  $S\subset \bbP^3$
be a  nodal cubic surface with a node
$P=(0,0,0,1)$. Then its equation is of the form:
\begin{equation}
\label{singcub}
   F_2(X_0,X_1,X_2)X_3+F_3(X_0,X_1,X_2)=0,
\end{equation}
where the $F_i$ are homogeneous of degree $i$ and $F_2=0$ defines a smooth
conic. Projection from $P$ is a
birational isomorphism $S-\rightarrow \bbP^2$ with inverse given by: $$
\bbP^2\,-\,\rightarrow S,\qquad
x=(x_0:x_1:x_2)\longmapsto (F_2(x)x_0,\,F_2(x)x_1,\,F_2(x)x_2,\,-F_3(x)). $$
It is a rational  map given by the linear system of cubics through
$B=(F_2=0)\cap(F_3=0)$. The inverse image of a point in $B$ is a line on
$S$. There are at most two  nodes on a line in $S$ which implies that
each point in $B$ has multiplicity at most $2$. In particular, $S$ has at
most $4$ nodes.

Let $S$ be a
nodal cubic surface and let $\tilde{S}\rightarrow S$ be the
desingularization of $S$. The fibre over a node is a $(-2)$-curve, i.e. a
smooth rational curve with selfintersection $-2$.
The rational map $S-\rightarrow \bbP^2$ defines a morphism
$\pi:\tilde{S}\rightarrow \bbP^2$ which is the composition of birational
morphisms
$$
\pi:\tilde{S}=\tilde{S}_0\to \tilde{S}_1\to\cdots \to \tilde{S}_6
 = \bbP^2,
$$
where each $\pi_i:\tilde{S}_{i-1}\to \tilde{S}_{i}, i= 1,\ldots,6,$
is the blow-down of an exceptional curve of the first kind
(a $(-1)$-curve for short).

Let
$E_i'\subset \tilde{S}_i$ be the exceptional curve of $\pi_i$ and put
$E_i = (\pi_{i-1}\circ\cdots\circ\pi_1)^{-1}(E_i')$.
Let $ e_i$ be the
divisor class of $E_i$ and
let $e_0$ be the divisor class of the pre-image of a line $l\subset
\bbP^2$ under $\pi$.
The classes $e_0,e_1,\ldots,e_6$ form an orthonormal basis in $$
H^2(\tilde{S},\bbZ)
=\Pic(\tilde{S})=\bbZ e_0\oplus\bbZ e_1\oplus\cdots\oplus \bbZ e_6 $$
in the sense that $e_0^2 = 1, e_i^2 = -1, i\ne 0, (e_i,e_j) = 0,
i\ne j$.
The canonical class $K_{\tilde{S}}$ of $\tilde{S}$ is equal to
$-3e_0+e_1+\cdots+e_6$.

The anti-canonical map $\tilde{S}\rightarrow \bbP^3$ maps $\tilde{S}$ onto
$S$ and contracts the $(-2)$-curves to nodes.
In particular, $K_{\tilde{S}}$ is orthogonal to the class of each
$(-2)$-curve. Such a class is, up to sign,
one of the following $36$:
$$
e_i-e_j,\qquad e_0-e_i-e_j-e_k,\qquad
2e_0-e_1-e_2-\cdots-e_6,
$$
with $1\leq i<j<k\leq 6$. Let $p_i=\pi(E_i)\in \bbP^2$. Then $e_i-e_j, i > j,$
is effective iff
$p_i$ and $p_j$ coincide, $e_0-e_i-e_j-e_k$ is effective iff the points
$p_i,p_j$ and $p_k$ are on a line and
$2e_0-e_1-e_2-\cdots-e_6$
is effective if and only if the six points $p_1,\ldots,p_6$ are on a conic.

It follows easily from considering equation \eqref{singcub} that
other nodes of $S$ appear only when the cubic defined by $F_3$ is
simply tangent to the conic defined by $F_2$. Equivalently, $S$ can be
obtained as the blow-up 6 points on a conic,
where among the points there could
be infinitely near points of order at most 2.

\subsection{Geometric markings}
A minimal resolution of a nodal cubic surface
is a Del Pezzo surface of degree 3. Recall that a Del Pezzo surface of
degree $d$ is a smooth surface $X$ with $-K_X$
nef and $K_X^2 = d > 0$. For $d\ge 3$, the anti-canonical linear system
$|-K_X|$ maps $X$ birationally to a surface of
degree $d$ in $\bbP^d$ with at most rational double points as singularities.
Notice that we do not assume that $-K_X$
is ample, in that case one should call $X$ a Fano surface.
It is known that a Del Pezzo surface admits a birational
morphism $\pi:X\to \bbP^2$ as in \ref{2.1}.
A choice of  such $\pi$
and its decomposition $\pi = \pi_6\circ\ldots\circ
\pi_1$ is called a
\emph{geometric marking} of $X$.
Two geometric markings $X =X_0\to X_1\to\ldots \to
X_6 = \bbP^2$ and
$X '=X_0'\to X_1'\to\ldots \to X_6' =
\bbP^2$ are called {\it isomorphic} if there exist
isomorphisms $\phi_i:X_i\to X_i', i= 0,\ldots,6,$ such that
$\pi_{i+1}'\circ\phi_i= \phi_{i+1}\circ\pi_{i+1}, i = 0,\ldots,5$.

\subsection{Lattice markings}\label{lmnc} The Picard lattice of
a Del Pezzo surface $X$ of degree $d$ is isomorphic to
$$
I_{1,9-d} = <1>\oplus <-1>^{9-d},
$$
the standard odd unimodular hyperbolic lattice with the standard
orthonormal basis $(\bfe_0,\ldots,\bfe_{9-d})$.
Let $k=-3\bfe_0+\bfe_1+\ldots+\bfe_{9-d}$. Let $k^\perp$ be the orthogonal
complement  of $\bbZ  k$ in $I_{1,9-d}$. Assume $d\le 6$. Then the
sublattice $k^\perp$ is isomorphic to $Q(E_{9-d})$, where $E_{9-d}$ is
the root lattice $E_{9-d}$ if $d = 1,2,3$, the root lattice $D_5$ if $d =
4$, the root lattice $A_4$ if $d = 5$, and the root lattice $A_2+A_1$ if
$d = 6$,  formed by  vectors $\bfe_0-\bfe_1-\bfe_2-\bfe_3,
\bfe_1-\bfe_2,\ldots,\bfe_{9-d+1}-\bfe_{9-d}.$
A \emph{lattice marking} of a Del Pezzo surface $X$ is an isometry
$$
\phi:I_{1,9-d}\longrightarrow \Pic(X), \qquad \textrm{such
that}\;\phi(k)=K_{X}.
$$
In particular, the restriction of $\phi$ to $k^\perp$ is an isometry
$k^\perp\to K_{X}^\perp$.

A geometric marking defines a lattice marking by $\phi(\bfe_i)=e_i$
with $e_i$ as in \ref{2.1}.

Let $W(X)$ be the subgroup of
the orthogonal group of $\Pic(X)$ generated by reflections in the
classes of the $(-2)$-curves on $X$.
Two lattice markings
$\phi,\phi':I_{1,{9-d}}\to \Pic(X)$
are called {\it equivalent} if there exists an element
$\sigma\in W(X)$
such that
$\phi = \sigma\circ\phi'$.

The proof of the following result can be found in \cite{Loo}.

\subsection{Proposition} {\it Let $X$ be a Del Pezzo surface.
Then there is a natural bijection between the isomorphism classes
of geometric markings and
equivalence classes of lattice markings on $X$.}

\subsection{The moduli space of marked smooth cubics.}
\label{Mcubm}
We denote by $\calM_{\cub}^m$ the moduli space of
marked smooth cubic surfaces. Its points correspond to
isomorphism classes of pairs
$(S,\phi)$, where $S$ is a smooth cubic surface
and $\phi$ is a lattice marking of $S$.
There is an isomorphism:
$$
\calM_{\cub}^m\longrightarrow
\left((\bbP^2)^6-\Delta\right)/\!/\SL(3),\qquad
(S,\phi)\longmapsto (p_1,\ldots,p_6)
$$
where the $p_i\in\bbP^2$ are the images of the lines with classes
$\phi(\bfe_i)\in \Pic(S)$ in the blow down $\bbP^2$ of $S$
and $\Delta$ is the set of 6-tuples of points where either two points
coincide, or three are on a line or all six are on a conic.
The inverse image of a 6-tuple consists of the surface $S$
obtained by blowing up the $p_i$ and the marking
is defined by putting $\phi(\bfe_i)$ equal to the class
of the exceptional divisor over $p_i$.

\subsection{The Cremona action on $\calM_{\cub}^m$}
\label{nodallocus}
The Weyl group $W(E_6)$ is the subgroup of $O(I_{1,6})$
which fixes the element $k\in I_{1,6}$. It
acts naturally on $\calM_{\cub}^m$
by composing a lattice marking with (the inverse of) an isometry in
$W(E_6)$:
$$
W(E_6)\longrightarrow \Aut(\calM_{\cub}^m),\qquad
\sigma\longmapsto [(S,\phi)\longmapsto (S,\phi\circ \sigma^{-1})].
$$
Equivalently, $W(E_6)$ acts via the Cremona action
on 6 ordered points in $\bbP^2$ (see \cite{DO}).
We will simply identify $W(E_6)$ with its image in
$\Aut(\calM_{\cub}^m)$ from now on.

The quotient of $\calM_{\cub}^m$ by $W(E_6)$ is the
moduli space of smooth cubic surfaces $\calM_{\cub}$.
Let $p_{\cub}$ be this quotient map:
$$
p_{\cub}:\calM_{\cub}^m\longrightarrow
\calM_{\cub}^m/W(E_6)\,\cong\, \calM_{\cub}.
$$

\subsection{The GIT compactification}
Geometric Invariant Theory provides a natural compactification
of the moduli space of cubic surfaces $\calM_{\cub}$:
\[
\overline{\calM}_{\cub} =
\bbP(H^0(\bbP^3,\calO_{\bbP^3}(3))^{\ss}/\!/\SL(4).
\]
The stable points in
$\bbP(H^0(\bbP^3,\calO_{\bbP^3}(3))$ are the isomorphism classes
of nodal cubic surfaces. Points in
\[
\calM_{\ncub} = \bbP(H^0(\bbP^3,\calO_{P^3}(3))^{\st}/\!/\SL(4)
\]
are thus
isomorphism classes of nodal cubics.
The strictly semi-stable points all map to one point in
$\overline{\calM}_{\cub}$. The complement of this point in
$\overline{\calM}_{\cub}$ is denoted by $\calM_{\ncub}$, the moduli
space of nodal cubic surfaces.

The explicit computations of invariants of cubic quaternary forms,
due to A. Cayley and G. Salmon \cite{Salmon},
gives an isomorphism
\begin{equation}\label{CS}
  \overline{\calM}_\cub\cong \bbP(1,2,3,4,5).
\end{equation}
The moduli space of nonsingular surfaces is isomorphic to
the complement of a hypersurface of degree 4
defined by the discriminant.
In particular, ${\calM}_{\cub}$ is affine.

\subsection{Moduli of marked nodal cubics}\label{modnod}
We can construct the moduli space of marked nodal cubic surfaces as
follows. Let $k(\calM_{\cub}^m)$ be the field of rational functions
of $\calM_{\cub}^m$. It is an extension, with Galois group
$W(E_6)$,
of $k(\calM_{\cub})=k(\overline{\calM}_{\cub})$.
Now we define $\overline{\calM}_{\cub}^m$ to be the normalisation
of $\overline{\calM}_{\cub}$ in the field $k(\calM_{\cub}^m)$.

By its definition, $\overline{\calM}_{\cub}^m$ is a normal projective
variety and, since $\calM_{\cub}^m$ is smooth (see sections
\ref{naruki} and \ref{gitmodel}),
we have
$$
\calM_{\cub}^m\hookrightarrow \overline{\calM}_{\cub}^m,
$$
the complement of $\calM_{\cub}^m$ will be called the boundary
of $\overline{\calM}_{\cub}^m$.
By construction, the Weyl group $W(E_6)$ acts on
$\overline{\calM}_{\cub}^m$ with quotient $\overline{\calM}_{\cub}$:
$$
\bar{p}_\cub:\overline{\calM}_{\cub}^m  \longrightarrow
\overline{\calM}_{\cub}=\overline{\calM}_{\cub}^m/W(E_6)
$$
and $\bar{p}_\cub=p_\cub$ on the subvariety $\calM_{\cub}^m$.
Finally we define the moduli space of marked nodal cubic surfaces
to be:
$$
{\calM}_{\ncub}^m:=\bar{p}^{-1}(\calM_{\ncub}).
$$
This moduli space is the complement of a finite set of points,
called the cusps, in $\overline{\calM}_{\cub}^m$
and the cusps are all in one $W(E_6)$-orbit.

Despite its abstract definition, the variety
$\overline{\calM}_{\cub}^m$ is rather well-known.
Below we present some other constructions of it, and we show that the
points in ${\calM}_{\ncub}^m$ correspond to isomorphism
classes of marked nodal cubic surfaces. We do not know whether
${\calM}_{\ncub}^m$ is the (coarse) moduli space of some functor.

\subsection{Naruki's model}\label{naruki}
In \cite{Naruki}, Naruki  constructs a smooth, projective
compactification of the moduli space
$\calM_{\cub}^m$ which he calls the cross-ratio variety.
Its boundary contains 40 divisors which can be blown down to 40
singular points of a variety $\calN$. The variety $\calN$
is normal (it is smooth outside the 40 singular points and $\calN$
is also normal in the singular points, where $\calN$ is locally
a cone over a Segre embedding of $(\bbP^1)^3$).
Hence the isomorphism on the open subsets,
isomorphic to $\calM_{\cub}^m$,
of $\calN$ and $\overline{\calM}_{\cub}^m$ extends to
a birational morphism
$$
\phi_\calN:\calN\longrightarrow \overline{\calM}_{\cub}^m.
$$
Naruki also shows that the
action of $W(E_6)$ on $\calM_{\cub}^m\;(\subset \calN)$
extends to a biregular action on $\calN$ with quotient
$\calN/W(E_6)= \overline{\calM}_{\cub}$.
This implies that
$\phi_\calN$
is a bijection. Hence, by Zariski's Main Theorem,
$\phi_\calN$ is an isomorphism.

From Naruki's description of $\calN$, see also \cite{CG},
one obtains that the forty singular points of $\calN$ map to the
cusps of $\overline{\calM}_{\cub}^m$. Moreover, the boundary of
$\overline{\calM}_{\cub}^m$ consists of $36$ divisors, each of
which is isomorphic to the Segre cubic threefold $\calS_3$,
best seen as a subvariety of $\bbP^5$:
$$
\calS_3:\qquad \sum_{i=1}^6x_i=0,\quad \sum_{i=1}^6x_i^3=0.
$$
The group $W(E_6)$ acts transitively on the set of
36 boundary divisors.

\subsection{Boundary divisors.}\label{Boundary divisors}
The 36 boundary divisors are parametrized by the 36 positive simple
roots of $E_6$.
Let $\alpha$ be a positive simple root,
so
$\alpha = e_i-e_j$, $e_0-e_i-e_j-e_k$, $2e_0-e_1-\ldots-e_6$,
with $i > j>k$,
note that indeed $\alpha\in k^\perp\subset I_{1,6}$.
If $(S,\phi)$ is a marked nodal cubic surface
such that $\phi(\alpha)$ is effective, then $S$ has a node,
cf.\ \ref{2.1}.
To each $\alpha$ we assign the divisor
$D_\alpha$ in $\calM_{\ncub}^m$, we write:
\[D_\alpha = \begin{cases}
     D_{ij}  & \text{if}\quad \alpha = e_i-e_j \\
       D_{ijk}& \text{if}\quad \alpha = e_0-e_i-e_j-e_k\\
D_0& \text{if}\quad \alpha = 2e_0-e_1-e_2-e_3-e_4-e_5-e_6.
\end{cases}
\]

Each $D_\alpha$ is isomorphic to the Segre cubic $\calS_3$
which is a projective model of the GIT-quotient
$(\bbP^1)^6/\!/\SL(2)$
\cite{DO}. This implies that the points of $\calM_\ncub^m$
correspond to marked nodal cubics. In fact, for the points in
$\calM_\cub^m$ this is obvious, if a point $x$
is in a boundary divisor $D$
then it corresponds to a sixtuple of points on $\bbP^1$ which gives
a geometrical marking of the corresponding nodal cubic, in fact
the boundary divisor corresponds a $(-2)$-curve on the nodal
cubic surface $\bar{p}(x)\in\calM_\ncub$ and together with
the sixtuple of points gives a geometrical marking of
the nodal cubic, see also the examples in \ref{lines}.

In particular, $D_\alpha$ parametrizes marked nodal cubic surfaces
$(S,\phi)$ for which $\phi(\alpha)$ is effective.
If $\phi(\alpha)$ is effective and $r_\alpha$ denotes the
reflection in
$W(E_6)$ defined by the root $\alpha$,
then the lattice marked nodal cubic surfaces $(S,\phi)$ and
$(S,\phi\circ
r_\alpha)$ are equivalent. This suggests that in the Cremona action
of $W(E_6)$ on $\calM_{\ncub}^m$ the reflection $r_\alpha$ acts
identically on $D_\alpha$. This is in fact the case
(\cite{Naruki}, p.\ 22).

The Segre cubic has 10 nodes $p_{1ij}$,
for example, $p_{125}=(1:1:-1:-1:1:-1)$,
corresponding to the  minimal
orbit of sixtuples $(p_1,\ldots,p_6)$ of points on $\bbP^1$ such
that $p_i=p_j=p_k$, $p_l =p_m=p_n$.
Identifying $\calS_3$ with $D_0$,
the nodes of $\calS_3$ are the cusps of
$\overline{\calM}_{\cub}^m$ lying on $D_0$.

A boundary divisor $D$ is invariant with respect to
a subgroup of $W(E_6)$ isomorphic to $S_6$, under the
isomorphism $D\cong \calS_3$ the $S_6$ acts by permutation
of the coordinates.
On $D_0$, which parametrizes 6 points on a conic, this $S_6$ acts
by interchanging the 6 points.

The image $\bar{p}(D)$ of a boundary divisor in
$\overline{\calM}_{\cub}$ is the locus of singular cubic surfaces.
It is defined by the vanishing of the discriminant invariant on the
space of cubic surfaces, which is of degree 32
in the coefficients of the cubic form. In the
isomorphism \eqref{CS} it corresponds to the hyperplane defined by
the unknown with weight 4. Thus it is isomorphic to $\bbP(1,2,3,5)$.
It is also known that  the quotient of the Segre cubic by $S_6$ is
isomorphic to $\bbP(1,2,3,5)$.
Hence the restriction of  the map
$\bar{p}:\overline{\calM}_{\cub}^m\to \overline{\calM}_{\cub}$
to a boundary divisor
$D$ is the quotient map $D\to D/S_6$.

\subsection{A GIT model}\label{gitmodel}
We sketch another construction of $\overline{\calM}_\cub$.
First we recall the explicit construction of the GIT-quotient
$X=(\bbP^2)^6/\!/\SL(3)$ given in
\cite{DO}. The graded ring of invariants
$$
R = \bigoplus_{n=0}^\infty \big(H^0((\bbP^2)^6,\otimes_{i=1}^6
\pi_i^*\calO_{\bbP^2}(n)\big)^{\SL(3)}
$$
is generated by elements $t_0,t_1,t_2,t_3,t_4$ of degree 1 and one
element $t_5$ of degree 2.
Here $\pi_i$ is the $i$-th projection from $(\bbP^2)^6$.
The relation between the generators is
$t_5^2+F_4(t_0,t_1,t_2,t_3,t_4) = 0$, where $F_4$ is a homogeneous
polynomial of degree 4.  Thus $X$ is isomorphic to a hypersurface of
degree 4 in the weighted projective space $\bbP = \bbP(1,1,1,1,1,2)$.

The quartic 3-fold $V(F_4)$ in $\bbP^4$ has 15 double lines $l_{ij}$
corresponding to minimal semi-stable orbits of points sets
$(p_1,\ldots,p_6)$ where
$p_i = p_j$.
Three lines
$l_{ij},l_{kl}, l_{mn}$, where $\{1,2,3,4,5,6\} = \{i,j\}\cup
\{k,l\}\cup\{m,n\}$, intersect at one point $P_{ij,kl,mn}$.
It represents the orbit of the point set $p_i =p_j$, $p_k =p_l$,
$p_m=p_n$.
It follows from the explicit equation of $F_4$ that  its local
equation at  $P_{ij,kl,mn}$ is given by $w^2+z_1z_2z_3=0$, where $w =
z_i=z_j=0$
is the local equation of one of the 3 double lines meeting at the
point. This implies that $X$ is given locally at the point
$P_{ij,kl,mn}' = (P_{ij,kl,mn},0)$ by the equation $uv+xyz = 0$.

Let $Z$ be the singular locus of $X$ and $\calI_Z$ its sheaf of
ideals. One considers the linear system $ |\calI_Z(3)| \subset R_3$.
A. Coble \cite{Coble} gives explicitly 40 elements of $|\calI_Z(3)|$
which span a $\bbP V\cong\bbP^9$ and shows that the birational
action of $W(E_6)$ on $X$ induces a linear action on $V$.
We construct the moduli space of
marked cubic surfaces as the image $Y$ of $X$ under the rational map
given by the linear system $\bbP V$.

First we blow up the ambient space $\bbP$ at the points
$P_{ij,kl,mn}'$. Let $E_{ij,kl,mn} \cong \bbP^4$ be the exceptional
locus. The proper inverse transform $X_1$ of $X$ intersects  each
$E_{ij,kl,mn}$ along the union of two hyperplanes
$H_{ij,kl,mn},H_{ij,kl,mn}'$ corresponding to the tangent cone of the
singular point. The proper inverse transforms of the lines $l_{ij}$ are double curves $C_{ij}$ on $X_1$.  Each
of the curves $C_{ij}, C_{kl}, C_{mn}$ intersects $E_{ij,kl,mn}$ at a
point. The three points span the plane $\Pi_{ij,kl,mn} =
H_{ij,kl,mn}\cap H_{ij,kl,mn}'$.  Next we blow up the 15 singular
curves $C_{ab}$ to get a  variety $X_2$. 
The proper inverse transform of the linear system $\bbP V$ in $X_2$
has base locus equal to the union of the proper transforms $\bar{\Pi}_{ij,kl,mn}$ of the
planes $\Pi_{ij,kl,mn}$. Each surface $\bar{\Pi}_{ij,kl,mn}$ is isomorphic to the blow-up of 3 points on the plane. The proper transforms of the lines joining three pairs of points are double curves of $X_2$. Next we blow up the surfaces $\bar{\Pi}_{ij,kl,mn}$ to get a nonsingular variety $X_3$. Now the proper inverse transforms of the hyperplanes
$H_{ij,kl,mn},H_{ij,kl,mn}'$ become separated and the proper inverse transform of the linear system $\bbP V$ has no base points.

Let $Y\subset \bbP^9$ be the image of $X_3$ under this linear system.
Observe first that $Y$ is a compactification of the
geometric quotient $\calM^m_\cub=U/\SL(3)$,
where $U=(\bbP^2)^6-\Delta$ as in \ref{Mcubm}.

Next we shall see its complement.
First of all we have 20 divisors $D_{ijk}'$ in $X$
representing 6-tuples of points where $p_i,p_j,p_k$ are collinear.
The sum of the two divisors $D_{ijk}'$
and $D_{lmn}'$, where $\{i,j,k\}\cup \{l,m,n\} = \{1,\ldots,6\}$, is
defined by a linear function $L_{ijk}=L_{lmn} \in R_1$ (see
\cite{DO}).
The corresponding hyperplane $V(L_{ijk})$ cuts out the quartic
$V(F_4)$ along a nonsingular quadric $Q_{ijk} = Q_{lmn}$.
The quadric contains 6 double lines
$l_{ij},l_{ik},l_{jk},l_{lm},l_{ln},l_{mn}$.
Let $D_{ijk}$ be the proper inverse transforms of $D_{ijk}'$ in $Y$.
Let $D_{ij}$ be the proper inverse transforms in $X_3$ of the
pre-images of the curves $C_{ij}$ in $X_2$. We have 15 such divisors.
Finally, let $D_0$ be the proper inverse transform of  $V(t_5)\cong
V(F_4)$ in $Y$.
It is easy to see that under the map $X_3\to Y$ the
proper inverse transforms of the quadrics $Q_{ijk}$ are blown down to
points $c_{ijk} = c_{lmn}$. Also let $c_{ij,kl,mn}, c_{ij,kl,mn}'$ be
the images in $Y$ of the hyperplanes $H_{ij,kl,mn}, H_{ij,kl,mn}'$.
Altogether we have 40 points which we call the cusps. The forty
cusps is the set of singular points of the variety $Y$.
So, we see that the complement of the image of $U/\SL(3)$ in $Y$ is
equal to the union of 36 divisors $D_{ijk},D_{ij},D_0$.

The Weyl group $W(E_6)$ acts on $Y$ interchanging the boundary
divisors. This makes them all isomorphic to each other.
This is easy to check.
The restriction of the linear system $\bbP V$ to the
quartic $V(F_4)$ is the map given by the partials of $F_4$. It maps
$V(F_4)$ to the dual variety known to be isomorphic to the Segre
cubic $\calS_3\subset \bbP^4$.
This shows that $D_0\cong \calS_3$.

One can check that the variety $Y$ is a normal proper
$W(E_6)$-variety containing the $W(E_6)$-variety
$\calM_{\cub}^m$ as an open subset. Thus there is a
birational morphism $f:Y\to \overline{\calM}_{\cub}^m$.
We claim that $f$ is an isomorphism.
Let $E$ be an irreducible component of the
exceptional locus of $f$.
It is contained in one of the 36  boundary divisors $D$.
However $D\cong \calS_3$  has $\Pic(D)\cong \bbZ$.
Nothing can be blown down on $D$. Thus we obtain that
$$
Y\cong \overline{\calM}_{\cub}^m.
$$

\subsection{Remark.}
In \cite{ACT}, $\calM^m_\cub$ is identified with an open subset
of a smooth ball quotient. In \cite{AF} Allcock and Freitag show,
using modular forms constructed via a Borcherds lift,
that this ball quotient embeds into a $\bbP^9$ and that the closure
of its image is isomorphic to the Satake compactification
of the ball quotient, the boundary consists of 40 singular points.
Freitag \cite{F} proved that ideal of the image of the
ball quotient is generated by explicitly given cubics and that it
is a normal variety.

Coble, in \cite{Coble}, defines a rational map $(\bbP^2)^6
\rightarrow \bbP^9$ which is
$\SL(3)$-invariant and hence factors over $\calM_\cub$.
It is easily seen to be a birational
isomorphism between $\calM_\cub$ and its image.
This map is moreover
equivariant with respect to the Cremona action of $W(E_6)$.
See also \cite{Y} where in particular
the restriction to a boundary divisor is worked out.
It is easy to verify that the image of $\calM_\cub$ lies in the
subvariety defined by the cubics.

In \cite{Geemen} the corresponding
rational functions on Naruki's variety
$\calN\cong \overline{\calM}_\cub^m$ are explicitly identified,
and also the 40 functions used by Coble are given.

Matsumoto and Terasoma \cite{MT} showed how to get this embedding
via an embedding of the complex ball into the Siegel space
(of genus $5$) followed by a map to $\bbP^9$ given by explicitly
determined theta constants.

\subsection{Moduli of $i$-nodal cubics}\label{inodal}
The irreducible components of the
locus of marked nodal cubics with $i$ nodes are parametrized by
unordered subsets
of $i$ orthogonal roots (up to sign) in $E_6$.
We denote by $D_{\alpha_1,...,\alpha_i}$ the
intersection of the divisors $D_{\alpha_1},..., D_{\alpha_i}$
corresponding to $i$ orthogonal roots $\alpha_1,..., \alpha_i$.

The stabilizer in $W(E_6)$ of such a locus
$D_{\alpha_1,...,\alpha_i}$ is the product of the subgroup
of order $2^i$, generated by the corresponding $i$ roots
(this subgroup acts trivially on the component),
the permutations on $i$ roots $\alpha_1,..., \alpha_i$
($\simeq S_i$)
and the subgroup generated by reflections in the roots
orthogonal to the $i$ simple roots.
The stabilizer modulo the subgroup of order $2^i$ is the group of
permutations of geometric markings on $S$.

In case $i=1$, the 30 roots $e_i-e_j$ are all orthogonal to the root
$\alpha=2e_0-e_1-\ldots-e_6$,
so we see that ${\bbZ}/2{\bbZ} \times W(A_5)\cong
{\bbZ}/2{\bbZ} \times
S_6$ acts on $D_\alpha$.
Thus we recover the fact that $W(A_5)\cong S-6$ acts on a boundary
divisor.

In case $i=2$, there are 12 roots $e_i-e_j$, $(3\leq i, j \leq 6)$
orthogonal to the two roots
$\alpha_1=2e_0-e_1-\ldots-e_6$ and
$\alpha_2 = e_1-e_2$.  They form
a root system of type $A_3$.  So
$({\bbZ}/2{\bbZ})^2 \cdot S_2 \times W(A_3) \simeq ({\bbZ}/2{\bbZ})^2 \cdot
S_2 \times S_4$
acts on $D_{\alpha_1,\alpha_2}$.

In case $i=3$, there are
two roots $\pm (e_5-e_6)$ orthogonal to the three roots
$\alpha_1=2e_0-e_1-\ldots-e_6$,
$\alpha_2 = e_1-e_2$ and $\alpha_3 = e_3-e_4$.  They form
a root system of type $A_1$.  So
$({\bbZ}/2{\bbZ})^3 \cdot S_3 \times W(A_1) \simeq ({\bbZ}/2{\bbZ})^3 \cdot
S_3 \times {\bf Z}/2{\bf Z}$ acts
on $D_{\alpha_1,\alpha_2, \alpha_3}$.

In case $i=4$, there are no roots orthogonal to the four roots
$\alpha_1=2e_0-e_1-\ldots-e_6$,
$\alpha_2 = e_1-e_2$, $\alpha_3 = e_3-e_4$ and $\alpha_4 = e_5-e_6$.
So $({\bbZ}/2{\bbZ})^4 \cdot S_4$ acts on $D_{\alpha_1,..., \alpha_4}$.

\subsection{Lines on a nodal cubic surface}\label{lines}
A nonsingular cubic surface contains 27 lines.
They represent the classes
$e_0-e_i-e_j$, $1\le i <j \le 6$, $e_i, 2e_0-e_1-\cdots-e_6+e_i$,
$i =1,\ldots,6$.

Assume now that $S$ has a node $s_0$. Projecting from $s_0$,
we see that $\tilde{S}$ admits a geometric marking
$\pi:\tilde{S}\to \bbP^2$ such that the images $p_i$ of the
$E_i$ (as in \ref{2.1})
lie on an irreducible  conic $C$. If $S$ has no more
nodes, the six points $p_i$ are distinct.
If there is one more node, we may
assume without loss of generality that $p_2$ is
infinitely near to $p_1$ (i.e. $E_2= E_1+C$,
where $C$ is  a $(-2)$-curve and the point $p_2$ corresponds
to the tangent direction of $C$ at $p_1$).
If $S$ has three nodes we can further
assume that $p_4$ is infinitely near to $p_3$
with the similar tangency condition. Finally if $S$ has 4 nodes we
can further assume that $p_6$ is infinitely near to
$p_5$. From this we easily deduce the following facts.

If $S$ has one node, there are 21 lines on $S$.
Six of them contain the
node, and are represented by
the exceptional curves  $E_i=\phi(\bfe_i)$, where $\phi$ is the
lattice marking corresponding to the geometric marking, we will
simply omit $\phi$ in what follows.
The remaining 15 lines have the classes
$\phi(e_0-e_i-e_j)$. The $(-2)$-curve $C$ has class
$\alpha_1=2e_0-(e_1+\ldots+e_6)$ and the classes
$e_i+\alpha_1=s_{\alpha_1}(e_i)$ also
represent the lines on the node. So the lines on the nodes
are limits of pairs of lines on a smooth cubic surface.

If $S$ has 2 nodes, there are 16 lines on $S$.
The $(-2)$-curves are $\alpha_1=2e_0-(e_1+\ldots+e_6)$ and $e_2-e_1$,
the orbits on the set of classes of 27 lines
of the group generated by
$s_{\alpha_1}$ and $s_{\alpha_2}$ correspond to the lines on $S$.
One line connects
the two nodes and represents the orbit
$\{e_1,e_2=e_1+\alpha_2,e_1+\alpha_1,e_2+\alpha_1\}$.
There are 4 lines passing through the node $s_0$
which representing the orbits
$\{e_i,e_i+\alpha_1\}$, $i=3,4,5,6$. Another 4 lines pass
through the second node.
They represent the orbits $\{e_0-e_2-e_i, e_0-e_1-e_i\}$, $i =
3,4,5,6$. The remaining 7 lines do not contain
nodes. They represent orbits with one element,
given by the classes $e_0-e_i-e_j, 3\le i < j \le 6$
and $e_0-e_1-e_2$.

If $S$ has 3 nodes, there are 12 lines.
There are 3 lines connecting pairs
of nodes. They represent the classes
$e_1, e_3, e_0-e_1-e_3$. There are 6 lines each containing one node.
They represent the classes $e_5,e_6, e_0-e_1-e_i,
e_0-e_3-e_i, i = 5,6$.
The remaining 3 lines do not contain nodes.
They represent the classes $e_0-e_1-e_2, e_0-e_3-e_4, e_0-e_5-e_6$.

If $S$ has 4 nodes there are 9 lines. Six of them connect 6 pairs of nodes.
They represent the classes
$e_1,e_3,e_5, e_0-e_1-e_3,e_0-e_1-e_5,e_0-e_3-e_5$. The remaining three
lines do not contain nodes and represent the
classes $e_0-e_1-e_2, e_0-e_3-e_4, e_0-e_5-e_6$.

\subsection{Pencils of conics}\label{conics}
A conic on a nodal cubic surface $S$ is cut out by a plane.
The residual component of the plane section is a line. The pencil of planes
through this line defines a pencil of
conics. Thus the number of pencils of conics is equal to the number of
lines.  The preimage of the pencil on
$\tilde{S}$ is a
conic bundle, i.e. a morphism $f: \tilde{S}\to \bbP^1$ with general fibre
isomorphic to $\bbP^1$.  A standard
computation  shows that singular fibres of $f$ are of the following three
types:

Type $I$: $F = E_1+E_2$, where $E_1,E_2$ are two $(-1)$-curves and $E_1\cdot
E_2 = 1$.

Type $II$: $F = E_1+E_2+R$, where $E_1,E_2$ are $(-1)$-curves, $R$ is a
$(-2)$-curve,
$E_1\cdot E_2 = 0, E_1\cdot R = E_2\cdot R = 1$.

Type $III$: $F = R_1+R_2+2E$, where $R_1,R_2$ are $(-2)$-curves, $E$ is a
$(-1)$-curve,
$R_1\cdot R_2 = 0, R_1\cdot E = R_2\cdot E = 1$.

The number of singular fibres is equal to 5 if we count
the fibres of type $II$ and $III$ with multiplicity 2.

The pre-image of the line $l$ corresponding to the pencil defines a
bisection $B$ of $f$.
There are three possible cases:

No nodes on $l$: $B$ is irreducible.

One node on $l$: $B= B_0+R$, where $B_0$ is a $(-1)$-curve, $R$ is a
$(-2)$-curve, $B_0\cdot R = 1$.
Each component of $B$ is a section of $f$.

Two nodes on $l$: $B = B_0+R_1+R_2$, where $B_0$ is a $(-1)$-curve, $R_1,
R_2$ are $(-2)$-curves,
$B_0\cdot R_1 =B_0\cdot R_2 = 1$.  The components $R_1$ and $R_2$ are
sections of $f$. The component $B_0$ is contained in a fibre.

Let $p_1,\ldots,p_s\in \bbP^1$ be the points such that the fibre
$f^{-1}(p_i)$ is singular. We assign to each point
$p_i$  the multiplicity  $m_i$ equal to 2 if the fibre is of type I and
equal to 4 otherwise. The divisor $D =
\sum_{s=1}^sm_ip_i$ will be called the \emph{discriminant} of the conic
pencil. Let $p_{s+1},p_{s+2}\in \bbP^1$
be the points such that the bisection $B$ ramifies over  these points. If
$B$ is reducible, we assume that $p_{s+1} =
p_{s+2} = q$, where $B$ has a singular point over $q$. The divisor $T =
p_{s+1} + p_{s+2}$ will be called the
\emph{bisection branch divisor}. Let us write the divisor
$D+T = \sum_{i=1}^sm_ip_i+p_{s+1}+p_{s+2}$ as $\sum_{i=1}^{s'}n_ip_i,$
where $s' \le s+2$.  We order the points in such a way that  $n_1\ge n_2\ge
\ldots\ge n_{s'}$.  The vector
$\bft = (n_1,\ldots,n_{s'})$ will be called the \emph{type  vector} of the
conic pencil.

The next table lists all
possible type vectors. Also we indicate the total number $r$ of nodes on
$S$, the number $e$ of
Eckardt points on $l$
(i.e. points where three lines meet).
\begin{table}[h]
\[
\begin{array}{rlllll}
{}& \bft&{\rm Singular\; fibres}&{\rm Kodaira \quad
fibres}&r&e \\
1) & (2 2 2 2 2 1 1)&5 I&5 IV,\ 2 II&0&0\\
\noalign{\smallskip}
2) & (3 2 2 2 2 1) &5 I&I_0^*,\ 4 IV,\ II&0&1\\
    \noalign{\smallskip}
3) & (3 3 2 2 2)&5 I&2 I_0^*,\ 3 IV&0&2 \\
    \noalign{\smallskip}
4) & (2 2 2 2 2 2) &5 I&6 IV&1&0
\\ \noalign{\smallskip}
5) & (4 2 2 2 1 1) &II,\ 3 I&IV^*,\ 3 IV,\ 2 II&1&0
\\ \noalign{\smallskip}
6) & (4 3 2 2 1) &II,\ 3 I&IV^*,\ I_0^*,\ 2 IV,\ II&1&1 \\
\noalign{\smallskip}
7) & (4 3 3 2) &II,\ 3 I&IV^*,\ 2I_0^*,\ IV&1&2
\\ \noalign{\smallskip}
8) & (4 2 2 2 2) &II,\ 3 I&IV^*,\ 4 IV&2&0\\
\noalign{\smallskip}
8^*) & (4 2 2 2 2)&5 I&IV^*,\ 4 IV&2&0\\
    \noalign{\smallskip}
9) & (4 4 2 1 1) &2 II,\ I&2 IV^*,\ IV,\ 2 II&2&0  \\
\noalign{\smallskip}
10) & (5 2 2 2 1) &III,\ 3 I&II^*,\  3 IV,\ II&2&0\\
\noalign{\smallskip}
11) & (4 4 3 1) &2 II,\ I&2 IV^*,\ I_0^*,\ II&2&1 \\
\noalign{\smallskip}
12) & (5 3 2 2)&III,\ 3 I&II^*,\ I_0^*,\ 2 IV&2&1 \\
\noalign{\smallskip}
13) & (4 4 2 2) &2 II,\ I&2 IV^*,\ 2 IV&3&0\\
\noalign{\smallskip}
13^*) & (4 4 2 2) &II,\ 3 I&2 IV^*,\ 2 IV&3&0\\
\noalign{\smallskip}
14) & (5 4 2 1) &III,\ II,\ I&II^*,\ IV^*,\ IV,\ II&3&0 \\
\noalign{\smallskip}
15) & (5 4 3) &III,\ II,\ I&II^*,\ IV^*,\ I_0^*&3&1 \\
\noalign{\smallskip}
16) & (4 4 4) &2 II,\ I&3 IV* &4&0\\
\noalign{\smallskip}
17) & (5 5 2) &2 III,\ I&2 II^*,\ IV&4&0 \\
\end{array}
\]
\caption{Pencils of conics}
\end{table}

The column ``Kodaira fibres'' will be explained
later in section \ref{ellfib}.

\subsection{Types of lines}\label{types} Let $l$ be a line defining the
pencil of conics

Case 1), 2), 3) on the above Table 1: $l$ is any line.

Case 4): $l$ is one of 6 lines containing  the node.

Case 5), 6), 7):  $l$ is one of 15 lines not passing through the node.

Case 8):  $l$ is one of 8 lines through exactly one node.

Case 8*):  $l$ is the unique line containing two nodes.

Case 9), 11): $l$ is one of 6 lines not containing a node and not meeting
the line of type 8*).

Case 10), 12): $l$ is the unique line not containing a node and meeting the
line of type 8*).

Case 13):  $l$ is one of 6 lines passing exactly through one node.

Case 13*):  $l$ is one of 3 lines passing through two nodes.

Case 14, 15):  $l$ is one of 3 lines not containing a node.

Case 16): $l$ is one of 6 lines passing through two nodes.

Case 17): $l$ is one of 3 lines not containing a node.

\section{Cubic surfaces and 2+5 points on the line} \label{sec2+5}

\subsection{The forms $(F_2,F_5)$.}\label{f5f2}
Let $S$ be a nodal cubic surface and let $l$ be a line on $S$.
Consider the
pencil of conics through the line $l$, cf.\ section \ref{conics}.
Let $D = \sum_{i=1}^sm_ip_i$ be its discriminant divisor and let
$T =p_{s+1}+p_{s+2}$ be the bisection branch
divisor.  Let $F_5(x_0,x_1)$ be a homogeneous  form of degree
5 defining
$D$ and let $F_2(x_0,x_1)$ be a homogeneous
form of degree 2 defining $T$.

It follows from section \ref{conics} that the following
properties are satisfied:
\begin{itemize}
\item [(i)] $F_2\ne 0$;
\item [(ii)] $F_5$ has at most double root;
\item [(iii)] $F_2$ and $F_5$ do not have common multiple roots.
\end{itemize}

A pair of binary forms $(F_5,F_2)$ satisfying  properties (i)-(iii) will
be called a {\it stable pair}.
Let $V(d)$ be the space of binary forms of degree $d$. A pair of nonzero
binary forms  $(F_5,F_2)$ defines a point in
$\bbP(V(5))\times \bbP(V(2))$.

\subsection{Proposition}\label{stable}
{\it A pair of nonzero binary forms
$(F_5,F_2)$ is stable if and only if it is a stable point with
respect to the diagonal action of $\SL(2)$ and the linearization
defined by the invertible sheaf
$ \calO_{\bbP(V(5))}(2)\boxtimes\calO_{\bbP(V(2))}(1)$.
A semi-stable point corresponds to a pair of forms
$(F_5,F_2)$ such that either $F_5$ has a root of multiplicity 3 or
$F_5$ and $F_2$ share a double root}.

\begin{proof} This easily follows from the Hilbert-Mumford numerical
criterion of stability and is left to the reader.
\end{proof}

\subsection{Line marked cubic surfaces}\label{moduli}
Let $(S,\phi)$ be a nodal cubic surface with a
geometric marking $\phi$ on
its minimal resolution and
let $l$ be a line on $S$.
The Weyl group $W(D_5)$, a subgroup of $W(E_6)$,
acts on markings of $S$ stabilizing the line $l$.
The quotient space
$$
\calM_{\ncub}\li = \calM_{\ncub}^m/W(D_5)
$$
is the moduli space of isomorphism classes of pairs $(S,l)$,
where $l$ is a line on $S$.

To a pair $(S,l)$ we associate the binary forms $F_2,F_5$
as in \ref{f5f2}.
This defines a map
\begin{equation}\label{bin}
\calM_{\ncub}\li\longrightarrow
\bigl(\bbP(V(2))\times \bbP(V(5))\bigr)^{\textup{s}}/\SL(2),\qquad
(S,l)\longmapsto [(F_2,F_5)],
\end{equation}
where   $(\bbP(V(2))\times \bbP(V(5)))^{\textup{s}}$ is the open subset
corresponding to stable pairs of  binary forms.

\subsection{Theorem.}\label{thm1}
{\it The map \eqref{bin}  is an isomorphism.}

\begin{proof}
We have to show how to reconstruct $(S,l)$ from the $\SL(2)$-orbit  of a
pair $(F_5,F_2)$. Let us first consider the case when $S$ is nonsingular.
We view the zeroes of binary forms as the tangent directions at a fixed
point $p\in \bbP^2$ and identify these with lines in $\bbP^2$ containing
$p$.
Given $(F_2,F_5)$, fix a  conic $Q$ not containing $p$ such that the lines
through $p$
defined by $F_2$ are tangents of $Q$.
Then a choice of 5 points $p_1,\ldots,p_5$ on the intersection of the
lines defined by $F_5$ with the conic,
no two lying on the same line, defines uniquely (up to isomorphism)
a cubic surface $S$ with a line $l$ corresponding to the conic.
It is obtained by blowing up $\bbP^2$ at the points $p_1,\ldots,p_5,p$.
Let $p_i'$ be the point on $Q$ such that $p_i,p_i',p$ are collinear.
Let us show that replacing $p_i$ with $p_i'$ leads to an
isomorphic pair $(S',l')$.

Note that changing $(p_1,\ldots,p_5)$ with $(p_1',\ldots,p_5')$ leads to
the same surface because  the points $(p_1,\ldots,p_5,p)$ and
$(p_1',\ldots,p_5',p)$ are projectively equivalent. This can be easily
seen by choosing projective coordinates such that $p = (0,0,1)$ and $Q =
V(x_0x_1-x_2^2)$. Then
$p_i = (1,a_i^2,a_i)$ and $p_i' = (1,a_i^2,-a_i)$.

Now it is enough to show that changing a pair  $(p_i,p_j)$ with
$(p_i',p_j')$ defines an isomorphic surface. For this we consider
the Cremona transformation with base points at $p,p_i,p_j$. It sends the
conic $Q$ through $p_1,\ldots,p_5$ to a conic $Q'$. Chose coordinates so that
$p = (0,0,1), p_i = (1,0,0), p_j = (0,1,0)$.  Then the sides of the
triangle of base points are mapped to the opposite vertices. The
transformation looks like
$(x_0,x_1,x_2)\mapsto (x_1x_2,x_0x_2,x_0x_1)$. Now we know that the
image of $p_i'$ lying on the line $\la p_i,p\ra$ is $p_j$ and the image of
$p_j'$ is $p_i$. The transformation leaves invariant any line $l$ through
$p$  not equal to a side of the triangle and induces a non-trivial
involution on it. In particular, $Q$ and $Q'$ have the same tangent
lines  through $p$ and have 2 common points $p_i,p_j$. It is easy to see
that they must coincide. Hence the point $p_k$ goes to $p_k'$ for $k \ne
i,j$. Composing with a projective transformation we map $p_i$ to $p_j'$,
$p_j$ to $p_i'$ and do not change $p_k, k\ne i,j$.

Thus we have constructed a map from
$(\bbP(V(2))'\times \bbP(V(5))^\textup{s})/\SL(2)$ to $\calM_{\cub}\li$. Let
us show
that it inverts the map \eqref{bin}. Choose a line $m$ not intersecting
$l$ and 5 skew lines $l_1',\ldots,l_5'$ intersecting $l$ and $m$. Let
$l_1,\ldots,l_5$ be the lines such that $l,l_i,l_i'$ are coplanar. Then we
can blow down the skew lines $l_1,\ldots,l_5,m$ to the points
$p_1,\ldots,p_5,p$, respectively. The image of $l$ is a conic $Q$
through the points $p_1,\ldots,p_5$. The image of the line $l_i'$ is the
line $\la p_i,p\ra$. The pencil of planes through the line $l$
corresponds to the pencil of lines through $p$. Thus the pair of binary
forms $(F_5,F_2)$ defined by $(S,l)$ corresponds to the tangent
directions at $p$. Clearly, the pair $(S',l')$ reconstructed from
$(F_5,F_2)$ by the previous construction is isomorphic to $(S,l)$.

Now assume that $S$ has nodes. Let us interpret the first column  of
Table 1 as the type of the pair $(F_5,F_2)$ in the following way.
\begin{itemize}
\item[(1)] $\bft = (2222211)$ corresponds to the case when
$F_5,F_2$ have no multiple roots and no roots in common;
\item[(2)] $\bft = (322221)$ corresponds to the case when
$F_5,F_2$ have no multiple roots and one common root;
\item[(3)] $\bft = (332222)$ corresponds to the case when
$F_5,F_2$ have no multiple roots and two common roots;
\item[(4)] $\bft = (222222)$ corresponds to the case when
$F_5$ has no multiple roots, $F_2$ has a double root which
is not a root of $F_5$;
\item[(5)] $\bft = (422211)$ corresponds to the case when
$F_5$ has one double root, $F_2$ has two roots which are
not roots of $F_5$;
\item[(6)] $\bft = (43221)$ corresponds to the case when
$F_5$ has one double root, $F_2$ has two roots, one of
them is a root   of $F_5$;
\item[(7)] $\bft = (4332)$ corresponds to the case when
$F_5$ has one double root, $F_2$ has two roots which are
roots of $F_5$;
\item[(8)] $\bft = (42222)$ corresponds to the case when
$F_5$ has one double root, $F_2$ has a double root;
\item[(8*)] $\bft = (42222)$ corresponds to the case when
$F_2$ has a double root which is a simple root of $F_5$;
\item[(9)] $\bft = (44211)$ corresponds to the case when
$F_5$ has two double roots, $F_2$ has two roots which are
not roots of $F_5$;
\item[(10)] $\bft = (52221)$ corresponds to the case when
$F_5$ has one double root which is a simple root of $F_2$,
the other root of $F_2$ is not a root  of $F_5$;
\item[(11)] $\bft = (4431)$ corresponds to the case when
$F_5$ has two double roots and its simple root is a root of
$F_2$, the other root of $F_2$ is not a root  of $F_5$;
\item[(12)] $\bft = (5322)$ corresponds to the case when
$F_5$ has one double root which is a simple roots of $F_2$
and another root of $F_2$ is a simple root of $F_5$;
\item[(13)] $\bft = (4422)$ corresponds to the case when
$F_5$ has one double root, $F_2$ has a double root which
is a simple root of $F_5$;
\item[(13*)] $\bft = (4422)$ corresponds to the case when
$F_5$ has two double roots, $F_2$ has a double root which
is not a root of $F_5$;
\item[(14)] $\bft = (5421)$ corresponds to the case when
$F_5$ has two double roots, one of which is a simple root of
$F_2$, the other root of $F_2$  is not a root  of $F_5$;
\item[(15)] $\bft = (543)$ corresponds to the case when
$F_5$ has two double roots, $F_2$ has two roots in common
with $F_5$, one of them is a double root of $F_5$;
\item[(16)] $\bft = (444)$ corresponds to the case when
$F_5$ has two double roots and its simple root is a double
root of $F_2$;
\item[(17)] $\bft = (552)$ corresponds to the case when
$F_5$ has two double roots which are simple  roots of $F_2$.
\end{itemize}
We reconstruct the surface $(S,l)$ in each case as the blow-up
of a point $p$ and 5 points $p_1,\ldots,p_5$ lying on the
intersection of lines $L_1,\ldots,L_5$ with a conic $Q$
through $p_1,\ldots,p_5$ not containing $p$. We may take two points
$p_i$ and $p_j$ collinear with $p$. This is
indicated by saying that $L_i = L_j$. The points $p_i,p_j$ could be
infinitely near in which case we indicate that the
corresponding lines coincide and touch $Q$. The conic $Q$ could also
consist of two lines $H_1$ and $H_2$.

\begin{itemize}
\item [1)] A smooth conic $Q$ and five lines $L_1,...,L_5$ through $p$.
\item [2)] A line $L_i$ is tangent to $Q$.
\item [3)] Two lines $L_i, L_j$ are two tangent lines of $Q$.
\item [4)] $Q$ splits into two lines $H_1 + H_2$, but none of the $L_i$
contains the point
$H_1 \cap H_2$.
\item [5)] Two lines $L_i$ and $L_j$ coincide.
\item [6)] A line $L_i$ is tangent to $Q$ and $L_j = L_k$.
\item [7)] Two lines $L_i, L_j$ are two tangent lines of $Q$ and $L_k = L_l$.
\item [8)] $Q$ splits into two lines and two lines $L_i, L_j$ coincide.
\item [8*)] $Q$ splits into two lines $H_1 + H_2$ and a line $L_i$
contains the point $H_1 \cap H_2$.
\item [9)] $L_i = L_j$ and $L_k = L_l$.
\item [10)] $L_i = L_j$ and this line is tangent to $Q$.
\item [11)] A line $L_i$ is tangent to $Q$ and $L_j = L_k, L_l = L_m$.
\item [12)] $L_i = L_j$ and $L_k$ are two tangent lines of $Q$.
\item [13)] $Q$ splits into two lines $H_1 + H_2$, a line $L_i$ contains
the point $H_1 \cap H_2$ and $L_j = L_k$.
\item [13*)] $Q$ splits into two lines $H_1,H_2$, no line contains the
point $H_1\cap H_2$ and
$L_i = L_j$, $L_k = L_l$.
\item [14)] $L_i = L_j$ is tangent to $Q$ and $L_k = L_l$.
\item [15)] $L_i = L_j$ and $L_k$ are tangent lines of $Q$ and $L_l = L_m$.
\item [16)] $Q$ splits into two lines $H_1 + H_2$, $L_i$ contains the point
$H_1 \cap H_2$ and $L_j = L_k, L_l = L_m$.
\item [17)] $L_i = L_j$ and $L_k = L_l$ are two tangent lines of $Q$.
\end{itemize}
As in the nonsingular case, we show that this defines the inverse map.
\end{proof}

\subsection{}
Since the variety $(\bbP(V(2))\times \bbP(V(5)))/\SL(2)$ is obviously birationally
isomorphic to
the quotient $ \bbP(V(5))/\bbC^*$ (by fixing first a binary form of degree
2),
we obtain the following:

\subsection{Corollary.}
{\it The moduli space $\calM_{\cub}\li$ is isomorphic to an open subset of a
toric variety. In particular, it is rational. }

\subsection{Remark.} It is well-known that the moduli space of cubic
surfaces is rational. However, as
far as we know, the rationality of the space $\calM_{\cub}\li$ was not known.
Note also that the moduli space $\calM_{\cub}\li$ is birationally isomorphic
to the universal surface over the moduli
space of Del Pezzo surfaces of degree 4.

\section{The $K3$ surface associated to a cubic surface} \label{sec2}
\label{K3asscub}

\subsection{} In the previous section we associated a pair of binary
forms $(F_2,F_5)$ to a nodal cubic surface $S$ with a line $l$.
We now use these binary forms to define a $K3$ surface $X_{S,l}$.

We will show that $X_{S,l}$ depends only
on the nodal cubic $S$ and that the lines on a generic
$S$ correspond to certain `standard' elliptic fibrations
(cf.\ section \ref{standard}, Corollary \ref{27lines}).
Finally we relate $X_{S,l}$ to $S$ using a cubic fourfold.

\subsection{Definition}\label{intro1}
Let $S$ be a nodal cubic surface and let $l$ be a line on $S$.
Let $F_2(x_0,x_1)$ be a homogeneous form of degree 2
and let $F_5(x_0,x_1)$ be a homogeneous
form of degree 5 associated to $(S,l)$ as in \ref{f5f2}.

To the pair $(S,l)$ we associate a
surface $X_{S,l}$ which is a nonsingular
minimal model of the double plane with the branch divisor
\begin{equation}\label{bbranch}
W:\qquad x_2(F_2(x_0,x_1)x_2^3+F_5(x_0,x_1)) = 0.
\end{equation}

It is easy to check that the properties (i)-(iii) in \ref{f5f2}
are equivalent
to the property that any singular point of the curve $W$ is
analytically equivalent to a singularity $f(x,y) = 0$
such that the surface singularity $z^2+f(x,y) = 0$ is a double
rational point. This implies that $X_{S,l}$ is a $K3$ surface.
The multiplication of $x_2$ by a primitive cube root of
unity induces an automorphism of $X_{S,l}$ of order 3.

\subsection{Elliptic fibration}\label{ellfib}
Consider the pencil of lines
$$
L(t_0,t_1):t_1 x_0-t_0 x_1 = 0
$$
in $\bbP^2$ passing through the point
$(0,0,1)$. Since a general line
$L(\lambda,\mu)$ intersects $W$ at four nonsingular points,
we obtain that the pre-image of the pencil of lines on
$X_{S,l}$ is an elliptic pencil. Thus we have an elliptic fibration
$$
f=f_l:X_{S,l}\longrightarrow \bbP^1.
$$
The singular fibres correspond to lines
$L(t_0,t_1)$ such that  $F_5(t_0,t_1) = 0$ or
$F_2(t_0,t_1) = 0$.  The proper transform of $W$ in the blow-up $V \cong
\bfF_1$ of the point $(0,0,1)$ is a curve
$\bar{W}$ in the linear system $|6f+4e|$, where $e$ is the exceptional
section and $f$ is a fibre. The pre-image $T$ of
the line $x_2 = 0$ is a component of $\bar{W}$. It is a section with the
divisor class $f+e$. The pre-image of a line
corresponding to a zero $(x_0,x_1)$ of $F_5$ is a fibre of $V\to \bbP^1$
over $(x_0,x_1)$ which intersects $B = \bar{W}-T$
with multiplicity $3$ at a point where $B$ intersects $T$. A line
corresponding to a zero of $F_2$ is a fibre which
intersects $B$ with multiplicity 3 at a point where $B$ intersects $e$. The
surface $X_{S,l}$ is isomorphic to a minimal
resolution of the double cover of $V$ branched along $\bar{W}$.

Now it is easy to describe the singular fibres of the elliptic fibration
$f:X_{S,l}\to \bbP^1$. For example,
in the case when $F_5$ and $F_2$ have no multiple roots and have no common
roots, the fibres over the zeroes of $F_2$
are  cuspidal cubics. The fibres over the zeroes of $F_5$ are reducible of
type $IV$ in Kodaira's notation. If $F_2$ has
a common zero with $F_5$, the fibre of $V\to \bbP^1$ becomes an irreducible
component of $B$. The corresponding fibre of
$f$ is of type $I_0^*$. If $F_2$ has a double root which is not a root of
$F_5$, then $B$ acquires a cusp. Instead of
two irreducible fibres of $f$ we obtain one reducible fibre of type $IV$.
If $F_5$ has a double root which is not a root
of $F_2$, then $B$ acquires a cusp at the curve $T$. The corresponding
fibre of $f$ is of  type $IV^*$.   It is not
difficult to describe the fibres in all possible cases.
Their Kodaira types are given in  Table 1.
Note that the irreducible singular fibres correspond to
zeroes of $F_2$ which are not zeroes of $F_5$.
Observe also that the pre-image of $T$ in $X_{S,l}$ is a section $s$ of the
elliptic fibration. The pre-image of $e$  is
a bisection $b$. If $B$ acquires a cusp at the exceptional section $e$ or
has a fibre component, then $b$ splits in two
disjoint sections.

\subsection{}\label{3auto} Let $l$ be a line on a nodal cubic surface $S$,
and let $m$ be another line disjoint from $l$.
Consider the rational map $T:l\times m\, - \to S$ defined by taking the
third intersection point of the line through
the points $(p,q)\in l\times m$ with $S$.
We denote by $L$ and $M$ the irreducible curves in $l\times m$ which
map onto the lines $l$ and $m$ in $S$ respectively under $T$.

\subsection{Lemma}\label{RL}
{\it The rational map $T$ extends to an isomorphism from
the blow-up $Z$ of $l\times m$ along $L\cup M$,
which is a set of 5 points (including infinitely near points)
to a minimal resolution  $\tilde{S}$ of $S$.
The curves $L$ and $M$ have bi-degree (2,1) and (1,2) respectively.
}

\begin{proof} This is just a straightforward computation. Choose
coordinates on $\bbP^3$ such that
$m:x_0 = x_1  = 0$ and $l:x_2 = x_3 = 0$ so that the equation of $S$ is
given by
\begin{equation}\label{EF}
\sum_{i,j=0}^1A_{ij}(x_2,x_3)x_ix_j+2\sum_{i=0}^1B_i(x_2,x_3)x_i = 0,
\end{equation}
where $A_{ij},B_i$ are homogeneous forms of degree 1 and 2, respectively.
Let $p = (a_0,a_1,0,0)\in l$,
$q = (0,0,a_2,a_3)\in m$.
The line $l'$ spanned by $p,\,q$ has parametric equation $(x_0,x_1,x_2,x_3)
= (sa_0,sa_1$
$ta_2,ta_3)$.
Plugging it in equation (\ref{EF}), we obtain
\[
st\left(s\sum_{i,j=0}^1A_{ij}(a_2,a_3)a_ia_j+
2t\sum_{i=0}^1B_i(a_2,a_3)a_i\right) = 0.
\]
Thus the rational map $T$ is given by the formula
\begin{equation}\label{EF2}
T(p,q) = (Ma_0,Ma_1,La_2,La_3),
\end{equation}
where
\begin{equation}\label{EEE}
M(p,q) = -2\sum_{i=0}^1B_i(a_2,a_3)a_i, \quad L(p,q)
=\sum_{i,j=0}^1A_{ij}(a_2,a_3)a_ia_j.
\end{equation}
It is easy to see that the base locus $Z$ of the linear system of divisors
of bi-degree (3,3) defining $T$ is
the complete intersection of the divisor $M = 0$ of bi-degree (1,2) and
$L= 0$ of bi-degree (2,1).
Local computations show that $Z$  is reduced and consists of
5 points  if and only if $S$ is smooth.
The rational map $T$ is obviously birational,
and defines a birational morphism $T':S'\to S$
of the blow-up $S'$ of $l\times m$ along $Z$.
It is clear that the proper images under $T$ of the divisors
$L = 0$ and $M = 0$ are the lines $l$ and $m$, respectively.
Comparing the Betti numbers of $S'$ and $\tilde{S}$, we see that they are
equal.
Thus $T'$ defines an isomorphism from $S'$ to $\tilde{S}$.
\end{proof}

\subsection{Remark.}  Assume $S$ is nonsingular. Then we obtain that $S$ is
isomorphic to the blow-up of 5 distinct points in
$\bbP^1\times \bbP^1$. The map $S\to \bbP^1\times \bbP^1$ is the blowing
down of 5 disjoint lines intersecting the
lines $l$ and $m$. This is of course well-known. Take any two skew lines on
$S$. It is known that there are exactly
five skew lines on $S$ which intersect $l,m$. The easiest way to see it is
to complete $l,m$ to a set of six skew lines
$n_1 = l,\,n_2 = m,\,n_3,\ldots,n_6$, then consider the blow-down $\pi:S\to
\bbP^2$ of these lines to points
$p_1,\ldots,p_6$ in the plane. The five skew lines are the proper inverse
transforms of the line spanned by $p_1,\,p_2$
and the four conics $C_i$ passing through all $p_j$'s except $p_i$ with
$3\leq i\leq 6$.  Blowing down the five lines,
we obtain $\bbP^1\times \bbP^1$. The images of the lines $l,m$ are the
curves of bi-degree $(2,1)$ and $(1,2)$. The
blowing down morphism $S\to \bbP^1\times \bbP^1$ which inverts $T$ is the
Cartesian product of the linear projections from the lines $l$ and $m$.

\subsection{The surface $X_{S,l,m}$.}\label{triplecover}
The divisor $W'= L+M$ on
$l\times m = \bbP^1\times\bbP^1$ is of
bi-degree (3,3).
Let us consider the cyclic triple cover $Y\to l\times m$
branched along $W'$.
It has singular points over the singular locus of $W'$.
If $L$ intersects $M$ transversally,
$Y$ has 5 double rational points of type $A_2$.
Let $X_{S,l,m}$ be a nonsingular minimal model of $Y$.

\subsection{Lemma.}\label{ellipY}
{\it Let
$$
f=f_{l,m}:X_{S,l,m}\longrightarrow m \cong \bbP^1
$$
be the composition of the blow down map $X_{S,l,m}\to Y$,
the triple covering $Y\rightarrow l\times m$ and
the second projection $l\times m\to m$.
Then $f$ is an elliptic fibration with a section
whose Weierstrass form is given by:
\begin{equation}\label{we}
y^2+x^3+F_5(t_0,t_1)^2F_2(t_0,t_1) = 0,
\end{equation}
where  the binary forms $F_2(t_0,t_1)$ and $F_5(t_0,t_1)$ coincide with
the binary forms $F_2$ and $F_5$ associated to $(S,l)$
in section \ref{f5f2}.}

\begin{proof}
For any general point $(t_0,t_1)\in
\bbP^1$, the fibre of $f$ over this point is isomorphic
to a plane cubic curve with the equation
\begin{equation}\label{cubiccurve}
x_2^3+(B_0(t_0,t_1)x_0+B_1(t_0,t_1)x_1)(A_{00}(t_0,t_1)x_0^2+
2A_{01}(t_0,t_1)x_0x_1+A_{11}(t_0,t_1)x_1^2) = 0.
\end{equation}
The cubic curve has an obvious automorphism of order $3$ defined by
multiplying $x_2$ by the third roots of unity.
As is well-known such a cubic can be reduced by a projective transformation
to the Weierstrass form
$$
y^2t+x^3+bt^3 = 0.
$$
The coefficient $b$ is the value of a certain $\SL(3)$-invariant $T$ on the
space of homogeneous
polynomials of degree 3 in 3 variables.  Using the explicit formula for $T$
(see \cite{Salmon2}, p. 192),
a direct computation shows that
\begin{equation}\label{eqnew1}
b = F_5(t_0,t_1)^2F_2(t_0,t_1),
\end{equation}
where
\[
\begin{array}{rcl}
F_5& =&
B_0(t_0,t_1)^2A_{11}(t_0,t_1)+B_1(t_0,t_1)^2A_{00}(t_0,t_1)-2A_{01}(t_0,t_1)B_0(t_0,t_1)B_1(t_0,t_1),\\
F_2 &=& A_{00}(t_0,t_1)A_{11}(t_0,t_1)-A_{01}(t_0,t_1)^2.
\end{array}
\]
Let $t_1x_2-t_0x_3 = 0$ be the pencil of planes through the line
$l:x_2 = x_3 = 0$.
Using the equation \eqref{EF} of $S$ we find that the pencil of conics
defined by the line $l$ has the equation
\begin{equation}\label{EFF}
A_{00}(t_0,t_1)x_0^2+2A_{01}(t_0,t_1)x_0x_1+A_{11}(t_0,t_1)x_1^2+
2B_0(t_0,t_1)x_2x_0+2B_1(t_0,t_1)x_2x_1 = 0.
\end{equation}
Its discriminant is equal to
$$
\det
\begin{pmatrix}A_{00}&A_{01}&B_0\\
A_{01}&A_{11}&B_1\\
B_0&B_1&0\end{pmatrix} = -F_5(t_0,t_1).
$$
The restriction of the member of the pencil corresponding to the
parameters $(t_0,t_1)$ to the line $l$ is
given by the binary form
\begin{equation}\label{pencil}
A_{00}(t_0,t_1)x_0^2+2A_{01}(t_0,t_1)x_0x_1+A_{11}(t_0,t_1)x_1^2 = 0.
\end{equation}
The discriminant of this binary form is equal to
$$
\det
\begin{pmatrix}A_{00}&A_{01}\\
A_{01}&A_{11}\end{pmatrix} = F_2(t_0,t_1).
$$
If $l$ does not contain nodes, the equation \eqref{pencil} defines a
base-point free pencil of divisors of degree 2
on $l$, and we see that $F_2 = 0$ describes the locus of points in the
parameter space of the pencil of conics where the
bisection $l$  ramifies. If $l$ contains a node, we may assume that its
coordinates are $(1,0,0,0)$. Then $A_{11} = 0$
and we get a pencil of divisors of degree 1 on $l$ with one base point. The
discriminant is equal to $-A_{01}^2$ and
describes one point with multiplicity 1 corresponding to the singular point
of the bisection $B$ defined by $l$.
Finally, if $l$ contains two nodes,  we may assume that $A_{11} = A_{00} =
0$. Then the pencil \eqref{pencil} cuts out
the fixed divisor with equation $A_{01}(t_0,t_1)x_0x_1 = 0$.
It is equal to zero when $A_{01}(t_0,t_1) = 0$.
These points correspond to fibre components of the bisection $B$
of the conic bundle.
The discriminant is again $-A_{01}(t_0,t_1)^2$.
\end{proof}

\subsection{Theorem}\label{independence}
{\it Let $S$ be a nodal cubic surface and let
$l$ be a line on $S$.
Then the isomorphism class of the $K3$ surface $X_{S,l}$
associated to a pair $(S,l)$  is independent
on the choice of the line $l$.}

\begin{proof} We compare the elliptic fibration $f_l$
on $X_{S,l}$ obtained from the pencil of lines through the singular
point $(0,0,1)$ of the branch curve $W$ and the elliptic fibration
$f_{l,m}$ on the triple cover $X_{S,l,m}$, where $m$ is a line
disjoint from $l$.
The fibre of $f_l$ corresponding to a general line
$t_1x_0-t_0x_1 = 0$, with $t_0=1$,
passing through the point $(0,0,1)$ is birationally
isomorphic to the curve
$$
z^2+x_2x_0^2(F_2(1,t_1)x_2^3+F_5(1,t_1)x_0^3) = 0.
$$
After the  change of variables $y = F_5z/x_0x_2^2$, $x = F_5x_0/x_2$  we
reduce this equation to the Weierstrass
form \eqref{we} from Lemma \ref{ellipY}.
This shows that the surfaces $X_{S,l}$ and $X_{S,l,m}$ have isomorphic
elliptic pencils.
Hence $X_{S,l}\cong X_{S,l,m}$.
Switching the roles of $l$ and $m$, we see that $X_{S,l}\cong X_{S,m}$.
It is easy to see that if two lines $l,m$ on $S$ are not skew,
then there exists a third line $n$ which is disjoint from $l$ and $m$,
so again $X_{S,l}\cong X_{S,n}\cong X_{S,m}$.
We conclude that $X_{S,l}$ does not depend on a choice of a line $l$.
\end{proof}

\subsection{Definition.} Let $S$ be a nodal cubic surface.
A $K3$ surface associated  to $S$ is a $K3$ surface $X_S$
isomorphic to the surface $X_{S,l}$ associated to a pair
$(S,l)$, where $l$
is a line on $S$ defined in section
\ref{intro1} or the surface $X_{S,l,m}$ associated
to a triple $(S,l,m)$, where $l,m$ is a pair of skew lines on $S$
defined in \ref{triplecover}.

As a corollary of the results above and those of the previous section
we have:

\subsection{Corollary}\label{modell}
{\it The moduli space $\calM_{\ncub}\li$
is isomorphic to the moduli space of elliptic $K3$ surfaces with
the Weierstrass form
$$
y^2+x^3+F_5(t_0,t_1)^2F_2(t_0,t_1) = 0,
$$
where $(F_5,F_2)$ is a stable pair of binary forms of degrees 5 and 2.}

\subsection{Cubic fourfolds.} Let $Y \subset \bbP^5$ be a  cubic
hypersurface containing  a plane $\Pi$.
Then the projection from a plane defines a structure of a quadric bundle
on  the blow-up of $Y$ along the plane $\Pi$. The
discriminant curve of the quadric bundle is  of degree 6. The double cover
of $\bbP^2$ branched along the discriminant
parametrizes the rulings of quadrics. If $Y$ is general enough, a minimal
nonsingular model of the double cover  is a $K3$
surface. All of this is well-known and can be found for example in
\cite{Voi}. Let $F(x_0,x_1,x_2,x_3) = 0$ be a nodal
cubic surface. Consider the cubic fourfold $Y$ defined by the equation
\begin{equation}\label{cubic}
F(x_0,x_1,x_2,x_3)+x_4x_5(x_4+x_5) = 0.
\end{equation}
Let $l\subset \{x_4 = x_5 = 0\}$ be a line on $S$ and $\Pi$ the plane
spanned by $l$ and a point $(0,0,0,0,0,1)$.
Then $\Pi$ is contained in $Y$, and the direct computation shows that the
corresponding $K3$ surface is isomorphic to the
$K3$ surface $X_{S,l}$.

Let $m$ be a line on $S$ disjoint from $l$ and let $\Pi'$ be the plane
contained in $Y$ spanned by
$m\subset \{x_4 = x_5 = 0\}$ and the point $(0,0,0,0,1,0)$. Consider the
rational map $T:\Pi\times \Pi'\, - \to Y$
defined by taking the third intersection point of the line through the
points $(p,q)\in \Pi\times \Pi'$ with $Y$. A
straightforward computation shows that the fundamental locus $B$ of $T$ is
a complete intersection of divisors of
bi-degree $(2,1)$ and $(1,2)$ and its minimal nonsingular model is
isomorphic to the $K3$ surface $X_S$. The surface $B$
is nonsingular if $S$ is nonsingular.

The lattice of transcendental cycles of $X_S$ and that of the cubic
fourfold $Y$ are isomorphic.
In fact, the blow-up $Y'$ of $Y$ along the union of two disjoint
planes is isomorphic to the blow-up of $\bbP^2\times \bbP^2$ along the $K3$
surface $X\cong X_{l,m}$. This gives an
isomorphism of Hodge structures
$$
H^4(Y',\bbZ)\cong H^4(\bbP^2\times\bbP^2,\bbZ)\oplus H^2(X,\bbZ)[1]. $$
This isomorphism is compatible with the
cup-product such that the two summands become orthogonal.
Here $H^2(X,\bbZ)[1]$ is identified with
$\xi\cdot\pi^*(H^2(X,\bbZ))$, where
$\pi:Y'\to \bbP^2\times \bbP^2$ is the
natural morphism of the blow-up and $\xi$ is
a cohomology class from $H^2(Y',\bbZ)$ which cuts out the tautological
class of the exceptional divisor isomorphic to
the projectivization of the normal bundle of $X$.
This implies that the sublattice consisting of algebraic
cycles in $H^4(Y',\bbZ)$ is isomorphic to
$H^4(\bbP^2\times\bbP^2,\bbZ)\oplus \Pic(X)[1]$.
Passing to the orthogonal complements we get the result.

\subsection{Cubic threefolds.} \label{3fold} We relate the $K3$ surface
$X_S$ to the
Matsumoto-Terasoma curve associated to $(S,l)$. Given a smooth cubic
surface $S$ in $\bbP^3$, we define,
following \cite{ACT}, the cubic threefold $V\subset \bbP^4$ to be the
triple cover of $\bbP^3$ branched along $S$.
So if
$$
S:\;F(x_0,x_1,x_2,x_3) = 0,
$$
then
$$
V:\; F(x_0,x_1,x_2,x_3) + x_4^3 = 0.
$$
Note that $S\subset V$ (the points of $V$ with $x_4=0$), hence a line
$l\subset S$ defines a line, also denoted by $l$,
in $V$. The projection of a cubic threefold away from a line on $\bbP^2$
defines the structure of a conic bundle on the
blow-up of $V$ along the line. The associated discriminant curve in
$\bbP^2$ is a plane quintic.  A straightforward
computation shows that the discriminant curve is a plane quintic
with the equation
$$
W':\; F_5(t_0,t_1)+t_2^3F_2(t_0,t_1)=0,
$$
where the $F_i$ are as in \ref{f5f2},
so $W'$ is a component of $W$.

\subsection{Remark.}\label{tB'}
Each smooth point $t$ of the plane quintic $W'$ defines two lines (the
components of the singular conic in
the fibre of $V\rightarrow \bbP^2$ over $t$). Thus there is a natural
double cover $C'\rightarrow W'$. This double cover
was studied by Matsumoto and Terasoma in \cite{MT}, the corresponding
double cover $C\rightarrow \bar{W}'$ of the
normalizations of these curves is ramified in two points, which are
identified in $C'$. The curve $C$ is isomorphic to
the affine curve (\cite{MT}, (3.1)):
$$
v^3 -xf(x^2) = 0,
$$
where $f$ is a polynomial of degree 5.
The Prym variety of the double cover $C\rightarrow \bar{W}'$ is a
$5$-dimensional principally polarized abelian variety which
is isomorphic to the intermediate Jacobian variety $P$ of
the cubic threefold $V$ (cf.\ \cite{MT}).
The Matsumoto-Terasoma curve $C$ has the following property.

\subsection{Proposition.} \label{visot} {\it Let $f:X_{S,l}\rightarrow
\bbP^1$ be the elliptic fibration as in
the subsection \ref{ellfib}.
The pull-back of $X_{S,l}$ along the base change $C\to \bbP^1$,
$(v,x)\mapsto x$,
is birationally equivalent to the product $C\times E$ where $E$ is the
elliptic curve with $j=0$:
$E\cong\bbC/(\bbZ+\bbZ\zeta_3)$. }

\begin{proof}
In \cite{MT} it is proved that
that $W=C/\iota$ where $\iota$
is the (Clemens-Griffiths) involution $\iota:(v,x)\mapsto (-v,-x)$.
Therefore the quotient curve is given by
$y^3=u^2f(u)$ where $u=x^2$ and $y=xv$.
This curve is birationally equivalent to $W'$. In fact, choosing
coordinates such that $F_2(y_0,y_1)=y_0y_1$
the equation of $W'$ is $y_2^3y_0y_1+F_5(y_0,y_1)$, hence
$y_2^3y_1+F_5(1,y_1)$ is an affine equation. Putting
$v=-y_1y_2$, $u=y_1$ we find the birational isomorphism with $f(u)=F_5(1,u)$.

The function field of $X_{S,l}$ is defined by
$s^2=y_0y_1+F_5(y_0,y_1)$.
The elliptic fibration is given by the
rational function $t=y_1/y_0$. Rewriting the equation we get:
$(s/y_0)^2=t+y_0^3F_5(1,t)$, equivalently, since $F_5(1,t)=f(t)$:
$$ Y^2=X^3+tf(t)^2\qquad (X=y_0f(t),\; Y=sf(t)/y_0). $$ Since on $C$ we
have $v^6=tf(t)^2$ we can write this as
$(sf(t)/y_0v^3)^2=(y_0f(t)/v^2)^3+1$, which is the equation $Y^2=X^3+1$
of the curve $E$.
\end{proof}

\subsection{Remark.}
According to Donagi and Smith \cite{DS}, the Prym map
$\calR_6\to \calA_5$ has degree 27 with the Galois group $W(E_6)$.
Identifying the branch points on $W$ and the ramification points
on $C$, we obtain the admissible double cover
$C'\rightarrow W'$ in $\calR_6$.
Thus we get $27$ `natural' pre-images of $P$ under the Prym map.
However, the Prym map has 2-dimensional fibre over
the intermediate Jacobian of a cubic threefold,
in fact any line in the threefold defines
an admissible double cover in $\calR_6$.

\section{The Picard lattice}\label{piclat}

In this section we compute the Picard lattice $\Pic(X_S)\subset
H^2(X_S,\bbZ)$ of the $K3$ surface $X_S$
associated to a nodal cubic surface and its orthogonal complement,
the lattice of transcendental cycles $T_{X_S}:=\Pic(X_S)^\perp$.

\subsection{Lattices.}\label{lattices}

Recall the following two lattices:
$$
U = \left(\bbZ^2,\begin{pmatrix}0&1\\1&0\end{pmatrix}\right),
\qquad A_2=\left(\bbZ^2,\begin{pmatrix}-2&1\\1&-2\end{pmatrix}\right). $$

The second cohomology group $H^2(X,\bbZ)$ equipped with the quadratic form
defined by the cup-product is an
even unimodular lattice of signature (3,19). It is isomorphic to the $K3$
lattice
$$ L = U^{\oplus 3}\oplus E_8^{\oplus 2},
$$
where $E_8 = \bbZ^8$ with the quadratic form defined by the opposite of the
Cartan matrix of the root system of
type $E_8$. In general, $A_m,D_n, E_k$ denote the root lattices of the
simple root systems of the corresponding symbol
(with the Cartan matrix multiplied by $-1$).

For any lattice $M$ we denote by $M(n)$ the lattice $M$ with the quadratic
form multiplied by $n$. Let $M$ be a
nondegenerate even lattice.  The dual abelian group $M^*$ contains $M$ as a
subgroup of finite index, the quotient
group $D(M) = M^*/M$ is called the {\it discriminant group} of $M$. It is
equipped with a quadratic form
$$ q:D(M)\to
\bbQ/2\bbZ,\qquad q(m^*+M) = t^{-2}(tm^*,tm^*)+2\bbZ,
$$
where $t\in\bbZ$ is such that $tm^*\in M$. We use the notation $\O(M)$
(resp.\ $\O(D)$) to denote the group of
automorphisms of $M$ (resp. $D(M)$) preserving the quadratic form. If $M$
is a primitive sublattice of a unimodular
lattice there is a natural isomorphism $D(M)\cong D(M^\perp)$.

\subsection{Lattices $M(\bft)$ and $T(\bft)$}\label{mt} Recall that a
choice of a line on a nodal cubic surface
$S$ defines an elliptic pencil on $f : X_S \to \bbP^1$.
Its type is
determined by the type  vector $\bf t$ of the conic bundle on $S$
corresponding to $l$, cf.\ \ref{conics}.
We call it the type vector of $(S,l)$ and the type
vector of the   elliptic fibration.
We will explain later that for any possible type vector $\bft$
there exists a pair $(S,l)$ of type $\bft$ such that the Picard
lattice of the $K3$ surface $X_{S}$ is of rank $12+2r+2e$,
where  $r$ is the number of nodes on $S$ and $e$ is the number
of Eckardt points on $l$.
We denote by $M(\bft)$ the smallest primitive
sublattice of $H^2(X_S,\bbZ)$ containing the sections and components
of fibres of the elliptic fibration defined by the line $l$.
Note that $\Pic(X_S)\cong M(\bft)$.  We will
compute the lattice $M(\bft)$ and its orthogonal complement $T(\bft)$ in
$H^2(X_S,\bbZ)$.

\subsection{Proposition.}\label{ptlatticev} {\it Assume that the
Mordell-Weil group MW$(f)$ is finite.  Then the lattices
$M(\bft)$ and
$T(\bft)$  are given in the following Table 2}:
\begin{table}[h]
\[
\begin{array}{rllllll}
{}& \bft&  M(\bft) & T(\bft) \\
1) & (2 2 2 2 2 1 1)\qquad& U \oplus A_{2}^{\oplus 5} & A_2(-1) \oplus
A_2^{\oplus 4}\\
\noalign{\smallskip}
2) & (3 2 2 2 2 1) & U \oplus D_{4} \oplus A_{2}^{\oplus 4} & A_2(-2)
\oplus A_2^{\oplus 3} \\
\noalign{\smallskip}
3) & (3 3 2 2 2)  & U \oplus D_{4}^{\oplus 2} \oplus A_{2}^{\oplus 3} &
A_2(-1) \oplus A_2(2)^{\oplus 2} \\
    \noalign{\smallskip}
4) & (2 2 2 2 2 2) & U \oplus E_{6} \oplus A_{2}^{\oplus 3} & A_2(-1)
\oplus A_2^{\oplus 3}
\\
\noalign{\smallskip}
5) & (4 2 2 2 1 1)  & U \oplus E_{6} \oplus A_{2}^{\oplus 3} & A_2(-1)
\oplus A_2^{\oplus 3}\\
\noalign{\smallskip}
6) & (4 3 2 2 1)  & U \oplus D_{4} \oplus E_{6} \oplus A_{2}^{\oplus 2} &
A_2(-2) \oplus A_2^{\oplus 2} \\
\noalign{\smallskip}
7) & (4 3 3 2)  & U \oplus D_{4}^{\oplus 2} \oplus E_{6} \oplus A_{2}
\qquad&
A_2(-2) \oplus A_2(2)\\
\noalign{\smallskip}
8) & (4 2 2 2 2)  & U \oplus E_{6}^{\oplus 2} \oplus A_{2} & A_2(-1) \oplus
A_2^{\oplus 2}\\
\noalign{\smallskip}
9) & (4 4 2 1 1)  & U \oplus E_{6}^{\oplus 2} \oplus A_{2} & A_2(-1) \oplus
A_2^{\oplus 2} \\
\noalign{\smallskip}
10) & (5 2 2 2 1)  & U \oplus E_{8} \oplus A_{2}^{\oplus 3} & A_2(-1)
\oplus A_2^{\oplus 2} \\
\noalign{\smallskip}
11) & (4 4 3 1)  & U \oplus E_{6}^{\oplus 2} \oplus D_{4} & A_2(-2) \oplus
A_2 \\ \noalign{\smallskip}
12) & (5 3 2 2)  & U \oplus E_{8} \oplus D_{4} \oplus A_{2}^{\oplus 2} &
A_2(-2) \oplus A_2 \\

\noalign{\smallskip}
13) & (4 4 2 2) & U \oplus E_{8} \oplus E_{6} \oplus A_{2} & A_2(-1) \oplus
A_2 \\ \noalign{\smallskip}
14) & (5 4 2 1)  & U \oplus E_{8} \oplus E_{6} \oplus A_{2} & A_2(-1)
\oplus A_2 \\ \noalign{\smallskip}
15) & (5 4 3)& U \oplus E_{8} \oplus E_{6} \oplus D_{4} & A_2(-2) \\
\noalign{\smallskip}
16) & (4 4 4) & U \oplus E_{8}^{\oplus 2} \oplus A_{2} & A_2(-1) \\
\noalign{\smallskip}
17) & (5 5 2)  & U \oplus E_{8}^{\oplus 2} \oplus A_{2} & A_2(-1) \\
\end{array}
\]
\caption{The Picard lattices}
\end{table}

\begin{proof} We will consider only the first two cases.
Let $f:X_S\to\bbP^1$ be the elliptic fibration of type $\bft =
(2 2 2 2 2 1 1)$ with Picard lattice
$\Pic(X_S)\cong M(\bft)$.  It follows from \ref{ellfib} that it has 5
reducible fibres of type $IV$ and a section $s$
defined by the line $x_2 = 0$. We will use the Shioda-Tate formula
\cite{Shi}:
\begin{equation}\label{shioda}
(\#\textup{MW})^2\cdot D(M(\bft)) = d_1\ldots d_k,
\end{equation}
where $\textup{MW}$ is the Mordell-Weil group and $d_1,\ldots,d_k$ are the
discriminants of the lattices generated
by components of reducible fibres not intersecting the zero section.  It
follows from \eqref{shioda}  that the
Mordell-Weil group $\textup{MW}$ is a torsion group of order $3^l$. We
claim that it  is trivial. Assume $\textup{MW}$
is not trivial. Then the translation by a nontrivial section defines an
automorphism of $X_S$ of order 3 which acts
trivially on $\Pic(X_S)^\perp$ and has at least 7 isolated fixed points on
$X_S$ (= the singular points of the 7
singular fibres). This contradicts the fact that the fixed locus of any
symplectic automorphism of $K3$ surface of order
3 is exactly 6 isolated points (Nikulin \cite{N3}, \S 5). Thus $f$ has a
unique section $s$. Now we use \eqref{shioda}
again and find that the discriminant of $M$ is equal to $3^5$. Since $M =
M(\bft)$ obviously contains the sublattice $U
\oplus A_2^{\oplus 5}$ of the same rank and discriminant (it is spanned by
the class of a fibre, the section, and
irreducible components of reducible fibres), it must coincide with it. The
discriminant group is then easy to compute.
Let $q_T$ be the discriminant form of $T$, then $q_T = -q_M$ (\cite{N1},
Prop.\ 1.6.1). We can easily see that $T$ and
$A_2(-1) \oplus A_2^4$ have the same discriminant form. It now follows from
Nikulin \cite{N1}, Cor.\ 1.13.3 that $T
\cong A_2(-1) \oplus A_2^4$.

Assume that the fibration is of type $(3 2 2 2 2 1)$. The product
$d_1\ldots d_k$ is equal to $2^23^4$.
The Shioda-Tate formula gives that either $\#\textrm{MW} = 1, 3$, or $
3^2$, or $6$.
Assume that $\textrm{MW}$ contains a nontrivial section $s'$. If $s'$ is of
order 2, it must leave invariant one
component in each fibre of type $IV$ and has one fixed point in it
different from the singular point of the fibre.
Altogether this gives 8 fixed points on reducible fibres and 1 fixed point
on the irreducible fibre. However, a
symplectic involution on a $K3$ surface has exactly 8 fixed points. Assume
$s'$ is of  order is 3. Then it leaves the
multiple component and non-multiple component of the fibre of type $I_0^*$
invariant. This gives $\ge 3$ isolated fixed
points on this fibre. Together with singular points of other fibres we get
$\ge 8$ isolated fixed points. This is
impossible.   So, the Shioda-Tate formula tells us that $D(M(\bft))$ is of
order $2^23^4$. The remaining arguments are
similar to the previous case.
\end{proof}

\subsection{The lattices $M$, $T$}\label{mt2} We set
$$
M:=U \oplus A_2^{\oplus 5} ,\qquad
T:= A_2(-1) \oplus A_2^{\oplus 4}.
$$
Since their discriminant groups are isomorphic and the quadratic forms are
the negative of each other,
they are orthogonal complements of each other in the unimodular lattice $L$
(see \cite{N1}). We set
\[D = D(M)\cong D(T).\]
These lattices correspond to the type $\bft = (2 2 2 2 2 1 1)$.

\subsection{An automorphism of order 3}\label{cauto}
As in section \ref{triplecover}, we choose two skew
lines on a nodal cubic surface $S$ and
consider the associated $K3$ surface $X =X_S\cong X_{S,l,m}$.
Recall that it is obtained as a minimal resolution of the triple
cyclic cover $Y$ of $\bbP^1\times \bbP^1$ branched along the union of two
divisors $L$ and $M$ of bidegree $(1,2)$ and
$(2,1)$. It is easy to describe the set of fixed points of the automorphism
$\sigma$ of $X$ defined by the triple cover.
We do it only in the case when $S$ is a nonsingular surface. Let
$q_1,\ldots,q_5$ be the intersection points of $L$ and
$M$. The cubic surface $S$ is obtained by blowing up the points $q_i$'s.
The surface $S$ is nonsingular if and and only
if no two points lie on a ruling, and no four points lie on a plane
section. An Eckardt point    on the line $l$
corresponds to a ruling which is tangent to $L$ at some point
$q_i$.

Assume that there are no Eckardt points on $l$. Consider the elliptic
fibration on $f:X\to \bbP^1$ corresponding to
the projection $\bbP^1\times \bbP^1\to \bbP^1$ such that $L$ is a section.
Its reducible singular fibres corresponds to
the ruling passing through the points $q_i$. Each fibre is of type $IV$.
Two components are the exceptional curves of
the resolution $X\to Y$ of a singular point of type $A_2$. The third
component is the proper transform of the ruling
passing through the corresponding point $q_i$. The bisection $b$ intersects
the latter component and one of the first
two components. The section $s$ intersects the other component coming from
the resolution of singularities.  The set of
fixed points of $\sigma$ is equal to the union of the section $s$, the
bisection $b$ and the singular points of the
reducible fibres.

In the case when $l$ contains one Eckardt point, the elliptic fibration
acquires one reducible  fibre of type $I_0^*$.
Other reducible fibres are of type $IV$. The bisection $b$ intersects the
multiple component $E_0$ of this fibre. The
section $s$ intersects a reduced component $E_1$. The fixed points of the
involution $\sigma$ is the union of the
section $s$, the bisection $b$, the point $E_0\cap E_1$, and the singular
points of fibres of type $IV$.  If $l$ has two
Eckardt points, we have two reducible fibres of type $IV$ and the set of
fixed points is described similarly to the
previous case.

\subsection{The involution $\tau$.}\label{P1'} Let $f:X \to \bbP^1$ be the
elliptic fibration with a section $s$ as in section \ref{cauto}. Let $\tau$
be the involution of $X$ defined by
the inversion $x\mapsto -x$ of each fibre.
Then $\tau$ switches the two components
of each singular fibre of type $IV$ which do not meet $s$ and preserves
each component of any singular fibre of
type $I_{0}^*$.

If $f$ has five singular fibres of type $IV$ and two singular fibres of
type $II$, then the fixed locus of $\tau$
is the union of $s$ and a smooth curve $C$ of genus $5$ which passes
through the singular point of each singular fibre.
If $f$ has four singular fibres of type $IV$, one of type $I_{0}^{*}$ and
one of type $II$, then the fixed locus of
$\tau$ is the union of $s$, the multiple component of the fibre of type
$I_{0}^{*}$ and a smooth curve of genus $3$. If
$f$ has three singular fibres of type $IV$ and two fibres of type
$I_{0}^{*}$, then the fixed locus of $\tau$ is the
union of $s$, two multiple components of singular fibres of type
$I_{0}^{*}$ and a smooth elliptic curve.

\subsection{Remark.}\label{infiniteAut}
The automorphism group of the $K3$ surface $X$ is infinite. For example,
consider the divisor consisting of the
$2$-section and the two components of a reducible singular fibre of $f$ not
meeting the section. It defines an elliptic
fibration on $X$ with a section which has two reducible singular fibres,
one is of type $I_3$ and another type $I_0^*$.
This elliptic fibration has a Mordell-Weil group of rank 4. Considering
translations by the sections of infinite order
we see that $\Aut(X)$ is an infinite group.

\subsection{Lemma}\label{sigmainv}
{\it Assume $S$ is nonsingular. Then
$$
H^{2}(X,\bbZ)^{\sigma^*}\subset \Pic(X),\qquad H^2(X,\bbZ)^{\sigma^*}
\cong M.
$$
The automorphism $\sigma$ acts trivially on the discriminant lattice
$D(H^2(X,\bbZ)^{\sigma^*})\cong D(M)$. }

\begin{proof}
Consider the elliptic fibration on $X$ defined in \ref{ellfib}. From
\ref{cauto} we know the description of
fixed points of $\sigma$. Assume first that all reducible fibres are of
type $IV$. Let $P$ be the sublattice of
$\Pic(X)$ spanned by the divisor classes of a fibre, of the section $s$ and
of the irreducible components of
fibres which do not intersect $s$. It is immediate that $P\cong M$ and
$\sigma$ acts identically on $P$. The fixed locus
$X^\sigma$ of the automorphism $\sigma$ consists of 5 isolated fixed points
(the singular points of the reducible
fibres) and two smooth rational curves (the section $s$ and the bisection
$b$). Applying the Lefschetz fixed point
formula we obtain that the trace of $\sigma^*$ on $H^2(X,\bbZ)$ is equal to
7. Thus the trace of $\sigma^*$ on $P^\perp$
is equal to $7-12 = -5$. This easily implies that the characteristic
polynomial of $\sigma^*$ on $P^\perp\otimes \bbC$
is equal to
$(t^2+t+1)^5$. Therefore $P^{\perp}\otimes \bbC$ does not contain non-zero
$\sigma^*$-invariant elements,
so $H^2(X,\bbZ)^{\sigma^*} = P\cong M$. Since $\sigma^*$ acts trivially on
$P \cong M$, it also acts trivially on
$D(P)\cong D(M)$.

Suppose now that $f$ contains a fibre $F = 2E_0+E_1+E_2+E_3+E_4$ of type
$I_0^*$. Assume that $E_1$ intersects the
section $s$. Then the divisor classes $E_0+E_2+E_3+E_4$ and $E_0$ are
$\sigma$-invariant and span a lattice of type
$A_2$. We define the lattice $P$ similar to the above by using this
contribution from a fibre of type $I_0^*$. The
remaining arguments are the same.
    \end{proof}

\section{The moduli space of $K3$ surfaces associated to a cubic surface}
\label{modulik3}

\subsection{} We first recall the basic facts about moduli of $K3$
surfaces. In the subsections before \ref{fixvm},
$M$ will be any even non-degenerate sublattice of signature $(1,t)$.

\subsection{Markings.}\label{defmark}
We recall the definition of a
$M$-marking of a $K3$ surface $X$ (see \cite{Do}). Fix a connected
component $V(M)^+$ of the cone
$V(M) = \{x\in M\otimes \bbR:(x,x) > 0\}$ and a subset $\Delta(M)^+$
of the set
$\Delta(M) = \{\delta\in M:(\delta,\delta) = -2\}$ such that
\begin{itemize} \item{} $\Delta(M)= \Delta(M)^+\coprod \Delta(M)^-$, where
$\Delta(M)^- =\{-\delta : \delta\in
\Delta(M)^+\}$, \item{} any $\delta\in\Delta(M)$ which can be written as a
nonnegative linear combination of elements
from $\Delta(M)^+$ belongs to $\Delta(M)^+$.
\end{itemize}
With these
choices, we define the subset \[ C(M)^+ = \{h\in
V(M)^+:\;(h,\delta) > 0\quad \text{for all}\quad \delta\in \Delta(M)^+\} \]
and we define $C(M)$ to be the closure of
$C(M)^+$ in $M\otimes\bbR$.

Now we define a $M$-marking of $X$ as a primitive lattice embedding
$\phi:M\to \Pic(X)$ such that
$C(X)^+ \cap \phi(M\otimes {\bbR}) \subset \phi(C(M)^+)$, where $C(X)^+$ is
the cone in
$\Pic(X)\otimes \bbR$ spanned by the pseudo-ample (i.e.\ nef and big)
divisor classes of $X$.

Note that the closure of $C(X)^+$ is the nef cone $C(X)$.
The closure $C(M)$ of $C(M)^+$ is the subset of the closure of $V(M)^+$
which consists of vectors $v$ such that
$(v,\delta)\ge 0$ for any $\delta\in \Delta(M)^+$. The marking $\phi$
embeds $C(X) \cap \phi(M\otimes {\bbR})$
in $\phi(C(M))$.  For any
$\delta\in \Delta(M)^+$ the image $\phi(\delta)$ is a divisor class $R$
with $R^2 = -2$. For any $v\in
C(M)$ the image $\phi(v)$ is a pseudo-ample divisor $D$ with $D^2 \ge
0$.  Since $R\cdot D = (\delta,v) > 0$, it
follows from Riemann-Roch that $R$ is effective. Note that $R$ is not
necessary the divisor class of an irreducible
curve (a $(-2)$-curve).

The marking is called \emph{ample} if $\phi(C(M)^+)\cap \Pic(X)^{+} \ne
\emptyset$, where $\Pic(X)^{+}$ is
the ample cone of $X$. It is easy to see that a marking $\phi$ is ample if
and only if the orthogonal complement of
$\phi(M)$ in $\Pic(X)$ does not contain the divisor classes of
$(-2)$-curves. In particular, any marking with $\phi(M) =
\Pic(X)$ is ample.

A pair $(X,\phi)$, where $\phi$ is a $M$-marking (resp. an ample
$M$-marking), is called a {\it $M$-polarized
$K3$ surface} (resp. {\it ample $M$-polarized $K3$ surface}). Two
$M$-polarized $K3$ surfaces $(X,\phi)$ and
$(X',\phi')$
are called isomorphic if there exists an isomorphism $f:X\to X'$ such that
$\phi = f^*\circ \phi'$.

\subsection{Moduli of $M$-polarized surfaces.}\label{modpol} It is known
(see \cite{Do}) that there exists a coarse
moduli space $\calM_{K3,M}$ of isomorphism classes of $M$-polarized $K3$
surfaces. Let us assume that $M$ admits a
unique (up to an isometry) embedding into the $K3$ lattice $L = U^{\oplus
3}\oplus E_8^{\oplus 2}$. Fix such an
embedding. Let $T$ be the orthogonal complement of $M$ in $L$. Any
$M$-marking $\phi$ of a $K3$ surface $X$ extends to
an isometry $\tilde{\phi}:L\rightarrow H^2(X,\bbZ)$ (a cohomology marking
of $X$). Extending $\tilde{\phi}$
$\bbC$-linearly, we get a one dimensional subspace
$\tilde{\phi}^{-1}(H^{2,0}(X))\subset T\otimes\bbC$ which is called
the \emph{period} of $(X,\tilde{\phi})$. $$ \begin{array}{rcl} M&\subset&L\\
\phi\downarrow&&\downarrow\tilde{\phi}\\
\Pic(X)&\hookrightarrow&H^2(X,\bbZ) \end{array} $$

The moduli space
$\calM_{K3,M}$ is isomorphic to the quotient ${\calD}_M/\Gamma_M$, where
${\calD}_M$ is the union of two copies of
a Hermitian symmetric domain of type IV corresponding to the inner product
vector space $T\otimes \bbR$ of signature
$(2,20-t)$, ${\calD}_M$ is a subset of the projective space
$\bbP(T\otimes\bbC)$. The group $\Gamma_M$ is the subgroup
of the orthogonal group $\O(L)$ of $L$ which leaves $M$ pointwise fixed. It
is also isomorphic to the subgroup of
$\O(T)$ which acts identically on the discriminant group $D(T) = T^*/T$.

The isomorphism classes of ample $M$-polarized $K3$ surfaces are
parametrized by an open subset of $\calM_{K3,M}$
whose complement is the image in $\calM_{K3,M}$ of the union of
hypersurfaces in ${\calD}_M$ defined by lines in
$T\otimes \bbC$ orthogonal to vectors $r\in T$ with $r^2 = -2$.

\subsection{The group $W(M)$.}\label{elfi} For any $\delta \in \Delta(M)$
we can define a reflection
$s_\delta\in \O(M)$ associated to $\delta$ by $s_\delta : v \mapsto v +
(v,\delta)\delta$. Let $W(M)$ be the subgroup of
$\O(M)$ generated by all $s_\delta$'s. The set $C(M)$ is a fundamental
domain for $W(M)$ in the closure of $V(M)^+$.
Thus for any $v\in M$ with $v^2 \ge 0$ there exists a $w\in W(M)$ such that
$(w(v),\delta)\ge 0, $ for any $\delta\in
\Delta(M)^+$.

Let $(X,\phi)$ be a $M$-polarized $K3$ surface. Then for any $v\in M$ with
$v^2\ge 0$ there is a $w\in W(M)$
such that $\phi(w(v))\in C(M)$. In particular, for any given embedding
$\phi : M \rightarrow \Pic(X)$, there is a $w\in
W(M)$ such that
$C(X)^+ \cap \phi(M\otimes {\bbR}) \subset (\phi\circ w) (C(M)^+)$, i.e.,
$\phi\circ w$ is a $M$-marking.

\subsection{Fixing $V(M)^+$ and $\Delta(M)^+$}\label{fixvm} The lattice $M$
from \ref{mt2} has a unique
(up to an isometry) primitive embedding in the $K3$ lattice $L$ \cite{N1}
and we identify $M$ with a primitive
sublattice of $L$ from now on. We fix a basis in $U$ formed by two
isotropic vectors
$f_1,\,f_2$ with $(f_1,f_2) = 1$ and a simple root basis $r_1, r_2$ in
$A_2$, i.e., $(r_1)^2 = (r_2)^2 = -2$
with $(r_1,r_2) = 1$. We define a basis of $M$ by taking $f_1$, $f_2$ in
$U$ and $r_1, r_2$ in each copy of $A_2$.

We define $V(M)^+$ by requiring that $f_1+f_2\in V(M)^+$. We define
$\Delta(M)^+$ as follows. Firstly,
$(-2)$-vectors with $(f_1+f_2, v) > 0$ belong to it. Secondly, if
$(f_1+f_2, v) = 0$, then $v\in \Delta(M)^+$ if and
only if it is a nonnegative combination of $f_2-f_1$ and the $r_i$'s in
each copy of $A_2$.

\subsection{Automorphisms of $L$.}\label{markfib} Let $\rho_o$ be the
isometry of $A_2$ defined by
\[
\rho_o (r_1) = r_2, \quad \rho_{o} (r_2) = -r_1 - r_2. \] Obviously
$\rho_o$ is of order 3, has no non-zero
fixed vectors and acts trivially on $D(A_2) = (A_2)^{*}/A_2$.
Let $\rho$ be the isometry of $T= A_2(-1) \oplus A_2^{\oplus 4}$ defined by
$\rho = (\rho_o)^{\oplus 5}$.
Then $\rho$ is of order 3, has no non-zero fixed vectors and acts trivially
on $D(T)$. Thus the isometry $(1_M, \rho)$
of $M \oplus T$ can be extended to the one of the $K3$ lattice $L$ (Nikulin
\cite{N1}, Corollary 1.5.2). For simplicity
we denote this isometry of $L$ by the same letter $\rho$.

\subsection{Period domains.}\label{The period domain.} The period domain
for $M$-polarized $K3$ surfaces is \[
\calD_M = \{\omega\in \bbP(T\otimes_{\bbZ}\bbC):(\omega, \omega) = 0,\quad
(\omega, \bar{\omega}) > 0\}.
\]
Let $\rho$ be the isometry of $T$ defined in \ref{markfib}. Let
%\begin{equation}\label{V+}
\[
T\otimes \bbC = V_+\oplus V_-
%\end{equation}
\]
be the decomposition of $T\otimes \bbC$ into the two 5-dimensional
eigenspaces of $\rho$ with eigenvalues
$\zeta_3 = e^{2\pi i/3}$ and $\zeta_3^{-1}$, respectively. Since
\[ (\omega,\omega) = (\rho(\omega),\rho(\omega)) = \zeta^2(\omega,\omega),
\] we see that $(\omega,\omega) = 0$ for all
$\omega\in V_{+}$, and similarly for $V_-$. Let
\[
\calB
= \{ \omega \in \bbP(V_{+}) : ( \omega, {\bar \omega}) > 0 \}
\;=\;\calD_M\cap \bbP(V_+).
\]
In a suitable basis of $V_+$ we have
$(\omega,\bar{\omega})=x_0\bar{x}_0-(x_1\bar{x}_1+\ldots+x_4\bar{x}_4)$.
Thus, if $(\omega,\bar{\omega})>0$, then
$x_0\neq 0$ and we can normalize $x_0=1$, hence
$\calB$ is a 4-dimensional complex ball:
\[ \calB\cong\,\{x = (x_1,\ldots,x_4)\in\bbC^4:\sum_{i}x_i\bar{x}_i < 1\}.
\] The 4-ball is a bounded
symmetric domain of type $I_{1,4}$.

\subsection{Discrete groups.}
We define the following four groups using the notation from \ref{markfib}:
{\renewcommand{\arraystretch}{1.5}
\[
\begin{array}{rcl}
\Gamma_M&=&\{g\in \O(L):\;g(m)=m,\quad\forall m\in M \},\\
\tilde{\Gamma}_{\rho} &=& \{g\in \O(L) : g\circ \rho =
\rho\circ g\},\\ \Gamma_{\rho} &=& \{g\in \O(T) : g\circ \rho = \rho\circ
g\},\\ \Gamma_{M,\rho} &=& \Ker( \Gamma_{\rho}
\to \O(D)). \end{array} \] }

\subsection{The Hermitian module.}\label{af} The isometry $\rho$ of $T$
gives $T$ the structure of a free module
$\Lambda$ of rank $5$ over the ring of Eisenstein integers $\bbZ[\zeta_3]$:
for any $a+b\zeta_3\in \bbZ[\zeta_3]$ and
any $x\in T$ we have \[(a+b\zeta_3)\cdot x = (a1_T +b\rho)(x).\] If
$r_i,\,r_i'$ is the simple root basis of the $i$-th
copy of $A_2$ with $\rho(r_i) = r_i'$, then $\zeta_3r_i=r_i'$ and any
element in this $A_2$ can be written as $r=zr_i$
with $z=a+b\zeta_3\in \bbZ[\zeta_3]$. Note that: $$
z\bar{z}=(a+b\zeta_3)(a+b\zeta_3^{-1})=a^2-ab+b^2=-(r,r)/2. $$
Therefore the quadratic form on $T$ is twice the real part of the
$\bbZ[\zeta_3]$-valued Hermitian form $H$, of
signature $(1,4)$, on the Eisenstein lattice $T$ with
$$
H(z,w)=z_0\bar{w}_0-(z_1\bar{w}_1+\ldots+z_4\bar{w}_4). $$ The group
$\Gamma_{\rho}$ is the unitary group $U(T)$ of
$T$ considered as a Hermitian lattice
over the ring of Eisenstein integers
(see \cite{ACT}, \cite{AF}).

\subsection{The discriminant group.}\label{discgrp}
The residue field $\bbZ[\zeta_3]/\sqrt{-3}\bbZ[\zeta_3]$
is isomorphic to $\bbF_3$ and $\zeta_3 $ maps to $ 1 \,\text{mod}\ 3$.
Thus $V =
\Lambda/\sqrt{-3}\Lambda$ acquires a natural structure
of a 5-dimensional vector space over $\bbF_3$ equipped with a
non-degenerate quadratic form.
We show that the discriminant group $D(T)$ is isomorphic to $V$.
Define a $\bbZ$-linear homomorphism
\begin{equation}
\label{fh}
h:\Lambda\longrightarrow T^*,\quad
h(x) = (x+2\rho(x))/3,
\end{equation}
where we
identify $\Lambda$ with $T$ as a $\bbZ$-module.
Then $h(\sqrt{-3}x) = h((1+2\zeta_3)x) = (1+2\rho)^2x/3 = -x\in T$.
This shows that $h$ factors through an isomorphism
\[
V = \Lambda/\sqrt{-3}\Lambda \longrightarrow D(T) = T^*/T.
\] The basis $(r_1,\ldots,r_5)$ of $\Lambda$
(as $\bbZ[\zeta_3]$-module) is an orthonormal basis with respect
to $H$.
Since $h(r_i)^2 = (r_i+2r_i')^2/9 =
-\frac{2}{3}, (h(r_i),h(r_j)) = 0, i\ne j$,
we obtain that
\[
h(x)^2 = -\frac{2}{3}x^2.
\]
In particular, if we identify $D(T)$ with $V$,
then the quadratic form on $D(T)$ obtained from the quadratic form
on $V$ by multiplying it by $-\frac{2}{3}$.

If $Q$ is the root lattice of type $E_6$, then  $Q/3Q$ inherits a
non-degenerate quadratic form such that $Q$ is isomorphic to $V$ as quadratic
spaces over $\bbF_3$. This defines an isomorphism of groups

\begin{eqnarray}
W(E_6) &\cong &\textup{SO}(V),\\ \notag
\O(D(T)) &\cong &\O(V)  \cong \{1,-1\}\times \textup{SO}(V).\notag
\end{eqnarray}
All of this is well-known and can be found, for example, in \cite{Bourbaki},
Chapter 6, \S 4, exercise 2.

\subsection{Proposition}\label{surjection} {\it Each of the natural maps
$$ \tilde{\Gamma}_{\rho} \longrightarrow \Gamma_{\rho} \longrightarrow
\O(D(T)) $$ is surjective. In particular,
\[
\Gamma_\rho/\Gamma_{M,\rho} \cong \O(D(T))\cong\{\pm 1\}\times W(E_6). \]
Moreover, any isometry in
$\Gamma_{M,\rho}$ can be extended to an isometry of $L$ which acts
trivially on $M$ defining an injective homomorphism
of groups \[
\Gamma_{M,\rho}\hookrightarrow \Gamma_{M}. \] }

\begin{proof}
For the surjectivity of the map $\Gamma_{\rho} \to \O(D(T))$ see \cite{ACT},
Lemma 4.5.
It is proven in Nikulin \cite{N1}, Theorem 1.14.2 that the natural map
$\O(M) \to \O(D(M))$ is surjective.
By Corollary 1.5.2 of loc.\ cit.\
this implies that the map $\tilde{\Gamma}_{\rho} \to \Gamma_{\rho}$ is
surjective. The inclusion
$\Gamma_{M,\rho}\rightarrow \Gamma_{M}$ follows from (Nikulin \cite{N1},
Corollary 1.5.2). \end{proof}

\subsection{Definition.}
An (ample) $(M,\rho)$-polarized $K3$ surface is an (ample) $M$-polarized
$K3$ surface $(X,\phi)$ such that
there is an extension $\tilde{\phi}:L\stackrel{}{\rightarrow}H^2(X,\bbZ)$
of $\phi$ which satisfies $$
\tilde{\phi}^{-1}(H^{2,0}(X))\;\in\;\calB\quad(\subset \bbP(T\otimes\bbC)).
$$ Two $(M,\rho)$-polarized
$K3$ surfaces $(X,\phi)$ and $(X',\phi')$ are said to be isomorphic if
there is an isomorphism $f:X\rightarrow X'$ such
that $\phi=f^*\circ \phi'$ and
$\tilde{\phi}^{-1}\circ f^*\circ \tilde{\phi}'\in \O(L)$ commutes with
$\rho\in \O(L)$.

\subsection{Lemma.}\label{k3autos}{\it
Let $(X,\phi)$ be an ample $(M,\rho)$-polarized $K3$ surface. Then $X$ has
an automorphism $\sigma$ of
order 3 such that $\sigma^*=\tilde{ \phi}\circ\rho\circ\tilde{\phi}^{-1}$
for an extension
$\tilde{\phi}:L\rightarrow H^2(X,\bbZ)$ of $\phi$. In particular, $\sigma$
acts trivially on $\phi(M)$
$( \subset Pic(X) )$.}

\begin{proof}
Choosing $\tilde{\phi}$ as in the definition of $(M,\rho)$-polarization,
the period of $X$ is fixed by $\rho$.
Since $(X,\phi)$ is amply polarized, $\Pic(X) \cap M^{\perp}$ contains no
$(-2)$-vectors. Moreover, the
$M$-polarization of $X$ is ample and $\rho$ acts trivially on $M$.
Therefore \cite{Na}, Theorem 3.10 shows that $X$ has
an automorphism $\sigma$ with $\sigma^*=\tilde{
\phi}\circ\rho\circ\tilde{\phi}^{-1}$. \end{proof}

\subsection{The moduli spaces $\calK3_{M,\rho}^m$ and
$\calK3_{M,\rho}$}\label{modulispaces} We know from
section \ref{modpol} that the moduli space of $M$-polarized $K3$ surfaces
is isomorphic to $\calD/\Gamma_M$. The
isometry $\rho$ acts naturally on $T_\bbC$ as is described in \ref{The
period domain.} and induces an automorphism of
order $3$ of the domain $\calD_M\subset \bbP(T_\bbC)$. It defines the union
of two balls $\calB_\pm = \calD_M\cap
\bbP(V_\pm)$. Complex conjugation switches the two balls $\calB_\pm$.
Obviously the group $\Gamma_{\rho}$ is the
stabilizer subgroup of $\calB = \calB_+$ in $\Gamma_M$. We set $$
\calK3_{M,\rho}^{m} =\calB/\Gamma_{M,\rho},\qquad
\calK3_{M,\rho} = \calB/\Gamma_{\rho}.
$$

The element $-I\in \Gamma_\rho$ acts trivially on $\bbP(T\otimes\bbC)$ and
thus on $\calB$, and $-I$ maps to
$-1\in \O(D)$. Thus $\O(D)/\{\pm 1\}\cong W(E_6)$
acts on $\calK3_{M,\rho}^m$ and
there is a natural map:
$$
\pi_M:
\calK3_{M,\rho}^{m}\longrightarrow \calK3_{M,\rho}\;\cong
\calK3_{M,\rho}^{m}/W(E_6). $$
For $r\in L$, let $r^\perp$ be the hyperplane in $\bbP(V_+)$ of lines
orthogonal to $r$, and let $H(r)$ be
its intersection with $\calB$. The {\it discriminant locus} is the subset
$\calH\subset\calB$ defined by:
$$ \calH = \bigcup_{r} H(r),
$$
where $r$ varies over the set of all $(-2)$-vectors in $T=M^\perp$. The
image of $\calH$ in $\calK3_{M,\rho}^m$
(resp. $\calK3_{M,\rho}$) will be denoted by $\Delta^m$ (resp. $\Delta$).

It follows from Lemma \ref{k3autos} that the quasi-projective variety
$\calK3_{M,\rho}^m\setminus\Delta^m$ is
the coarse moduli space of ample $(M,\rho)$-polarized $K3$ surfaces. We
will refer to $\calK3_{M,\rho}^{m}$
as the moduli space of $(M,\rho)$-polarized $K3$ surfaces.

\subsection{Remark.} \label{fibpim}
If $[(X,\phi)],\,[(X',\phi')]\in \calK3_{M,\rho}^{m}$ are in
the same fibre of $\pi_M$, then the $K3$ surfaces $X$ and $X'$
are isomorphic.
This follows from the surjectivity of the map
$\tilde{\Gamma}_\rho\rightarrow \Gamma_\rho$ and the Torelli theorem
for $K3$ surfaces. Let $\alpha\in \O(D(M))$.
As we already noticed in the proof of Proposition \ref{surjection},
we can lift $\alpha$ to an isometry $\tilde{\alpha}$ of $M$.
Composing it with some element of $W(M)$ which acts identically on
$D(M)$, we may assume that $\tilde{\alpha}$ leaves $\Delta(M)^+$
invariant. Now $\alpha$ acts on $[(X,\phi)]\in \calK3_{M,\rho}^m$
by $[(X,\phi)]\mapsto [(X,\phi\circ\tilde{\alpha}^{-1})]$.
This describes the action of $\O(D(M))$ on $\calK3_{M,\rho}^m$.
If $\phi(M) = \Pic(X)$, then $\O(D(M))$ acts transitively on the
markings of $X$. Thus we can interpret a general point of
$\calK3_{M,\rho}$ as the isomorphism class of a $K3$ surface which
admits an ample $(M,\rho)$-marking.

\subsection{}\label{Generic $K3$ surfaces} Recall that the subspaces $V_+$
and $V_-$ (see \ref{The period domain.})
are defined over $\bbQ (\zeta_3)$ where $\zeta_3$ is a primitive cube root
of unity. Let $K$ be the extension field of
$\bbQ (\zeta)$ obtained by adjoining all primitive $6l$-th roots of unity
for which the value of the Euler function
satisfies $\varphi (6l) \leq 10 = \text{rank}(T)$. The only possible values
of $l$ are as follows: $l = 1,2,3,4,5$. We
consider the union $\calW$ of hyperplanes of $\bbP(V_+)$ defined over $K$.
A non-singular cubic surface $S$ is called
{\it generic} if the period of the associated
$K3$ surface $X_S$ is contained in the complement of $\calW$. For example,
a cubic surface with an Eckardt point is
    not generic (we shall show in \ref{eck} that
the period of $X_S$ is contained in the hyperplane orthogonal to some
vector $r \in T$).

\subsection{Lemma}\label{Lem}
{\it Assume that $S$ is a generic cubic surface and let $X_S$ be the
associated $K3$ surface. Then the image of the
natural map $$ \Aut(X_S)\longrightarrow \O(T)
$$
is a cyclic group of order $6$ generated by $\tau$ and $\sigma$ $($ For
$\tau$, $\sigma$, see $\ref{cauto}, \ref{P1'})$.
In particular the image of the natural map $$ \Aut(X_S)\longrightarrow
\O(D(T))
$$
is $\{ \pm 1\}$.}

\begin{proof}
The proof is similar to the one given in \cite{BP}, Lemma 2.9. It is known
that the image $G$ is a cyclic group
(Nikulin \cite{N3}, Theorem 3.1). Let $m$ be the order of $G$.
If $g \in \Aut(X_S)$ is a generator of $G$, then $g^*\omega_X = \zeta_m
\cdot \omega_X$ where $\omega_X$ is a
nowhere vanishing holomorphic 2-form on $X=X_S$ and $\zeta_m$ is a
primitive $m$-th root of unity. Since $\tau^* \omega_X
= -\omega_X$ and $\sigma^* \omega_X = \zeta_3 \omega_X$, $m$ is divisible
by $6$. Since $g^*$ is defined over $\bbQ$,
the eigenspaces of $g^*$ are defined over $\bbQ(\zeta_m)$. If $m > 6$, then
an eigenspace is a non-trivial subspace of
$V_+$. This contradicts the assumption of genericity of $S$.
$\sigma^*$ acts trivially on $D(T)$ and
$\tau^*$ acts as $-1$. Hence the second assertion follows. \end{proof}

\subsection{Corollary}\label{weyl}
{\it The
map $\pi_M:\calK3_{M,\rho}^m\to \calK3_{M,\rho}$ is a Galois cover with the
Galois group isomorphic to $W(E_6)$.}

\begin{proof} As we explained in \ref{modulispaces} the group
$\O(D(T))/\{\pm 1\}\cong W(E_6)$ acts on $\calK3_{M,\rho}$
with quotient isomorphic to $\calK3_{M,\rho}$.
The isotropy subgroup of $[(X,\phi)]$ is isomorphic to the image of
$\Aut(X)$ in $D(\phi(M)^\perp)/\{\pm 1\}$.
By the previous lemma it is trivial for a generic surface $X$.
\end{proof}

\subsection{Nef divisors}\label{ample}
Let $(X,\phi)$ be an ample $M$-polarized $K3$ surface. Then $X$ has an
automorphism $\sigma$ of order 3
(\ref{k3autos}).
For any $v \in M$ with $v^2 \ge 0$ there is a $w \in W(M)$ such that
$\phi(w(v)) \in C(M)$.
If $\phi(w(v))$ is not nef, then there is a smooth rational curve $R$ with
$(R, \phi(w(v))) < 0$.
Since $\phi(M)^{\perp} \cap \Pic(X)$ does not contain $(-2)$-vectors, $R =r
+ r'$ where $r \in M^*$,
$r' \in T^*$ and $r^2 < 0, (r')^2 < 0$. Since $r^2 + (r')^2 = R^2 = -2$,
$r^2 = -2/3$ or $-4/3$.
Since $\sigma$ is an automorphism, $(R, \sigma (R)) \ge 0$. Hence $(3r)^2 =
(R + \sigma (R) + \sigma^2 (R))^2 \ge -6$.
Thus $r^2 = -2/3$. Then $r$ defines a reflection
$$
s_r : x \longmapsto x + 3(x,r)r
$$
which acts trivially on $T$. Obviously $(R, \phi(s_r(w(v))) ) > 0$. If
necessary, by
using these reflections successively, we may assume that $\phi(w(v)) \in
C(X)$, i.e.,
$\phi(w(v))$ is nef. In particular, any primitive isotropic vector $f$ in
$M$ defines, uniquely, a nef divisor in
$\Pic(X)$. As is well-known a primitive nef divisor $F$ with $F^2=0$
defines an elliptic fibration with the
cohomology class of a fibre equal to $F$ (\cite{PS}, \S 3, Cor.3).

\subsection{Elliptic fibrations.}\label{standard} Let $(X,\phi)$ be an
ample $M$-polarized $K3$ surface. With the
definitions from \ref{fixvm}, we have $f_1\in C(M)$ and $f_1$ is obviously
isotropic and primitive. Therefore,
$\phi(f_1)\in \Pic(X)$ defines an elliptic fibration on $V$ (cf.\
\ref{ample}) which we denote by $$
\Phi_\phi:X\longrightarrow \bbP^1
$$
and we call it the \emph{standard elliptic fibration}. Since
$\phi(f_2-f_1)\cdot\phi(f_1) = ( f_2-f_1,f_1) = 1$,
the divisor class $\phi(f_2-f_1)$ is an effective class with $D^2 = -2$.
Let $D$ be the effective representative of this
class written as a sum $\sum n_iR_i$, where $R_i$ are irreducible curves.
Since $D$ intersects any fibre $F$ with
multiplicity 1, we see that one of the components $R_i$, say $R_1$, is a
section of the fibration. We also have $n_1 =
1$ and $R_i\cdot F = 0$ for $i > 1$. By the Hodge Index Theorem, $R_i^2 <
0$ for $i > 1$. By the adjunction formula, all
$R_i$'s are $(-2)$-curves and the $R_i$'s, $i\ne 1$, are contained in
fibres of the fibration. This easily implies that
$R_1$ is determined uniquely by $\phi(f_2-f_1)$. We shall denote the
section corresponding to $R_1$ by $s$.
We remark that $R_1$ is obtained from $D$ by applying suitable reflections
corresponding to $R_i$ ($i > 1$).  Thus, up to isometries, we may assume that
the classes $f_1$ and $f_2 - f_1$ define an elliptic fibration
$\Phi_\phi$ with a section $s$.

The images under $\phi$ of the simple root bases $\{r_i,r_i'\}, i =
1,\ldots,5,$ of each copy of $A_2$ are
effective divisor classes $R_i,R_i'$ on $X$ which are orthogonal to $F$ and
to the section $s$. As above we can show
that each such divisor class is a sum of $(-2)$-curves contained in a
fibre. Thus
$X$ has at least 10 smooth rational curves contained in fibres of
$\Phi_\phi$.

\subsection{Lemma.}\label{singularfibre} {\it Let $(X,\phi)$ be an ample
$(M,\rho)$-polarized $K3$ surface,
let $\sigma$ be an automorphism of order three as in \ref{k3autos} and let
$\Phi_\phi$ be the standard elliptic
fibration on $X$.

Then $\sigma$ preserves $\Phi_\phi$ and fixes pointwisely
its section $s$ and a smooth bisection $b$.
Moreover, the singular fibres of $\Phi_\phi$ are of the following types:
$$
(II, II, IV, IV, IV, IV, IV),\quad (II, IV, IV, IV, IV, I_0^*),\quad (IV,
IV, IV, I_0^*, I_0^*).
$$
In each case the fibration has exactly 5 reducible fibres. }

\begin{proof}
Let $X^\sigma$ be the fixed locus of the automorphism $\sigma$. Since
$\sigma$ can be locally linearized,
$X^\sigma$ is a smooth closed subset of $X$. It is easy to see that the
trace of $\rho$ in its action on $L\cong
H^2(X,\bbZ)$ is equal to 7. Applying the Lefschetz fixed point formula, we
obtain that the Euler characteristic of
$X^\sigma$ is equal to 9. Since $\sigma$ acts identically on $\phi(M)$, it
preserves the section $s$ and the divisor
class of a fibre of $\Phi_\phi$. Let us show that $\sigma$ fixes the
section $s$ pointwisely, or, equivalently, leaves
invariant each fibre of $\Phi_\phi$. Assuming otherwise, we obtain that
$X^\sigma$ is contained in fibres of
$\Phi_\phi$. Thus any irreducible one-dimensional component of $X^\sigma$
has the Euler characteristic equal to 0 (if it
is nonsingular fibre) or $2$ (if it is a component of a reducible fibre),
the smoothness of the fixed point set excludes
nodal cubics. Let $l$ be the number of irreducible one-dimensional
components of $X^\sigma$ different from a fibre, and
let $k$ be the number of isolated fixed points. Then $2l+k = \chi(X^\sigma)
= 9$. Since $\sigma$ has exactly two fixed
points on $s$, it leaves invariant the two fibres $F_1,F_2$ passing through
these points. Obviously the curves
$R_i,R_i'$ (see \ref{standard}) are contained in the union $F_1\cup F_2$.
In particular, the number of irreducible
components of the divisor $F_1+F_2$ is greater than or equal to $12$. Since
a Dynkin diagram of type ADE admits a
non-trivial automorphism of order 3 only in the case $D_4$, the
automorphism $\sigma$ acts identically on the set of
irreducible components of a fibre $F_i$ unless it is of type $I_0^*$. Note
that either $F_1$ or $F_2$ is not of type
$I_0^*$ because $F_1 + F_2$ has at least 12 components. Assume that both of
the $F_i$'s are not of this type. We apply
the Lefschetz fixed point formula to the cell complex $F_i$. Let $n_i$ be
the number of irreducible components of $F_i$.
The Lefschetz number of $\sigma|F_i$ is equal to $n_i$ if $F_i$ is of type
$I_n$ and to $n_i+1$ otherwise. Let $l_i$ be
the number of one-dimensional rational components of $X^\sigma$ contained
in $F_i$ and let $k_i$ be the number of
isolated fixed points of $\sigma$ contained in $F_i$. We have $2l_i+k_i \ge
n_i$, hence $9 = 2l+k\ge
2l_1+k_1+2l_2+k_2\ge n_1+n_2\ge 12$, a contradiction. Assume that one of
the fibres, say $F_1$ is of type $I_0^*$. Then
$2l_2+k_2\ge n_2\ge 12-5 = 7$. The automorphism $\sigma$ has a fixed point
on the non-multiple component $E$ of $F_1$
which is intersected by $s$. The multiple component $E_0$ of $F_1$ is
$\sigma$-invariant. If $\sigma$ is the identity on $E_0$, then $l_1,k_1\ge
1$, and $2l_1+k_1\ge 3$.
If $\sigma$ does not acts identically on $E_0$, it has 2 fixed points on
it. In both cases it is easy to see that
$2l_1+k_1\ge 3$ again. Thus we get $2l_1+k_1+2l_2+k_2\ge 3+n_2\ge 3+7 =
10$, again a contradiction.

Now we know that $\sigma$ preserves every fibre of $\Phi_\phi$, so that the
general fibre has a non-trivial
automorphism of order 3 over the function field of the base. This implies
that the $j$-function of the fibration is
constant. In particular, the singular fibres must be of additive type $II$,
$III$, $IV$, $IV^*$, $II^*$, $III^*$,
$I_n^*$. Each nonsingular fibre has exactly 3 fixed points of
$\sigma$, one lies on the section $s$, and the pairs of others lie on a
bisection $b$ (which could be the union of
two sections). The bisection $b$ is a part of $X^\sigma$ and hence smooth.

Let $\pi:X'\to X$ be the blow-up of the
$0$-dimensional part of $X^\sigma$. We know that $\sigma$ is not symplectic
(i.e. does not leave invariant a non-zero
holomorphic 2-form on $X$). This easily shows that it lifts to an
automorphism $\sigma'$ of $X'$ with $X'{}^{\sigma'}$
purely one-dimensional. Let $\bar{X}'$ be the quotient surface
$X'/(\sigma')$. It is a smooth surface. Let $C$ be a
smooth rational curve on $X$ such that $\sigma(C) = C$ but $\sigma|C$ is
not the identity. Then $\sigma$ has two fixed
points $p,q$ on $C$. If $p,q$ are isolated fixed points of $\sigma$ on $X$,
then the proper inverse transform $C'$ on
$X'$ has self-intersection $-4$. Since $C'$ is equal to the pre-image of
some curve on $\bar{X}'$ and $-4$ is not
divisible by 3, we get a contradiction. Similarly, if $p,q$ belong to
one-dimensional part of $X^\sigma$, we get $C'{}^2
= -2$ and again get a contradiction. Thus, one fixed point is an isolated
fixed point of $\sigma$ and another one
belongs to the one-dimensional part of $X^\sigma$.

As we have already observed before, $\sigma$ acts identically on the set of
irreducible components of any fibre,
unless it is of type $I_0^*$. In particular, all intersection points of
components are fixed. The previous observation
about the intersection of $X^\sigma$ with fibres easily excludes the
possibility for a reducible fibre of $\Phi_\phi$ to
be of type $III$, $I_n^* (n\ne 0)$, $III^*$. In the case of $I_0^*$,
$\sigma$ preserves the multiple component $E$ and
permutes the three simple components $E_1, E_2, E_3$ not meeting the
section. Notice that any $\sigma$-invariant
irreducible component of a fibre not intersecting the section $s$ must
belong to $\phi(M)\cap \phi(U)^\perp =
\phi(A_2^5)$. The fixed part of $D_4 = <E, E_1, E_2, E_3>$ under $\sigma^*$
is $<E, E+E_1+E_2+E_3> \cong A_2$. Since
$E_6$ and $E_8$ can not be embedded into $A_2^5$, singular fibres of type
$IV^*$, $II^*$ do not appear.

Using that the Euler characteristics of the fibres add up to 24, it remains
to show that we have exactly 5 reducible
fibres. Since a fibre of type $I_0^*$ or
$IV$ contributes one copy of $A_2$ in $A_2^5\cong \phi(M)\cap
\phi(U)^\perp$, there must be five of them. The lemma is
now proven. \end{proof}

\section{A complex ball uniformization}\label{sec5}\label{cbu}

\subsection{From $K3$'s to cubics.}\label{construction} We are going to
construct a map $$
G:\,\calK3_{M,\rho}^m\setminus\Delta^m\longrightarrow \calM_{\cub}^m, $$
where $\calM_{\cub}^m$ is the moduli space of
marked smooth cubic surfaces,
i.e., smooth cubic surfaces with an ordered set of six skew
lines $L_1,\ldots,L_6$.

Let $[(X,\phi)]\in \calK3_{M,\rho}^m\setminus\Delta^m$ be an ample
$(M,\rho)$-polarized $K3$ surface.
We use the notation of Lemma \ref{singularfibre} and its proof.
For simplicity we consider the case where $\Phi_\phi$ has two singular
fibres of type $II$ and five singular fibres
of type $IV$. The construction for the other two cases is similar. It
follows from the proof of  lemma
\ref{singularfibre} that on each reducible fibre $\sigma$ has one fixed
point, the point of intersection of the three
components. The bisection $b$ intersects two components, and the section
$s$ intersects the third one. Let $X'$ be the
blow-up of the five isolated fixed points of $\sigma$ as in the proof of the
lemma. The quotient
$\bar{X}'$ of $X'$ by the action of $\sigma$ is a smooth rational surface
and the images of the components of the fibers
of type $IV$ are $(-1)$-curves in $\bar{X}'$. The marking $\phi$ gives an
ordering of the 2 components in each fibre
which meet the bisection $b$, and we blow down the first one in each of the
5 fibres as well as the component in the
fibre which meets the section. The result is a smooth rational surface $S$
which has $(-1)$-curves $L_1,\ldots, L_5$
  the images of the remaining components in the type $IV$ fibres
(these are numbered by the marking $\phi$) as
well as the $(-1)$-curve $m$ which is the image of the section $s$. These
six curves do not intersect and thus can be
blown down to get a smooth rational surface with $b_2 =1$, hence this
surface must be $\bbP^2$. Therefore $S$ is a cubic
surface and the six $(-1)$-curves define a marking on $S$. It is easy to
see that this marked cubic surface $S$ depends
only on the isomorphism class of $(X,\phi)$. We may now define:
$$
G:\,[(X,\phi)]\longmapsto (S,\,L_1,\ldots,L_5,L_6 = m). $$ Note that the
$2$-section $C$ maps to a line $l$ in $S$
which is skew with $m$ and does meet $L_1$, $\ldots$, $L_5$. By the
uniqueness of the triple cover (Theorem \ref{independence}) we have that
$X\cong X_{S,l,m}$ and, by construction (see \ref{k3autos})
$\sigma^*=\tilde{\phi}\circ\rho\circ\tilde{\phi}^{-1}$ for some extension
$\tilde{\phi}:L\rightarrow H^2(X,\bbZ)$ of $\phi$.

\subsection{Theorem.} \label{main}
{\it
The map $G$ defines a $W(E_6)$-equivariant isomorphism $$
G:\calK3_{M,\rho}^m\setminus\Delta^m
\stackrel{\cong}{\longrightarrow} \calM_{\cub}^{m}. $$
}

\begin{proof}
We first construct the inverse map
$$
G^{-1}:\calM_{\cub}^{m}\longrightarrow \calK3_{M,\rho}^m\setminus\Delta^m.
$$
Given $(S,L_1,\ldots,L_6)\in \calM^{m}_{\cub}$,
let $m=L_6$ and let $l$ be the (unique) line which meets $L_1$, $\ldots$,
$L_5$ but not $m$ (if we blow down the $L_i$
to points $x_i\in \bbP^2$, $l$ maps to the conic on $x_1,\ldots,x_5$).

Let $X_{l,m}$ be the $K3$ surface associated to $(S,l,m)$ and let
$f:X_{l,m}\rightarrow \bbP^1$ be the elliptic
fibration from subsection \ref{ellfib}.
We define a
marking $\phi_{l,m}:M\rightarrow \Pic(X_{l,m})$ as in the proof of Lemma
\ref{sigmainv} by fixing an order on
the set of reducible fibres and the order on the set of components of
fibres of type $IV$ which do not intersect the
section $s$. Thus
$\phi(f_1)$ is the class of a fibre of $f$ and $\phi(f_2)$ is the sum of
the class of a fibre and the class of
the section (see \ref{standard}).
The image of $r_1$ in the $i$-th copy of
$A_2\subset M$ is the first component of the
$i$-th fibre if it is of type $IV$, and it is the divisor class
$E+E_1+E_2+E_3$ if the $i$-th fibre is of type $I_0^*$
(see the notation in the proof of Lemma \ref{singularfibre}).

The $K3$ surface $X_{l,m}$ is a triple cyclic covering of $S$ with an
automorphism $\sigma$.
We proved in Lemma \ref{sigmainv} that $\sigma^*$ acts identically on
$\phi(M)$ and has the trace $-5$ on
$\phi(M)^\perp$. This implies that $\sigma^*$ has no
eigenvectors in
$\phi(M)^\perp\otimes\bbQ$, and hence $\phi(M)^\perp$ is a free module of
rank 5 over the ring of
Eisenstein integers $\bbZ[\zeta_3]$. In particular, the maps $\sigma^*$
glue to a locally constant map on the
local system with fibers $H^2(X_{l,m},\bbZ)$.
The construction of the map $G$ is such that if
$(S',L'_1,\ldots,L'_6)=G(X,\phi)$ for some $(X,\phi)$, then
    $\rho = \tilde{\phi}^{-1}\circ
\sigma_{S'}^*\circ\tilde{\phi}$ where
$\tilde{\phi}:L\to H^2(X_{l,m},\bbZ)$ is a cohomology marking of $X$ such
that $\tilde{\phi}|M = \phi$ and
$\tilde{\phi}(T) = \phi(M)^\perp$. As $\sigma^*$ is locally constant we
conclude that there is an extension
$\tilde{\phi}_{l,m}$ of the marking $\phi_{l,m}$ such that $\rho =
\tilde{\phi}_{l,m}^{-1}\circ\sigma^*\circ\tilde{\phi}_{l,m}$. This shows
that $G^{-1}[(S,L_1,\ldots,L_6)]:=
[(X_{l,m},\phi)]$ belongs to $\calK3_{M,\rho}^m\setminus\Delta^m$. It is
obvious that $G^{-1}$ is the inverse of $G$.

We show that $G^{-1}$ is $W(E_6)$-equivariant, then $G=(G^{-1})^{-1}$ is
obviously equivariant as well.
The group $W(E_6)$ acts on $\calM_{\cub}^{m}$ in the standard way via
symmetries of the set of lines and
$W(E_6)=\Gal(\calM_{\cub}^{m}/\calM_{\cub})$. Let $\mu:
\Gal(\calM_{\cub}^{m}/\calM_{\cub})\to
\Aut(\calK3_{M,\rho}^m\setminus \Delta^m)$ be the action defined
via the isomorphism
$G^{-1}$, obviously $\mu$ is injective. Let
$S\in \calM_{\cub}$, the main result of the section 3
(Theorem \ref{independence}) was that $X_{l,m}$ is
independent of the choice of the lines
$l$, $m$ in $S$, hence $\mu(g)$ is a covering transformation of
$\calK3_{M,\rho}^m\setminus\Delta^m\rightarrow
\calK3_{M,\rho}\setminus\Delta$ for any $g\in W(E_6)$.
Thus we have an injection:
$$
\mu:W(E_6)\,\cong\,
\Gal(\calM_{\cub}^{m}/\calM_{\cub}) \longrightarrow
\Gal(\calK3_{M,\rho}^m/\calK3_{M,\rho}).
$$
Since $\Gal(\calK3_{M,\rho}^m/\calK3_{M,\rho})\cong W(E_6)$ (see \ref{weyl}),
$\mu$ is an isomorphism.
\end{proof}

\subsection{}\label{W(D_5)}
The moduli space of cubic surfaces ${\calM}_{\cub}$ is the
quotient of $\calM_{\cub}^m$ by $W(E_6)$.
Let $W(E_6)_l\subset W(E_6) \subset \Aut(\Pic(S))$
be the subgroup which fixes the class of a line $l$ on $S$. It is
well-known that $W(E_6)_l\cong W(D_5)$, which is
the semi-direct product of $(\bbZ/2)^4$ and $S_5$.

The action of $S_5\subset W(D_5)$ on a marking $(L_1,\ldots,L_6=l)$ of a
cubic surface is by permuting the first
$5$ lines. The group $W(D_5)$ is generated by these permutations and an
element $c_{123}$ of order two which acts as
the standard Cremona transformation on
$\bbP^2$ defined by the points $p_1,p_2$ and $p_3$ where $\pi:S\to\bbP^2$
is the blow down of
the $L_i$ and  $p_i=\pi(L_i)$. Thus $c_{123}$ maps
$L_1$ to $L_1'$, the strict transform of the line on $p_2$ and $p_3$, and
it fixes $L_4,L_5$ and $L_6$.
It also permutes the $2\cdot 5$ lines on $S$ which meet $l$. Let $l_i$ be
the line which maps to the line through $p_i$
and $p_6$ and let $m_i$ be the conic through all $6$ points except $p_i$.
Then $c_{123}$ fixes the $l_i$ and $m_i$
except for permuting $l_4\leftrightarrow m_5$ and $l_5\leftrightarrow m_4$.
This implies that an element in $W(D_5)$
permutes the indices and exchanges an even number of $l_i$ with an even
number of $m_i$.

\subsection{}\label{D5marking}
Recall  from Proposition \ref{surjection} that
$$
\Gamma_\rho/\Gamma_{M,\rho}\cong O(D)\cong W(E_6)\times\{\pm 1\} $$
acts on the discriminant lattice $D=D(T)\cong\bbF_3^5$. The subgroup of
$O(D)$ which consists of
isometries preserving an unordered basis (up to signs) of $D(T)$ is
isomorphic to $W(D_5)\times\{\pm 1\}$.
This provides us with a natural copy of $W(D_5)$ in
$\Gamma_\rho/\Gamma_{M,\rho}$.
Let $\Gamma_{M,\rho}'$ be the inverse image in $\Gamma_\rho$ of this
subgroup.
The
group $\Gamma_{M,\rho}'$ acts on $\calK3_{M,\rho}^m$ by changing the
markings without changing the standard
elliptic fibration defined by the marking.
Since $W(D_5)$ is a maximal
subgroup of $W(E_6)$ we
see that any
$w\in W(E_6)\setminus W(D_5)$ does not preserve the isomorphism class of
the standard elliptic fibration. This implies the following corollaries:

\subsection{Corollary.} \label{cor2.2}
{\it
Let $\calM_{\cub}$ be the moduli space of cubic surfaces. There are
isomorphisms
\[
(\calB\setminus\calH)/\Gamma_{M,\rho} \cong \calK3_{M,\rho}\setminus \Delta
\cong
\calM_{\cub}.
\]
Let $\calM_{\cub}\li$ be the
moduli space of cubic surfaces with a line. There are isomorphisms \[
(\calB\setminus\calH)/\Gamma_{M,\rho}'\;\cong\; (\calK3_{M,\rho}^m\setminus
\Delta^m)/W(D_5) \;\cong\; \calM_{\cub}\li. \]
as well as a birational isomorphism
\[
\calB/\Gamma_{M,\rho}'\simeq \calM_{\cub}\li \]
where $\Gamma_{M,\rho}'$ is the inverse image of $W(E_6)_l\times\{\pm 1\}
\subset W(E_6)\times\{\pm 1\}\cong
\Gamma_\rho/\Gamma_{M,\rho}$ in $\Gamma_\rho$.}

\subsection{Corollary.}\label{27lines}
{\it Assume that $S$ is a generic cubic surface. Then $X_S$ has exactly $27$
( = the index of $W(D_5)$ in $W(E_6)$) non-isomorphic standard elliptic
fibrations. }

\section{The geometry of the discriminant locus.}
\label{discriminant}

\subsection{} \label{Ext}
Here we will give a geometric interpretation of the points in
$\calK3_{M,\rho}^m$ belonging to the discriminant locus
$\Delta^m$.
We know that each such point represents the isomorphism class of a
non-amply $M$-polarized $K3$ surface $(X,\phi)$.
For such a surface there is a $(-2)$-vector $r$ in $\phi(M)^\perp\cap
\Pic(X)$. This implies that $\rho$ (cf.
\ref{markfib}) can not be represented by an automorphism of $X$. Let $R$ be
the sublattice of $\Pic(X)$ generated by all
$(-2)$-vectors in $\phi(M)^\perp\cap \Pic(X)$.
Then $R$ is a negative definite lattice generated by $(-2)$-vectors,
i.e., a root lattice.
Hence $R$ is an orthogonal direct sum
$$
R = R_1 \oplus \cdot \cdot \cdot \oplus R_r,
$$
where $R_i$ is an indecomposable root lattices of type
$A_m, D_n, E_k$.
Obviously $\rho$ preserves $R$.
Since $\rho$ has no non-zero fixed vectors in $R$,
$\rho$ preserves each $R_i$. Thus
$R_i$ is an indecomposable root lattice with an isometry of order 3
without non-zero fixed vectors.
In the following we shall show that $R_i \cong A_2$ and $r \le 4$
(see \ref{level}).

\subsection{Lemma.}\label{lemma-dscri}
$R_i \cong A_2$ {\it for any} $i$.

\begin{proof}
First of all, note that the rank of $R_i$ is even because it has an
isometry of order 3 without non-zero fixed vectors.
Since the rank of $\Pic(X) \le 20$, $R_i$ is isometric to $A_{2n}$,
$D_{2n}$, $E_6$ or $E_8$ ($n \le 4$).
Let $K$ be a primitive sublattice of $H^2(X,\bbZ)$ generated by
$M$ and $R$. Let $l(K)$ be the number of minimal generator of the
3-elementary subgroup of $K^*/K$.
Then $K^*/K \cong (K^{\perp})^*/K^{\perp}$ and
$l(K) = l(K^{\perp}) \le \text{rank}(K^{\perp})$.
Using this observation and the fact $l(M) = 5$, we can easily see that $R$
is isometric to
$D_{4}$, $A_{2}^{\oplus n}$ ($1 \le n \le 4$) or $E_6$ (for example if $R =
E_8$, then $K = M \oplus E_8$ and $l(K) =
5$. This contradicts to the fact $l(K^{\perp}) \le \text{rank}(K^{\perp})
= 2$).  Next we shall show that $R$ is not isometric
to $D_4$. In this case $K = M \oplus D_4$ and
the elliptic fibration defined by $M$-marking has five singular fibres of
type $IV$ and one of type $I_0^*$.  This contradicts the fact that the
Euler number of $K3$ surface is 24.  By the same argument, the case
$R=E_6$ does not occur.
\end{proof}

\subsection{}\label{S}
We remark that all $R_i$ are 3-elementary, i.e., $R_i^*/R_i \cong
(\bbZ/3\bbZ)^l$ for some non-negative integer
$l$ and $\rho$ acts trivially on $R_{i}^*/R_{i}$.

Let $$
T'=(\phi(M)\oplus R)^\perp,\qquad S=(T')^\perp\qquad(\subset
H^2(X,\bbZ)). $$
Thus $S$ is the smallest primitive
sublattice of $H^2(X,\bbZ)$ containing $\phi(M)\oplus R$.
By definition, the lattice $T' \cap \Pic(X)$ contains no $(-2)$-vectors.

\subsection{Lemma.} \label{lemma3.2}
{\it
Let $(X,\phi)$ be an $(M,\rho)$-polarized $K3$ surface. Let $S$, $R$, $T'$
be as above.
Then $S, T'$ are $3$-elementary lattices, and $\rho$ acts trivially on
$(T')^*/T'$. Moreover
$X$ has an automorphism $\sigma'$ of order three such that
$S=H^2(X,\bbZ)^{(\sigma')^*}$. }

\begin{proof}
We have a chain of lattices:
$$\phi(M)\oplus R \subset S \subset S^* \subset (\phi(M)\oplus R)^*$$ and
$S^*/S \cong (S^*/(\phi(M)\oplus R))/(S/(\phi(M)\oplus R))$. Since $M$ and
$R$ are 3-elementary,
$S$ is a 3-elementary lattice, i.e., $S^*/S \cong (\bbZ/3\bbZ)^l$. Since
$\rho$ acts trivially on
$(\phi(M)\oplus R)^*/(\phi(M)\oplus R) \cong \phi(M)^*/\phi(M) \oplus
R^*/R$, $\rho$ acts trivially on $S^*/S$.
Since $T'$ is the orthogonal complement of $S$ in unimodular lattice
$H^2(X,\bbZ)$, $T'$ is 3-elementary and $\rho$ acts
trivially on $(T')^*/T'$ (see Nikulin \cite{N1}, Proposition 1.6.1). Hence
the isometry $(1_{S}, \rho \mid T')$ can be
extended to an isometry $\rho'$ of $H^{2}(X, \bbZ)$ (Nikulin \cite{N1},
Corollary 1.5.2). Then $\rho'$ is represented by
an automorphism $\sigma'$ of $X$ (see \cite{Na}, Theorem 3.1). \end{proof}

The following fact was first observed by Vorontsov \cite{Vor}.

\subsection{Lemma}\label{Vor}
{\it We keep the same assumption as in Lemma \ref{lemma3.2}. Define a
non-negative integer $l(T')$ by: $(T')^*/T' \cong (\bbZ/3\bbZ)^{l(T')}$.
Then } $$
\text {rank}(T')\ge 2 l(T').
$$

\begin{proof}
Let $x \in T'$. Since
$$
(x, \rho' (x)) = (\rho' (x), (\rho')^2(x)) = (\rho' (x), -x - \rho' (x)),
$$ we get $2(x, \rho (x)) = -(x,x)$.
Hence $x$ and $\rho'(x)$ generate a sublattice $A_2(m)$, where $m = (x,x)$.
   From this we can find a sublattice
$K = A_2(m_1) \oplus \cdot \cdot \cdot \oplus A_2(m_k)$ of $T'$ of finite
index. Moreover we have $(T')^*/T' \cong
((T')^*/K)/(T'/K)$. If $m_i$ is not divisible by 3, the contribution from
$A_2(m_i)$ to $l(T')$ is at most 1. In case
$m_i$ is divisible by 3, the fixed part under $\rho'$ in
$A_2(m_i)^*/A_2(m_i)$ is $\bbZ/3\bbZ$. Since $\rho$ acts
trivially on $(T')^*/T'$, the contribution from $A_2(m_i)$ is at most 1.
This implies the assertion. \end{proof}

\subsection{Lemma}\label{A2}
{\it We keep the same notation as in Lemma \ref{lemma3.2}. Then $R \cong
A_2^{\oplus r}$ and $l(S) = 5-r$.}

\begin{proof}
Let
$$R = R_1 \oplus \cdot \cdot \cdot \oplus R_r$$ be the orthogonal
decomposition of $R$ into indecomposable root
lattices $R_i$. We know that $R_i$ is isomorphic to $A_2$ (Lemma
\ref{lemma-dscri}). Obviously $R_i^*/R_i$ is
$\bbZ/3\bbZ$. Since $S^*/S \cong (S^*/(\phi(M)\oplus R))/(S/(\phi(M)\oplus
R))$, we have $l(T') = l(S) \ge (l(M)+r) - 2r
= 5-r$. On the other hand, it follows from Lemma \ref{Vor} that $10 - 2r
\ge  {\rm rank} (T') \ge 2 l(T')$. Hence $l(S) = 5 - r$.
\end{proof}

Let us summarize the previous lemmas by stating the following:

\subsection{Theorem}\label{level}
{\it Let $(X,\phi)\in \calK3_{M,\rho}^m$. Then $X$ admits an automorphism
$\sigma'$ of order 3 such that
$H^2(X,\bbZ)^{(\sigma')^*}=S$,
the smallest primitive sublattice of $\Pic(X)$ which contains $\phi(M)$ and
the sublattice $R$ generated by all
$(-2)$-vectors in $\phi(M)^\perp\cap\Pic(X)$. The sublattices $\phi(M)$ and
$R$ are orthogonal to each other and the
lattice $R$ is isomorphic to $r\;(\le 4)$ copies of the lattice $A_2$. The
number $r$ will be called the
\emph{degeneracy rank} of $(X,\phi)$
}.

\medskip
The degeneracy rank of $(X,\phi)$ is  equal to the number of
nodes of the associated nodal cubic surface
(see \ref{conics}). This is easy to see from Table 2 by computing the
quotient of  $M(\bft)$ by $M = U\oplus A_2^{\oplus 5}$ and comparing the
result with the value of $r$ in Table 1.
The next theorem generalizes Lemma \ref{singularfibre}.

\subsection{Theorem}\label{T}
{\it Let $[(X,\phi)] \in \calK3_{M,\rho}^m$.
Then the $M$-marking $\phi$ of $X$
defines an elliptic fibration. Its  singular fibres are given in
the column Kodaira fibres of Table 1 from above.
The Picard lattice $S_X$ and  its lattice of transcendental
cycles $T_X$ can be found in the corresponding rows of Table 2
(under the assumption in Proposition \ref{ptlatticev}).
The degeneracy level is given in the column $r$ in Table 1. }

\begin{proof} By the same arguments as in \ref{ample}, \ref{standard}, the
$M$-marking on $X$ defines an elliptic
fibration with a section. The proof of the assertion about possible
combinations of singular fibres is very similar to
the proof of Lemma \ref{singularfibre} and is omitted. The description of
the transcendental lattice follows from the
following easy facts:
$$
q_{E_6} = -q_{A_2}, \quad q_{A_2(-1)} = -q_{A_2},\quad q_{A_2} \oplus
q_{A_2} = q_{A_2(-1)} \oplus q_{A_2(-1)},
\quad q_{A_2(-2)} = q_{D_4} \oplus q_{A_2}
$$
and Theorem 1.14.2 from \cite{N1}.
\end{proof}

\subsection{The Eckardt locus}\label{eck} Let $[(X,\phi)]\in
\calK3_{M,\rho}^m\setminus \Delta^m$. We know
that the corresponding marked cubic surface $(S,L_1,\ldots,L_6)$ has an
Eckardt point on the unique line $l$
intersecting $L_1,\ldots,L_5$ if and only if the standard elliptic
fibration $\Phi_\phi$ on $(X,\phi)$ has a fibre of
type $I_0^*$. In that case $\phi(M)\ne \Pic(X)$, but for general $S$ with
such property, the orthogonal complement
$\phi(M)^\perp_{\Pic(X)}$ of $\phi(M)$ in $\Pic(X)$ is isomorphic to
$A_2(2)$. In fact if $F = 2E_0+E_1+\ldots+E_4$ is the fibre of type $I_0^*$
and $E_4$ meets the section,
then $\phi(M)^\perp_{\Pic(X)}$ is spanned by $E_1-E_2$ and $E_2-E_3$.

The involution $\tau$ (cf.\ \ref{P1'}) defined by the elliptic fibration
also acts on $\phi(M)$, via $\iota=\tau^*$,
in a different way. If all fibres are of type
$IV$, then the action of $\iota$ on $\phi(M) \cong U\oplus A_2^5$ permutes
the simple root basis
in each copy of $A_2$. Let $N = \phi(M)^\iota$ be the sublattice of the
invariant  elements, then \[ N \cong U\oplus A_1^5. \]
However, if one of the fibres is of type $I_0^*$, then
$\phi(M)^{\iota}\cong U\oplus A_2\oplus A_1^4$.
The orthogonal complement of $\phi(N)$ in $\phi(M)^{\iota}$ is spanned by
the class of the divisor $E_1+E_2+E_3$.
Also $r = [E_1]\in \phi(N)_L^\perp$ but not in $\phi(M)$.

For any $(-2)$-vector $r\in N^\perp\setminus T \subset L$ consider the
hyperplane $r^\perp$ in $\bbP(V_+)$
    of lines orthogonal to $r$. Let $H(r)_\iota$ be the intersection of this
hyperplane with the ball $\calB\subset
\bbP(V_+)$. Let $\calH_\iota$ be the union of the hyperplanes $H(r)_\iota$.
If an ample $(M,\rho)$-marked surface
$(V,\phi)$ has a fibre of type $I_0^*$ in its standard elliptic fibration
$\Phi_\phi$, then its period belongs to
$\calH_\iota$. Let $\Delta_\iota^m$ (resp. in $\Delta_\iota$) be the image
of $\calH_\iota$ in $\calK3_{M,\rho}^m$
(resp. in $\calK3_{M,\rho}$).  In this notation we have

\subsection{Theorem} {\it Under the isomorphism $\calM_{\cub}\cong
\calK3_{M,\rho}\setminus \Delta$,
the image of the locus of smooth cubic surfaces with Eckardt points (the
Eckardt locus) is mapped to
$\Delta_\iota \setminus (\Delta \cap \Delta_\iota)$.}

\subsection{} It is well-known that any nonsingular cubic surface contains
45 tritangent planes,
i.e.\ plane sections which split into the union of three lines. A marking
of a cubic surface defines an
order on the
set of tritangent planes. Let $\calE_i$ be the locus of points in
$\calM_{\cub}^m$ corresponding to marked cubic
surfaces which contain an Eckardt point in the $i$-th tritangent plane. The
Weyl group $W(E_6)$ acts on $\calM_{\cub}^m$
and  permutes the loci $\calE_i$'s transitively. Let $(S,L_1,\ldots,L_6)$ be
a marked cubic surface and let $M_i$ be the
line on $S$ which meets $L_i$ and $L_{i+3}$ for $i=1,2,3$ but none of the
other $L_j$. The $M_i$ lie in a tritangent
plane and they meet in a point if and only if the points
$p_1,\ldots,p_6\in\bbP^2$ obtained by blowing down the $L_i$ are such that
the three lines
$\langle p_i,p_{i+3}\rangle$ (the images of the $M_i$), intersect at some
point $q$. Let $\calE_j$ be the corresponding component of the Eckardt locus
in
$\calM_{\cub}^m$. Its pre-image $Z$ in $(\bbP^2)^6$ consists of 6-tuples of
points $(p_1,\ldots,p_6)$ such that the lines $\langle p_i,p_{i+3}\rangle, i
= 1,2,3$ intersect. The assigning the intersection point $q$ defines a
surjective map from $Z$ to $\bbP^2$  whose fibres, as is easy to see,
are irreducible.
This shows that $Z$, and hence $\calE_j$ is irreducible.
The image of each $\calE_i$
in $\calM_{\cub}$ is then an irreducible hypersurface.

The irreducibility of the Eckardt locus in \ref{eckloci}
  follows also from  our ball
uniformization of $\calM_{\cub}$.
We follow the proof given in \cite{AF}.

\subsection{Lemma}\label{D}
{\it Let $D=T^*/T$ be the discriminant group of $T$ as in \ref{lattices}
and let $N=M^{\iota}$.
The group $W(E_6) = \O(D)/\{\pm1\}$ acts transitively on the subsets of
$(D-\{0\})/\{\pm 1\}$ of vectors of the same
norm.  There are three such subsets.
\begin{enumerate}
\item [(i)]
The set of vectors of norm $0$ has $40$ elements. Each non-zero isotropic
vector is represented by
$(e +2\rho(e))/3,$ where $e \in T$ is a primitive isotropic vector. \item
[(ii)] The set of vectors of norm $-2/3$ has
$36$ elements. Each $(-2/3)$-vector is represented by a vector $(r + 2\rho
(r))/3$ in $T^*$ with $r \in T$, $r^2= -2$
and $(r,\rho (r) ) = 1$. \item [(iii)] The set of vectors of norm $-4/3$
has $45$ elements. Each $(-4/3)$-vector in
$D(T)$ is represented by $r''$ where $r = r' + r'' \in N^{\perp} \setminus
T$ is a $(-2)$-vector and $r'$, $r''$ is the
projection of $r$ into
$(N^{\perp} \cap M)^*$, $T^*$ respectively.
\end{enumerate} }

\begin{proof}
If we consider $T$ as a free Hermitian module $\Lambda$ over
$\bbZ[\zeta_3]$ (see \ref{af}), then \cite{ACT,AF}
define an isotropic vector, a short vector and a long vector as a vector
with Hermitian square equal to $0,-1,-2$,
respectively. The images of these vectors in $T^*$ with respect to the
isomorphism $h:\Lambda \to T^*$ \eqref{fh} are
vectors with square $0,-\frac{2}{3},-\frac{4}{3}$, respectively. It is
proven in \cite{AF}, Proposition 2.1 that there
are exactly three $\Gamma_\rho$-orbits of the images of these vectors in
$D(T)$. Their cardinality is 40, 36 and 45,
respectively. This gives three orbits of $\O(D(T))$ in
$D(T)$ of the same cardinality.
The assertions (i) and (ii) follow from the explicit formula for the
isomorphism $h$ \eqref{fh}. To prove (iii),
we consider an ample $(M,\rho)$-polarized $K3$ surface $X$ whose standard
elliptic fibration acquires a fibre of type
$I_0^*$. Let $\tilde{\phi}:L\to H^2(X,\bbZ)$ be a cohomology marking with
$\tilde{\phi}|M = \phi$. In the notation of
\ref{eck}, we may assume that the image of the first copy of $A_2$ of $M$
in $\Pic(X)$ is spanned by
$E_0$ and $E_0+E_1+E_2+E_3$. Let
$r =\tilde{\phi}^{-1}([E_1])$. Then $r\in N^\perp\setminus T$ and $r' =
\frac{1}{3}(r+\rho(r)+\rho^2(r)) =
\frac{1}{3}\tilde{\phi}^{-1}(E_1+E_2+E_3)\in (M\cap N^\perp)^*$. We easily
check that $r'{}^2 = -\frac{2}{3}.$
Then $r'' = r-r'\in T^*$ and $(r'')^2 = - \frac{4}{3}$.
\end{proof}

\subsection{Moduli interpretation}\label{inter}
Consider the three $\Gamma_\rho$-orbits of vectors from $T^*$:
\begin{enumerate} \item[(1)] $\frac{1}{3}(e+2\rho(e)),$ where
$e$ is a primitive isotropic vector in $T$;
\item[(2)] $\frac{1}{3}(r+2\rho(r)),$ where $r$ is a $(-2)$-vector
in $T$ (this corresponds to a short root in $\Lambda$);
\item[(3)] $r''$ equal to the projection of a $(-2)$-vector
$r\in N^\perp\setminus T$ (this corresponds to a long root in
$\Lambda$).
\end{enumerate}
Each vector $v\in T^*$ defines a
hyperplane $v^\perp$ in $\bbP(V_+)$ of lines orthogonal to $v$.
So, we have three $\Gamma_\rho$-orbits of such hyperplanes
corresponding to vectors from the above list. It is shown in
\cite{AF} that there is a bijective correspondence between the
$\Gamma_{M,\rho}$-orbits of these vectors and their images in
$D(T)$. Thus each $\Gamma_{\rho}$-orbit consists of $40,36,45$
$\Gamma_{M,\rho}$-orbits, respectively.

\subsection{The boundary divisors.}\label{boundary}
We know that the discriminant $\calH$ is equal to the union of
hyperplanes $H(r) = r^\perp\cap\calB$, where $r$ is a $(-2)$-vector
from $T$. For any $x\in V_+$, we can easily see that $(r,x) = 0$
if and only if $(r + 2\rho(r),x) = 0$.
This shows that the hyperplane corresponding to a vector of
type (2) in \ref{inter}
is one of the hyperplanes $H(r)$. Thus the discriminant locus
$\Delta^m$ in $\calK3_{M,\rho}^m$ consists of 36 hypersurfaces
$\Delta^m_{\alpha}$ ($\alpha \in D/\{\pm1\}$ with norm $-2/3$)
which are permuted transitively by $W(E_6)$.
The discriminant locus $\Delta$ in $\calK3_{M,\rho}$ is
irreducible.

It is well-known that the stabilizer of each $\Delta^m_{\alpha}$ in
$W(E_6)$ is $S_6 \times {\bbZ}/2{\bbZ}$
(see \ref{Boundary divisors}).
Let $\alpha_1,..., \alpha_r$ be mutually orthogonal $r$ $(-2/3)$-vectors
in $D(T)$ ($1\leq r \leq 4$).
These vectors correspond to
a sublattice $R = A_2^r$ in $T$.  Let
$$\Delta^m_{\alpha_1,...,\alpha_r} =
\Delta^m_{\alpha_1} \cap \cdot\cdot\cdot
\cap
\Delta^m_{\alpha_r}.$$
It  parametrizes
the marked $K3$ surfaces whose periods are
orthogonal to $R$.

We fix an orthogonal basis $\{ \alpha_i\}$ of $D$ such that
$q_T(\alpha_i) = -4/3$.  This defines an isomorphism of
quadratic forms
$$D \simeq {\bbF}_3^5$$
where the quadratic form $q$ on ${\bbF}_3^5$ is given by
$$q(0,...,0,1,0,...,0) = -4/3.$$
Recall that the stabilizer of a basis of $D$ in $W(E_6)$ is $W(D_5) \simeq
({\bbZ}/2{\bbZ})^4 \cdot S_5$.

On the other hand, recall that the stabilizer $G_k$ in $W(E_6)$ of
$k$ mutually orthogonal $(-2/3)$-vectors in $D$ is
$S_6 \times {\bbZ}/2{\bbZ}$,
$ ({\bbZ}/2{\bbZ})^2 \times S_4 \times  {\bbZ}/2{\bbZ}$,
$({\bbZ}/2{\bbZ})^3 \times S_3 \times  {\bbZ}/2{\bbZ}$,
$({\bbZ}/2{\bbZ})^4 \cdot S_4$ for $k =1,2,3,4$ respectively (see
\ref{Boundary divisors}).

Now we list  the $W(D_5)$-orbits of mutually orthogonal $k$
$(-2/3)$-vectors in $D$:
\medskip

(i)    There are 36 $(-2/3)$-vectors in $D$ which are divided into two
orbits.
One consists of 16 vectors containing $(1,1,1,1,1)$ and another consists of
20 vectors
containing $(1,1,0,0,0)$.

The stabilizer in $W(D_5)$ of $(1,1,1,1,1)$ is $S_5$, and that of
$(1,1,0,0,0)$ is $({\bbZ}/2{\bbZ})^3 \cdot (S_2 \times  S_3)$.
Note that the sum of indices of these groups in $G_1$ is $12 + 15 = 27$.
The orbit of cardinality 20 corresponds to markings such that the marked line
does not contain the node.
For example, if the line  corresponds to $e_6$ under a geometric marking
defined by $(e_1,\ldots,e_6)$, then the effective class corresponding to the
node could be either of type  $e_i-e_j, 1\le i <j < 6$ or
$e_0-e_i-e_j-e_k, 1\le i < j < k < 6$.
\medskip

(ii)  There are four types of  mutually orthogonal pairs of
$(-2/3)$-vectors in $D$:
$$\{(1,1,1,1,1), (1,-1,0,0,0) \},
\quad \{(1,1,1,1,1), (-1,1,1,1,1)\},$$
$$
\{(1,1,0,0,0), (0,0,1,1,0) \},
\quad \{(1,1,0,0,0), (-1,1,0,0,0)\}.$$
The stabilizer in $W(D_5)$ of the first one is $S_2 \times S_3$, the second
is
$ ({\bbZ}/2{\bbZ}) \times S_4$, the third is
$ ({\bbZ}/2{\bbZ})^2 \cdot (S_2 \times S_2) \cdot S_2$, and
the fourth is $ ({\bbZ}/2{\bbZ})^4 \cdot  (S_2 \times S_3)$.
The sum of indices of these groups in $G_2$ is $16 + 4 + 6 +1 = 27$.
\medskip

(iii)  There are three types of  mutually orthogonal triples of
$(-2/3)$-vectors in $D$:
$$\{(1,1,1,1,1), (1,-1,0,0,0) , (0,0,1,-1,0) \},
\quad \{(1,1,1,1,1), (-1,1,1,1,1),  (0,-1,1,0,0)\},$$
$$
\{(1,1,0,0,0), (-1,1,0,0,0), (0,0,1,1,0) \}.$$
The stabilizer in $W(D_5)$ of the first one is
$(S_2 \times S_2) \cdot S_2$, the second is
$ ({\bbZ}/2{\bbZ}) \cdot  (S_2 \times S_2)$, and the third is
$ ({\bbZ}/2{\bbZ})^3 \cdot (S_2 \times S_2) $.
The sum of indices of these groups in $G_3$ is $12 + 12 + 3 = 27$.
\medskip

(iv)  There are two types of  mutually orthogonal 4-tuples of
$(-2/3)$-vectors in $D$:
$$\{(1,1,1,1,1), (-1,1,1,1,1),  (0,-1,1,0,0), (0,0,0,-1,1)\},$$
$$
\{(1,1,0,0,0), (-1,1,0,0,0), (0,0,1,1,0), (0,0,-1,1,0) \}.$$
The stabilizer in $W(D_5)$ of the first one is
$ ({\bbZ}/2{\bbZ}) \cdot  (S_2 \times S_2) \cdot S_2$, and the second is
$ ({\bbZ}/2{\bbZ})^4 \cdot (S_2 \times S_2) \cdot S_2$.
The sum of indices of these groups in $G_4$ is $24 + 3 = 27$.
\medskip

Let $\Delta_k$ be the image in $\Delta$ of  all
$\Delta^m_{\alpha_1,...,\alpha_k}$
where $\{\alpha_1,...,\alpha_k\}$ is a set of mutually orthogonal $k$
$(-2/3)$-vectors in $D$.  Then
the discriminant locus $\Delta$ in $\calK3_{M,\rho}$ has the following
stratification:
$$\Delta = \Delta_1 \cup \Delta_2 \cup \Delta_3 \cup \Delta_4.$$

Similarly $\Delta' = \Delta^m/W(D_5)$ in $\calK3_{M,\rho}^m/W(D_5)$
has the following stratification
$$\Delta' = \Delta'_1 \cup \Delta'_2 \cup \Delta'_3 \cup \Delta'_4;$$
$$\Delta'_1 = \Delta_1^{(1)} \cup \Delta_1^{(2)};$$
$$\Delta'_2 = \Delta_2^{(1)} \cup \Delta_2^{(2)}
\cup \Delta_2^{(3)} \cup \Delta_2^{(4)};$$
$$\Delta'_3 = \Delta_3^{(1)} \cup \Delta_3^{(2)}
\cup \Delta_3^{(3)};$$
$$\Delta'_4 = \Delta_4^{(1)} \cup \Delta_4^{(2)}.$$
Here $\Delta_k^{(r)}$ is the image of $\Delta^m_{\alpha_1,...,\alpha_k}$
where $\{\alpha_1,...,\alpha_k\}$ is as follows:
\medskip

In case $\Delta_1^{(1)}$, $\{\alpha_1\} = \{ (1,1,1,1,1)\}$;
\medskip

In case $\Delta_1^{(2)}$, $\{\alpha_1\} = \{(1,1,0,0,0)\}$;
\medskip

In case $\Delta_2^{(1)}$, $\{\alpha_1, \alpha_2\} = \{(1,1,1,1,1),
(1,-1,0,0,0) \}$;
\medskip

In case $\Delta_2^{(2)}$, $\{\alpha_1, \alpha_2\} = \{(1,1,1,1,1),
(-1,1,1,1,1)\}$;
\medskip

In case $\Delta_2^{(3)}$, $\{\alpha_1, \alpha_2\} = \{(1,1,0,0,0),
(0,0,1,1,0) \}$;
\medskip

In case $\Delta_2^{(4)}$, $\{\alpha_1, \alpha_2\} = \{(1,1,0,0,0),
(-1,1,0,0,0)\}$;
\medskip

In case $\Delta_3^{(1)}$, $\{\alpha_1,\alpha_2,\alpha_3\} = \{(1,1,1,1,1),
(1,-1,0,0,0) , (0,0,1,-1,0)  \}$;
\medskip

In case $ \Delta_3^{(2)}$, $\{\alpha_1,\alpha_2,\alpha_3\} = \{
(1,1,1,1,1), (-1,1,1,1,1),  (0,-1,1,0,0) \}$;
\medskip

In case $\Delta_3^{(3)}$, $\{\alpha_1,\alpha_2,\alpha_3\} = \{ (1,1,0,0,0),
(-1,1,0,0,0), (0,0,1,1,0) \}$;
\medskip

In case $\Delta_4^{(1)}$, $\{\alpha_1,...,\alpha_4\} = \{ (1,1,1,1,1),
(-1,1,1,1,1),  (0,-1,1,0,0), (0,0,0,-1,1)  \}$;
\medskip

In case $ \Delta_4^{(2)}$, $\{\alpha_1,...,\alpha_4\} = \{ (1,1,0,0,0),
(-1,1,0,0,0), (0,0,1,1,0), (0,0,-1,1,0) \}.$
\medskip

Now we conclude

\subsection{Remark}\label{extension}
{\it  For each $k$,
the degree of the natural map $\Delta'_k \to \Delta_k$
is $27$.
}

\subsection{Eckardt loci.}\label{eckloci} If $v$ is of type (3) in
\ref{inter} the hyperplane
$v^\perp\cap \calB$ is equal to the hyperplane $\calH(r)_\iota$ defined in
\ref{eck}. Thus we obtain that the image of
the Eckardt locus $\Delta_\iota^m$ in $\calK3_{M,\rho}^m$ consists of 45
irreducible hypersurfaces. The Eckardt locus
$\Delta_\iota$ in $\calK3_{M,\rho}$ is irreducible. This shows that the
Eckardt locus in $\calM_{\cub}$ is irreducible
(as promised).

\subsection{Cusps.}\label{Satake}
For a non-zero isotropic vector $e$ in $T$ we define a totally isotropic
sublattice $$
I(e):=\langle e,\,\rho (e)\rangle\qquad(\subset T). $$ Then $\bar{\calB}
\cap (\bbP(I(e) \otimes \bbC)$
is a cusp of $\calB$ (i.e. a rational boundary component), and any cusp of
$\calB$ corresponding to a parabolic subgroup
of $\Gamma_\rho$ is obtained in this manner. Thus we obtain that the
Satake-Baily-Borel compactification of
$\calK3_{M,\rho}^m = \calB/\Gamma_{M,\rho}$ (resp. $\calK3_{M,\rho} =
\calB/\Gamma_{\rho}$) is obtained by adding 40
cusps (resp. one cusp).

\section{Extension of the isomorphism to the boundaries}\label{extiso}

The purpose of this section is to extend the isomorphisms
$\calK3^m_{M,\rho}\setminus\Delta^m\to \calM_{\cub}^m$ to a
$W(E_6)$-equivariant isomorphism
$$
\calK3^m_{M,\rho}\longrightarrow \calM_{\ncub}^m.
$$
First we will prove the following.

\subsection{Theorem}\label{extensionthm}{\it
The isomorphism $f:(\calK3^m_{M,\rho}\setminus \Delta^m)/W(D_5)
\longrightarrow \calM_{\cub}\li$
extends to an isomorphism:
$$
\calK3^m_{M,\rho}/W(D_5)\stackrel{\cong}{\longrightarrow} \calM_{\ncub}\li.
$$
}

\begin{proof} It follows easily from  Theorem \ref{T} that the standard
elliptic fibration defined by $(X,\phi)\in \calK3^m_{M,\rho}$ has
Weierstrass form \eqref{we}. Let $\calE$ be the moduli space of such
elliptic fibrations. This defines a natural map from $\calK3^m_{M,\rho}$
to $\calE$ which obviously factors through a map
$\calK3^m_{M,\rho}/W(D_5)$ to $\calE$.
By Corollary \ref{modell}, we have
a natural isomorphism $\calE \cong \calM_{\ncub}\li$.
It follows from the
construction that the composition of these maps restricted to the
complement of  $\Delta^m$ coincides with the isomorphism from Corollary
\ref{cor2.2}.
   Thus it suffices to show that the map
$\calK3^m_{M,\rho}/W(D_5)\to \calE$ is an isomorphism. To construct the
inverse we have to define a marking on an elliptic $K3$ surface $X$ from
$\calE$. We may assume that $X = X_{S,l}$ for some  nodal cubic surface
$S$, and the elliptic fibration is defined by the pencil of planes
through $l$. We have a line $m$ skew to $l$ such that $W(D_5)$
corresponds to its stabilizer in $W(E_6)$.

Recall that  the elliptic $K3$ surface $X_{S,l}$ is isomorphic to a minimal
resolution of the  double cover of $\bbP^2$ branched along the curve
$$x_2(F_2(x_0,x_1)x_2^3+F_5(x_0,x_1)) = 0$$
and the elliptic pencil $f:X_{S,l}\to \bbP^1$ is defined by the pencil of
lines through the point $p = (0,0,1)$. We shall define a $M$-marking of
the elliptic $K3$ surface $X_{S,l}
\simeq X_{S,l,m}$.

First note that $f$ has a section corresponding to the line
$x_2 = 0$ or $m$.  This section and the class of a fibre of $f$ define the
hyperbolic plane $U$.  Thus it suffices
to give an embedding of $A_2^5$ into $\Pic(X_{S,l})$.  If $f$ has a
singular fibre of type $IV$ or of type $I_0^*$,
the triple cover construction of $X_{S,l,m}$ determines an embedding of
$A_2$ in $\Pic(X_{S,l})$ as in the smooth case.
If $f$ has a singular fibre of
type $IV^*$ or of type $II^*$, then we take and fix embeddings
$$\phi_1 : A_2 \to \Pic(X_{S,l}), \quad \phi_2 : A_2^2 \to \Pic(X_{S,l}),
\quad \phi_3 : A_2^2 \to \Pic(X_{S,l})$$
as follows:  $\phi_1$, $\phi_2$ sends simple roots to effective
$(-2)$-curves on the singular fibre
of type $IV^*$ which do not meet the section, and  $\phi_3$ sends simple
roots to effective $(-2)$-curves on the
singular fibre of type $II^*$ which do not meet the section.
Now we define a $M$-marking of $X_{S,l}$ according to types of
lines on $S$ listed  in \ref{types} (they determine the $\SL(2)$-orbit of
the stable pair $(F_5,F_2)$ and the degenerate fibres of $f$).
\medskip

Cases 1, 2, 3:  In these cases $S$ is smooth and the marking has
already been defined in the proof of Theorem \ref{main}.
Recall that in the case 1),
the subgroup of $W(E_6)$
which preserves the elliptic fibration
is the group $W(D_5) \simeq
({\bbZ}/2{\bbZ})^4\cdot S_5$, cf.\ \ref{W(D_5)}.

\medskip

Case 4:  in this case, five roots of $F_5$ gives five singular fibres of
type $IV$.  Thus we have a $M$-marking
$$\phi : M \to \Pic(X)$$
determined by the triple cover construction
by sending $A_2^5$ to the sublattice generated by ten components of these
five
singular fibres not meeting the section.
The multi-section splits into two disjoint section and hence
the stabilizer of $\phi$ in $W(D_5)$ is $S_5$.
\medskip

Case 5:  in this case, three simple roots of $F_5$ give three singular
fibres of type $IV$ and
one multiple root of $F_5$ gives a singular fibre of type $IV^*$.
Using $\phi_2$, we have a $M$-marking $\phi$. The stabilizer of $\phi$ in
$W(D_5)$ is
$({\bbZ}/2{\bbZ})^3\cdot (S_2 \times S_3)$.
\medskip

Cases 6, 7:  these cases are mixed one of case 2) and case 5).
\medskip

Case 8: this case is the mixed one of case 4) and case 5).   The stabilizer
of $\phi$ in $W(D_5)$ is
$S_2 \times S_3$.
\medskip

Case 8*:  in this case, the common root of $F_2$ and $F_5$ gives a singular
fibre of type $IV^*$
and the remaining simple roots of $F_5$ give four singular fibres of type
$IV$.
Using $\phi_1$, we have a $M$-marking $\phi$.
The stabilizer of $\phi$ in $W(D_5)$ is ${\bbZ}/2{\bbZ} \times S_4$.
\medskip

Case 9:  this case is the mixed one of  case 5).   The stabilizer of
$W(D_5)$ is
$({\bbZ}/2{\bbZ})^2\cdot (S_2 \times S_2)\cdot S_2$.
\medskip

Case 10: the common root of $F_2$ and $F_5$ gives a singular fibre of type
$II^*$
and the remaining three roots of $F_5$ give three singular fibres of type
$IV$.
Using $\phi_3$, we have a $M$-marking $\phi$.
The stabilizer of $\phi$ in $W(D_5)$ is $({\bbZ}/2{\bbZ})^4\cdot (S_2
\times S_3)$.
\medskip

Case 11:  this case is the mixed one of case 2) and case 5).
\medskip

Case 12:  this case is the mixed one of case 2) and case 10).
\medskip

Case 13:  this case is the mixed one of case 4) and case 5).  The
stabilizer of $\phi$ in $W(D_5)$ is
$(S_2 \times S_2) \cdot S_2$.
\medskip

Case 13*:  this case is the mixed one of case 5) and case 8*).  The
stabilizer of $\phi$ in $W(D_5)$ is
${\bbZ}/2{\bbZ} \cdot (S_2 \times S_2)$.
\medskip

Case 14:  this case is the mixed one of case 5) and case 10).  The
stabilizer of $\phi$ in $W(D_5)$ is
    $({\bbZ}/2{\bbZ})^3 \cdot (S_2 \times S_2)$.
\medskip

Case 15:  this case is the mixed one of case 2), case 5) and case 10).
\medskip

Case 16:  this case is the mixed one of case 5) and case 8*).  The
stabilizer of $\phi$ in $W(D_5)$ is
${\bbZ}/2{\bbZ} \cdot (S_2 \times S_2) \cdot S_2$.
\medskip

Case 17:  this case is the mixed one of case 10).  The stabilizer of $\phi$
in $W(D_5)$ is
$({\bbZ}/2{\bbZ})^4 \cdot (S_2 \times S_2) \cdot S_2$.
\medskip

Thus we have defined a $M$-marking $\phi$ of $X_{S,l}$
modulo the action of $W(D_5)$. It
gives the period of $(X_{S,l}, \phi)$ in $\calB$
as in the case of smooth cubic surfaces.
Recall that $O(E_6)/\{\pm 1\} \simeq W(E_6)$ and $O(E_8) \simeq W(E_8)$.
Since a reflection in $O(M(\bft))$ acts trivially on the
discriminant group of
$M(\bft)$ (for $M(\bft)$ see \ref{mt}),
it can be extended to an isometry
of $L$ acting trivially on $T(\bft)$.
Now we can easily see that the image of
the period of $(X_{S,l}, \phi)$ in $\calK3^m_{M,\rho}/W(D_5)$ is independent of
the choice of $\phi_1, \phi_2, \phi_3$.
This defines the inverse map $\calE \to \calK3_{M,\rho}^m/W(D_5)$ and
proves the theorem.
\end{proof}

\subsection{Remark.}
If $S$ has a node, the period defines a point
$[X_{S,l}]$ in $\Delta'$.  More precisely:
\medskip

In the case 4,  two singular fibres of type $II$
degenerate to a singular fibre of type $IV$ and
the multi-section splits into two disjoint sections $s_1, s_2$.
The projection of the class of $s_1$ in $D$ is of type $(1,1,1,1,1)$.
Hence $[X_{S,l}] \in \Delta_1^{(1)}$.
\medskip

In the cases 5, 6, 7, two singular fibres of type $IV$ degenerate to a
singular fibre of type $IV^*$.
The projection of some component of the singular fibre of type $IV^*$ is of
type $(1,1,0,0,0)$.
Hence  $[X_{S,l}] \in \Delta_1^{(2)}$.
\medskip

The case 8 is the mixed one of cases 4 and 5.  Hence $[X_{S,l}] \in
\Delta_2^{(1)}$.
\medskip

In the case 8*, two singular fibres of type $II$ and a fibre of type $IV$
degenerate to
a singular fibre of type $IV^*$, and the multi-section splits into two
disjoint sections $s_1, s_2$.
The projections of $s_1$ is of type $(1,1,1,1,1)$ and that of a component
of the fibre of type $IV^*$
is of type $(1,0,0,0,0)$.  Combining these two $(-2/3)$- and
$(-4/3)$-vectors, we have a vector of
type $(-1,1,1,1,1)$.  Hence $[X_{S,l}] \in \Delta_2^{(2)}$.
\medskip

The cases 9, 11 are the mixed ones of case 5.  Hence $[X_{S,l}] \in
\Delta_2^{(3)}$.
\medskip

In the cases 10, 12, two singular fibres of type $IV$ and one fibre of type
$II$
degenerate to a singular fibre of type $II^*$.  The projections of some two
components of
the fibre of type $II^*$ are of type $(1,1,0,0,0), (1,-1,0,0,0)$.  Hence
$[X_{S,l}] \in \Delta_2^{(4)}$.
\medskip

The case 13 is the mixed one of cases 5 and 8.  Hence $[X_{S,l}] \in
\Delta_3^{(1)}$.
\medskip

The case 13* is the mixed one of cases 5 and 8*.  Hence $[X_{S,l}] \in
\Delta_3^{(2)}$.
\medskip

The cases 14, 15 are the mixed one of cases 5 and 10.  Hence $[X_{S,l}] \in
\Delta_3^{(3)}$.
\medskip

The case 16 is the mixed one of cases 5 and 13*.  Hence $[X_{S,l}] \in
\Delta_4^{(1)}$.
\medskip

The case 17 is the mixed one of cases 10 and 14.  Hence $[X_{S,l}] \in
\Delta_4^{(2)}$.

\subsection{Theorem}
{\it The isomorphism $\calK3_{M,\rho}^m\setminus \Delta^m\cong
\calM_{\cub}^m$ extends to a
$W(E_6)$-equivariant isomorphism
$$\calK3_{M,\rho}^m \cong \calM_{\ncub}^m.$$
Passing to the quotients it defines an isomorphism
$$\calK3_{M,\rho} \cong \calM_{\ncub}.$$}

\begin{proof} The isomorphism
$\calK3_{M,\rho}^m/W(D_5)\cong \calM_{\ncub}^m = \calM_{\ncub}^m/W(D_5)$
constructed in Theorem \ref{extensionthm} lifts to a
$W(E_6)$-equivariant isomorphism $\calK3_{M,\rho}^m \cong
\calM_{\ncub}^m$. In fact, this is true for open Zariski subsets defined
by nonsingular cubic surfaces, hence each of the varieties is the
normalization of the quotient in the  field of rational functions
$\bbC(\calK3_{M,\rho}^m)= \bbC(\calM_{\ncub}^m)$. Now we have an
isomorphism $\alpha$ of varieties which defines a birational isomorphism
of $W(E_6)$-varieties. Obviously, it is an isomorphism of
$W(E_6)$-varieties  (for each $g\in W(E_6)$ the maps $g\circ \alpha$ and
$\alpha\circ g$ coincide on an open Zariski subset, hence coincide
everywhere).
\end{proof}

\subsection{Corollary} {\it The isomorphism
$$(\calB\setminus \calH)/\Gamma_{M,\rho} \cong  \calM_{\cub}$$
from Corollary \ref{cor2.2} extends to an isomorphism
$$
\calB/\Gamma_{M,\rho} \cong \calM_{\ncub}.
$$}

\subsection{Remark}
The isomorphism $\calM_{\ncub}\cong \calK3_{M,\rho}$ can be extended to the
compactification obtained by adding one
strictly-semistable point. The image of this point goes to the cusp of the
ball quotient.
The corresponding $K3$ surface is isomorphic to the double cover of
$\bbP^2$  ramified along the sextic curve $$
t_2(L_1(t_0,t_1)^3 L_2(t_0, t_1)^2 + t_2^3 L_2(t_0,t_1)^2) = 0, $$
where $L_1, L_2$ are independent linear forms. This sextic appears as a
semistable sextic in Shah \cite{Sha},
Theorem 2.4, Group II, (2).
Its double cover is a Type II degeneration of $K3$ surfaces, i.e.
corresponding to a point on an 1-dimensional
rational boundary component of the period domain of polarized $K3$ surfaces
of degree 2 (= a bounded symmetric domain of
type IV and of dimension 19). The 1-dimensional rational boundary
components of a bounded symmetric domain of type IV
bijectively correspond to the set of totally isotropic primitive
sublattices of rank 2 of its underlying lattice of
signature $(2,r)$. In our situation, $\rho$-invariant totally isotropic
primitive sublattices of rank 2 of $T$
correspond to the set of cusps of $\calB$. Thus the semistable points of
type $(6,6)$ correspond to the boundary of the
Satake's compactification of $\calB/\Gamma_{M,\rho}$.

The strictly semistable
cubic surface defined by
\begin{equation}\label{sscubic}
X_3^3-X_0X_1X_2=0
\end{equation} (cf.\ \cite{ACT} (4.6))
has three double rational points of type $A_2$ and has only three lines
which lie in one $\Aut(S)$-orbit.
This defines three planes in the cubic fourfold $X$ defined by
$X_5^3+X_4^3+X_3^3-X_0X_1X_2=0$
(one such plane is $\Pi:\,X_2=X_3=X_4+X_5=0$) and projection away from such
a plane defines a quadric bundle
structure on $X$. The discriminant curve is easily computed and is a sextic
as above.

It follows from Proposition \ref{stable} that the pair
$(F_5,F_2) = (L_1^3L_2^2,L_2^2)$ represents a semi-stable but not stable
point in $\bbP(V(5))\times \bbP(V(2))$ whose orbit is closed in the set
of semi-stable points. The corresponding point in
$\bigl(\bbP(V(5))\times \bbP(V(2))\bigr)^{\textup{ss}}/\!/\SL(2)$
compactifies $\bigl(\bbP(V(5))\times
\bbP(V(2))\bigr)^{\textup{s}}/\SL(2)$. Thus we see that
$\calM_{\ncub}\li$ admits a one-point compactification corresponding to
the surface \eqref{sscubic} together with its unique (up to automorphism)
line.

\subsection{Configurations of 7 points in $\bbP^1$}
\label{configurations}

Recall from Theorem \ref{thm1} that we have  a natural isomorphism
$$\calM_{\ncub}\li \cong \bigl(\bbP(V(5))\times
\bbP(V(2))\big)^\textup{s}/\SL(2),$$
where $\bigl(\bbP(V(5))\times \bbP(V(2))\big)'$ is
the open subset corresponding to stable pairs of binary forms
$(F_5,F_2)$.  Consider the product $(\bbP_1)^7$ as the product
$(\bbP^1)^5\times (\bbP^1)^2$. We have an isomorphism
$$\psi:(\bbP^1)^7/S_5\times S_2 \to \bbP(V(5))\times \bbP(V(2)).$$
Let $p: (\bbP^1)^7\to \bbP(V(5))\times \bbP(V(2))$ be
the composition of the quotient map and $\psi$ and
$$\calL = p^*\bigl(\calO_{\bbP(V(5))}(2)\boxtimes
\calO_{\bbP(V(2))}(1)\bigr) \cong
\boxtimes_{\substack{i=1}}^{\substack{5}}\calO_{\bbP^1}(2)\otimes
(\calO_{\bbP^1}(1)\boxtimes \calO_{\bbP^1}(1)).$$
Since the stability is preserved under the action of finite groups,
we see that semi-stable (stable) points in $\bbP(V(5))\times \bbP(V(2))$
with respect to the action of $\SL(2)$ and the linearization defined by
the invertible sheaf $\calO_{\bbP(V(5))}(2)\boxtimes
\calO_{\bbP(V(1))}(1)$ correspond to semi-stable (stable) points in
$(\bbP_1)^7$ with respect to the diagonal action of $\SL(2)$ and the
linearization defined by the line bundle $\calL$. Let
$$P_1(2^5,1,1) = \bigl((\bbP_1)^7\bigr)^s/\SL(2).$$
We have
$$\bigl(\bbP(V(5))\times \bbP(V(2))\big)^\textup{s}/\SL(2)
\cong P_1(2^5,1,1)/S_5\times S_2.$$
We know that
$\calM_{\ncub}^{\textup{l}}= \calM_{\ncub}^m/W(D_5)$. The group $W(D_5)$ is
equal to
the semi-direct product $(\bbZ/2\bbZ)^4\ltimes S_5$. Here
$S_5$ is the subgroup of $W(D_5)$ which acts on markings on nonsingular
surfaces by permuting  the divisor classes $e_1,\ldots,e_5. $ It
stabilizes the divisor class $2e_0-e_1-\ldots-e_5$ of a line $l$. The
subgroup $H = (\bbZ/2\bbZ)^4$ is generated by the conjugates of the
product of two commuting reflections
$s_{e_0-e_1-e_2-e_6}\circ s_{e_1-e_2}$. Let  $l_i'$ be the lines
representing the classes $e_0-e_i-e_6$. Then $H$ acts by switching even
number of $l_i$'s with $l_i'$'s. The proof of Theorem \ref{thm1} shows
that the map $\calM_{\ncub}\li \to \bigl(\bbP(V_5)\times
\bbP(V_2)\big)'/\SL(2)$ induces a $S_5$-equivariant isomorphism
$$\calM_{\ncub}^m/H \cong P_1(2^5,1,1)/S_2.$$

\subsection{Monodromy groups.}
According to Deligne and Mostow \cite{DM}, the variety $P_1(2^5,1,1)$ is
isomorphic to the quotient of a
complex 4 ball by a reflection subgroup $\Pi'$ corresponding to
hypergeometric function defined by the multi-valued
form
$$
\omega= z^{-1/6}[(z-1)(z-a_1)(z-a_2)(z-a_3)(z-a_4)]^{-1/3}dz. $$ They also
show that
$\Pi'$ and $S_2$ generate a reflection subgroup $\Pi$ such that the ball
quotient is isomorphic to $P_1(2^5,1,1)/S_2$.
As shown in \ref{visot}, $X$ is the minimal model of a quotient $(C\times
E)/(\bbZ/6\bbZ).$
This correspondence gives us an isogeny between our group $\Gamma_{\rho}$
and $\Pi$.

\section{Half twists.}
\label{halftwist}
\subsection{}
To a smooth cubic surface $S$ one can associate a principally polarized
Hodge structure of rank 10 and weight 1,
it is $H^1(P,\bbZ)$ where $P$ is the intermediate Jacobian of the cubic
threefold $V$ (cf.\ \ref{3fold}) associated to
$S$. In \cite{ACT}, see also \cite{MT}, it is shown that this Hodge
structure, with its automorphism of order three,
determines $S$.

The automorphism of order three defines the structure of a free
$\bbZ[\zeta]$-module on $H^1(P,\bbZ)$. It defines
eigenspaces $H^{1,0}(P)_\chi$ and $H^{1,0}(P)_{\bar\chi}$ of dimension $4$
and $1$ respectively. This allows one to
define a weight two Hodge structure $W$, with Hodge numbers $(1,8,1)$, and
with the same underlying lattice
$W=H^1(P,\bbZ)$ as follows:
$$
W^{2,0}=H^{1,0}(P)_{\bar\chi},\quad
W^{1,1}=H^{1,0}(P)_{\chi}\oplus H^{0,1}(P)_{\bar\chi},\quad
W^{0,2}=H^{0,1}(P)_{\chi},
$$
in fact it is easy to check that $W^{p,q}=\overline{W^{q,p}}$. The
automorphism of order three of $H^1(P,\bbZ)$
preserves this decomposition, hence also $W$ has an automorphism of order
three. The polarization $E$ on
$H^1(P,\bbZ)$ defines a $\bbQ[\zeta_3]$-valued Hermitian form $H$ on
$H^1(P,\bbZ)\cong \bbZ[\zeta_3]^5$
(cf.\ \cite{ACT}) with imaginary part $E$. The real part $Q$ of $H$ is a
polarization of $W$. The lattice $(W,Q)$
is of type $A_2^4\oplus A_2(-1)$.
The polarized Hodge structure $(W,Q)$ is the (negative) half twist of
$(H^1(P,\bbZ),E)$ (\cite{vG1}).

\subsection{} The lattice $(W,Q)\cong A_2^{\oplus 4}\oplus A_2(-1)$ has a
unique (up to an isometry) embedding in
the $K3$ lattice $L$ and the automorphism of order three on $ W$ extends to
an automorphism of order three on the
$K3$ lattice. The polarized Hodge structure $(W,Q)$ is invariant under this
automorphism and defines a $K3$ surface
with an automorphism of order three. So the half twist of $H^1(P,\bbZ)$
provides a purely Hodge theoretic approach to
the $K3$ surfaces which were constructed as triple covers of cubic surfaces
in this paper.

\end{document}